\documentclass[11pt,a4paper]{amsart}
\usepackage{latexsym}
\usepackage{amssymb}
\usepackage{bbm,amssymb,amsthm,amsmath}
\usepackage{graphicx}
\usepackage{multirow}
\usepackage{caption}
\usepackage{mathtools} %za duge strelice
\usepackage{tikz}	% tikz to creat diagrams
\usetikzlibrary{shapes,arrows}
\usetikzlibrary{positioning}

\tikzstyle{type1} = [rectangle, rounded corners, draw=none, 
text width=3cm, text centered, inner sep=3pt, minimum height=2em]
\tikzstyle{type2} = [rectangle, rounded corners, draw=none, 
text width=2.5cm, text centered, inner sep=3pt, minimum height=2em]
\usepackage{xcolor}

\calclayout

\makeatletter	
\@namedef{subjclassname@2020}{%
	\textup{2020} Mathematics Subject Classification}
\makeatother
\theoremstyle{plain}
\newtheorem{theorem}{Theorem}[section]	% brojevi tvrdnji imaju i broj odjeljka
\newtheorem{lemma}[theorem]{Lemma}
\newtheorem{corollary}[theorem]{Corollary}

\theoremstyle{definition}
\newtheorem{definition}[theorem]{Definition}
\newtheorem{example}[theorem]{Example}

\theoremstyle{remark}
\newtheorem{remark}[theorem]{Remark}

\numberwithin{equation}{section}	% da brojevi formula imaju i broj odjeljka

\let\mib=\boldsymbol

\def\C{{\mathbb{C}}}
\def\R{\mathbb{R}}    
  
\def\N{\mathbb{N}}   
\def\dup#1#2#3#4{{}_{#1\!}\langle\, #2 , #3 \,\rangle_{#4}}
\def\dupp#1#2#3#4{{\vphantom{\big\langle}}_{#1\!}\bigl\langle\, #2 , #3 \,\bigr\rangle_{#4}}
\def\lR{\mathcal{R}}
\def\lS{\mathcal{S}}

\def\lJ{\mathcal{J}}
\def\lA{\mathcal{A}}
\def\lL{\mathcal{L}}
\def\lM{\mathcal{M}}
\def\lD{\mathcal{D}}
\def\malpha{{\mib \alpha}}
\def\mbeta{{\mib \beta}}
\def\mgamma{{\mib \gamma}}

\def\mnu{{\mib \nu}}
\def\mzeta{{\mib \zeta}}
\def\mpi{{\mib \pi}}

\def\mx{{\bf x}}
\def\mxi{{\mib \xi}}
\def\meta{{\mib \eta}}
\def\mnul{{\bf 0}}

\def\my{{\bf y}}
\def\mz{{\bf z}}

\def\mB{{\bf B}}
\def\mA{{\bf A}}
\def\Rd{{\mathbb{R}^{d}}}   
\def\Rr{{\mathbb{R}^{r}}}

\def\Sdmj{{\rm S}^{d-1}}
\def\supp{{\rm supp\,}}

\def\Nnul{{\mathbb{N}_0}}
\def\Nnb{{\Nnul\cup\{\infty\}}}

\def\kobb{{\rm K_{0,\infty}}}
\def\kob{{{\rm K_{0,\infty}}(\Rd)}}

\def\kb{{{\rm K_{\infty}}}(\Rd)}
\def\Rdz{\mathbb{R}^d_*}
\def\Cr{\mathbb{C}^r}

\def\mmp{{\bf p}}
\def\Dupp#1#2{\left\langle#1,#2\right\rangle}
\def\Mx#1#2{{\rm M_{#1}}(#2)}

\def\Kz#1{\mathrm{B}[\vnul,#1]}
\def\Ko#1{\mathrm{B}(0,#1)}

\def\Cp#1{{{\rm C}(#1)}} 
\def\pC#1#2{{{\rm C}^{#1}(#2)}}

\def\W#1#2#3{{{\rm W}^{#1,#2}(#3)}}  
\def\Wloc#1#2#3{{{\rm W}^{#1,#2}_{{\rm loc}}(#3)}} 

\def\jhds#1{{#1}_\kobb}
\def\jhd{\nu_\kobb}
\def\jhdv{\mnu_\kobb}
\def\hd{\nu_H}
\def\pkd{\nu_{sc}}

\def\rnul{r_0} %u definiciji preslikavanja \lJ (\rnul={\sqrt{3}\over 3})
\def\vnul{\mathsf{0}}
\def\ve{\mathsf{e}}
\def\kvz{\mathrm{S}^d_{[0,r_1]}}
\def\kvo{\mathrm{S}^d_{(0,r_1)}}

\def\kvn{\mathrm{S}^d_{r_1}}
\def\kvb{\mathrm{S}^d_{0}}

\def\Int{\operatorname{Int}}
\def\cl{\operatorname{Cl}}
\def\rest#1{_{\displaystyle |_{#1}}}

\def\vu{\mathsf{u}}
\def\vv{\mathsf{v}}
\def\vg{\mathsf{g}}
\def\vf{\mathsf{f}}
\def\ph{\varphi}
\def\eps{\varepsilon}

\begin{document}
%\today

\title{One-scale H-distributions and variants}

\author{N.~Antoni\'c}\address{Nenad Antoni\'c,
	Department of Mathematics, Faculty of Science, University of Zagreb, Bijeni\v{c}ka cesta 30,
	10000 Zagreb, Croatia}\email{nenad@math.hr}

\author{M.~Erceg}\address{Marko Erceg,
	Department of Mathematics, Faculty of Science, University of Zagreb, Bijeni\v{c}ka cesta 30,
	10000 Zagreb, Croatia}\email{maerceg@math.hr}



\begin{abstract}
H-measures and semiclassical (Wigner) measures were introduced in early
1990s and since then they have found numerous applications in problems involving 
$\mathrm{L}^2$ weakly converging sequences. 
Although they are similar objects, neither of them is a generalisation of the other,
the fundamental difference between them being the fact that 
semiclassical measures have a characteristic length, while 
H-measures have none. 
Recently introduced objects, the one-scale H-measures, generalise both of them, thus
encompassing properties of both.

The main aim of this paper is to fully develop this theory to the 
$\mathrm{L}^p$ setting, $p\in(1,\infty)$, by constructing {\sl one-scale H-distributions},
a generalisation of one-scale H-measures and, at the same time, of H-distributions, a generalisation of
H-measures to the $\mathrm{L}^p$ setting. 
We also address an alternative approach to $\mathrm{L}^p$ extension of semiclassical measures
via the Wigner transform, introducing new type of objects (semiclassical distributions).

Furthermore, we derive a localisation principle in a rather general form, suitable for problems with a characteristic length,
as well as those not involving a specific characteristic length, providing some applications. 
 
\end{abstract}

\maketitle

%\clearpage

\tableofcontents

\section{Introduction}

\subsection{Overview of microlocal defect objects}

\emph{H-measures} or \emph{microlocal defect measures} 
were introduced independently by Luc Tartar \cite{Tar}
and Patrick G\'erard \cite{GerMdm} in early 1990s. 
They represent a generalisation of defect measures \cite{Lio1} in the sense that, 
besides the space variables, they depend on the dual variables as well. 
More precisely, an H-measure is a Radon measure on the cospherical bundle $\Omega\times\Sdmj\subseteq T^\ast\Rd\simeq \Rd\times\Rd$
over a domain $\Omega\subseteq\Rd$, and it is associated to a weakly converging 
sequence in $\mathrm{L}^2_\mathrm{loc}(\Omega)$.
The main advantage of H-measures over defect measures (and Young measures) can be seen when studying 
oscillating sequences (the oscillation effects, together with concentration effects, are the typical 
cause for non-compactness). 
Indeed, the information about direction of oscillation 
cannot be captured either by the defect measures or the Young measures, while the H-measures 
proved to be quite successful for that task. 

Let us illustrate this on a simple example of a plain wave:
\begin{equation}\label{eq:plain_wave}
u_n(\mx) = \varphi(\mx) e^{2\pi i \frac{\mx}{\eps_n}\cdot\mathsf{k}} \,,
\end{equation}
where $\ph\in\mathrm{L}^2_\mathrm{loc}(\Rd)$, $\mathsf{k}\in\Rd\setminus\{0\}$, 
and $\eps_n\to 0^+$.
This sequence weakly converges in $\mathrm{L}^2_\mathrm{loc}(\Rd)$ to $0$
(but not strongly, except in the trivial case $\ph=0$).
The associated defect measure is the limit of a stationary sequence $|u_n|^2=|\ph|^2$ in
the space of (unbounded) Radon measures with respect to the weak-$\ast$ 
topology. Thus, it is equal to $|\ph|^2\lambda$, i.e.~the measure having 
density $|\ph|^2$ with respect to the Lebesgue measure $\lambda$ on $\Rd$. 
Here we can see that the defect measure is independent of the direction of
oscillations $\mathsf{k}$.
On the other hand, it can be shown \cite{AEL,Tar} that the H-measure defined by this sequence is 
\begin{equation}\label{eq:plain_wave_Hm}
|\ph|^2\lambda\otimes\delta_{\frac{\mathsf{k}}{|\mathsf{k}|}} \,,
\end{equation}
where
$|\ph|^2\lambda$ and $\delta_{\frac{\mathsf{k}}{|\mathsf{k}|}}$ (the Dirac measure
at point $\mathsf{k}/|\mathsf{k}|$)
are
measures in the physical space (variable $\mx$) and the dual space (variable $\mxi$),
respectively. 
Hence, the direction of oscillation is inherent in the H-measure.

In the theory of partial differential equations one often deals with weakly 
converging sequences (see e.g.~\cite{Evans}), allowing H-measures and their 
(anisotropic) variants
to find numerous applications in the field. 
Let us mention the investigation on
propagation of oscillations and on concentration effects in
solutions of partial differential equations \cite{Asym, ALjfa, EI18, IvecLazar},
the applications in the control theory \cite{Burq,DLR, LZ}, 
the velocity averaging results \cite{EMM20,GerMdm,LM2}, the existence of traces and solutions 
to nonlinear evolution equations \cite{AM_jhde, EM20,HKMP,pan_jhde,pan_arma},
as well as explicit formulae and bounds in homogenisation \cite{ABJ,ALjmaa,Lazar,Tar}.
Most of these applications depend on the so-called localisation principle, which 
is closely related to the generalisation of compactness by compensation method 
to variable coefficients \cite{GerMdm,Tar}.

\begin{figure}[h!]
	\centering
	\begin{tikzpicture}[auto, scale=0.9, every node/.style={transform shape}]
	\node at (0,0) [type1] (jHd) {\textbf{one-scale\\ H-distributions}\\
		(Theorem \ref{thm:jhd_exist})};
	\node[type1, below left=0.7cm and -0.2cm of jHd] (Hd) {H-distributions \cite{AM}};
	\node[type1, below right=0.7cm and -0.2cm of jHd] (scd) {\textbf{semiclassical\\
			distributions}\\ (Definition  \ref{def:semiclass-distr})};
	\node[type2, left=1cm of Hd] (mcf) {microlocal\\ compactness forms \cite{rindler}};
	\node[type2, below=0.7cm of mcf] (Ym) {Young\\ measures \cite{Young}};
	
	\node[type1, below=4.5cm of jHd] (jHm) {one-scale\\ H-measures\\
		\cite{AEL, tar_book}};
	\node[type1, below right=0.7cm and -0.2cm of jHm] (scm) {semiclassical 
		measures (aka Wigner measures)\\ \cite{GerMsc, LP93}};
	\node[type1, below=1.8cm of jHm] (mHm) {multi-scale H-measures \cite{Tar2}};
	\node[type1, below left=0.7cm and -0.2cm of jHm] (Hm) {H-measures\\
		(aka microlocal defect measures)\\ \cite{GerMdm,Tar}};
	\node[type2, left=1cm of Hm] (dm) {defect\\ measures \cite{Lio1}};
	\node[type1, below=0.7cm of Hm] (pHm) {parabolic, ultra-parabolic, fractional
		and adaptive H-measures\\ \cite{ALjfa,EI17, EMM20,LM2,MI,pan_aihp}};
	
	\path
	(jHd) edge [thick,-stealth]  (Hd)
	(jHd) edge [thick,-stealth] (scd)
	(jHd) edge [thick,-stealth] (jHm)
	(Hd) edge [-stealth] (Hm)
	(Hd) edge [dashed,-] node[above] {\cite{AMP}} (mcf)
	(mcf) edge [-stealth] (Ym)
	(jHm) edge [-stealth] (Hm)
	(jHm) edge [-stealth] (scm)
	(Hm) edge [-stealth] (dm)
	(Ym) edge [-stealth] (dm)
	(mHm) edge [-stealth] (jHm)
	(pHm) edge [-stealth] (Hm)
	(scd) edge [thick,-stealth] (scm);
	
	\draw [color=blue,thick,dash dot, rounded corners](-1.6,0.7) rectangle (4.8,-9.9);
	\node[label={\textcolor{blue}{Variants with a characteristic length}}] at (2,0.7) {};
	\draw [color=red,thick, rounded corners](-4.8,-5) rectangle (5.2,-13.4);
	\node[label=below:{\textcolor{red}{Variants on ${\rm L}^2$ space}}] at (3,-13.4) {};
	\draw [color=gray,thick, dotted, rounded corners](-5.2,-10) rectangle (-1,-13.8);
	\node[type2] at (-6.5,-12) {\textcolor{gray}{Variants with 
			anisotropic scaling in variables}};
	\end{tikzpicture}
	\caption{Diagram of (microlocal) defect tools of interest. In bold we put
		new objects that are introduced in this paper. Pointing arrows suggest
		that the object in the tail is a generalisation of the object in the head
		of the arrow, while a dashed line means that objects are in a certain 
		relation. For each variant at least one reference is presented 
		(e.g.~the first source on the topic or/and later references with a 
		more detailed exposition).}
	\label{fig:diagram}
\end{figure}
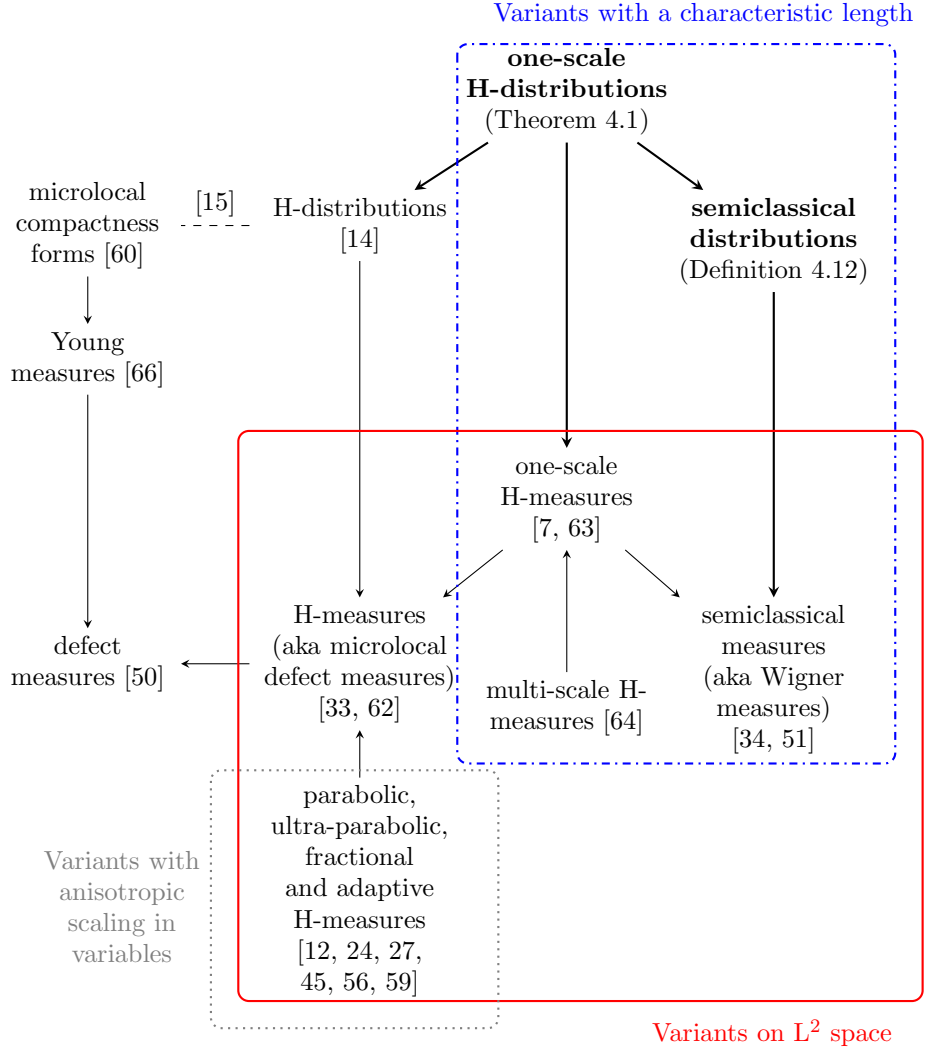

In the example above we can see that the H-measure \eqref{eq:plain_wave_Hm}
does not distinguish between sequences \eqref{eq:plain_wave} with different frequencies
$\frac{1}{\eps_n}$. In order to have such sensitivity we need to incorporate 
a scale into the definition of microlocal defect objects.
That is the case with \emph{semiclassical measures};
the Radon measures on the cotangential bundle $\Omega\times\Rd$
introduced by G\'erard \cite{GerMsc}, and later baptised \emph{Wigner measures}\/ by Pierre-Louis Lions
and Thierry Paul \cite{LP93}.
Since they depend upon a characteristic length $(\omega_n)$, $\omega_n\to 0^+$ in the real line,
they are more suitable in situations where such a characteristic length 
naturally appears, often related to highly oscillating problems for 
partial differential equations. Let us just mention the homogenisation limit \cite{GMMP, LP93},
the microlocal energy density for (semi)linear wave equation \cite{FrInt,Goce},
the semiclassical limit of Schr\"odinger equations (propagation of singularities
and derivation of effective equations)
\cite{CFKMS, CFKM20, FKGL},
the mean field limit \cite{AN,LNR14,MNO19} (for the relation with the bosonic quantum 
de Finetti theorem see \cite[Section 3]{AFP16}), 
as well as some other problems related to quantum mechanics.

However, the incorporation of a scale into the microlocal defect object as a byproduct brings 
a new issue. Namely, if the characteristic length $(\omega_n)$ 
of a semiclassical measure is chosen inappropriately, we can lose 
a significant amount of information. For example, if $\lim_n\frac{\omega_n}{\eps_n}=+\infty$,
the semiclassical measure associated to \eqref{eq:plain_wave} is equal to 0
(zero measure). There are examples where this drawback is even more severe
(cf.~\cite[Example 4]{EL18}). 
This in particular implies that, in contrast to H-measures, 
a zero semiclassical measure does not necessarily guarantee the strong convergence of the 
associated sequence (the so-called $(\omega_n)$-oscillatory property needs to be satisfied as well
\cite{GMMP,EL18}).

The above discussion shows that H-measures and semiclassical measures are in a general relation (neither is a 
generalisation of the other) and either has some advantages and disadvantages.
Recently, Tartar introduced \emph{one-scale H-measures} \cite[Chapter 32]{tar_book} 
(see also \cite{AEL,Tar2}) as a true extension of H-measures and semiclassical 
measures. Moreover, the localisation principles, both for H-measures and semiclassical 
measures, can be derived from the localisation principle for one-scale H-measures
\cite{AEL}.
The main part in their construction is a proper choice of the domain for 
dual variables. 
For a fine tuning with characteristic lengths the set has to be \emph{thick enough}, 
but we need to allow for directions to be \emph{detected} also at the origin and infinity.
This can be achieved with a radial compactification of $\Rd\setminus\{0\}$, denoted by
$\kob$, which is homeomorphic to the $d$-dimensional spherical shell (cf.~\cite{AEL}).  
Hence, one-scale H-measures are defined on $\Omega\times\kob$.

Heretofore mentioned microlocal defect objects are suitable only for problems expressed 
in the $\mathrm{L}^2$ framework.
In order to overcome this limitation, \emph{H-distributions} were introduced \cite{AM}, 
which are an extension of H-measures to the $\mathrm{L}^p$ context
($1<p<\infty$).
Due to absence of positivity (which was trivially present in the $\mathrm{L}^2$ setting),
H-distributions are in general not measures, but Schwartz' distributions 
\cite{Schwartz}. 
However, in \cite{AEM} a more precise
description is given through anisotropic distributions. 
Another approach for studying $\mathrm{L}^p$ weakly converging sequences was based 
on the extension of generalised Young measures obtaining \emph{microlocal compactness forms}
\cite{rindler},
which turned out to be closely related to H-distributions \cite{AMP}.
The H-distributions and their variants found their application in developing velocity averaging 
results \cite{LM3} and the $\mathrm{L}^p$ version of compactness by compensation 
\cite{AM,MM,rindler}.

In this paper our intention is to close the picture in the $\mathrm{L}^p$
setting, as it was already done for $p=2$ in \cite{AEL}. More precisely, we shall develop the \emph{one-scale
H-distributions} as the $\mathrm{L}^p$ counterpart to one-scale H-measures,
which, at the same time, generalise H-distributions. Moreover, we can recognise 
the $\mathrm{L}^p$ counterpart of semiclassical measures, the objects which can 
still be achieved by means of the Wigner transform (see Figure \ref{fig:diagram}).
Besides the purely academic motivation to uplift the existing theory to the
$\mathrm{L}^p$ setting, we hope that these new objects might find 
more practical applications, like in refining and generalising
the procedure of second microlocalisation \cite{CFK00},
analysing (semiclassical) wave front sets 
(cf.~\cite[chapters 5 and 8]{Zworski}, \cite{Zanelli}),
and studying homogenisation (as well as semiclassical) limits of 
partial differential equations, all in ${\rm L}^p$ spaces.
For the latter case we emphasise the situations where the nonlinearity 
in partial differential equations forces one to work in the ${\rm L}^p$ framework.

%\smallskip

\subsection{Notation}

%{\color{blue}

Before we proceed, let us briefly recall the notation which we shall use in the paper.
For denoting derivatives of (complex) functions defined on $\Rd$, we use multiindices $\malpha\in\N_0^d$, their length
(or 1-norm) being $|\malpha|=\sum \alpha_i$. However, for vectors or points in $\Rd$ we use $|\mx|$ to denote the
Euclidean norm (the 2-norm), and similarly on $\C^d$. In other cases we use indices to make precise which (semi)norm we are using.
We shall use the projection $\mpi(\mxi)=\mxi/|\mxi|$ defined on $\Rdz:=\Rd\setminus\{\vnul\}$, and taking values on
the unit sphere $\Sdmj\subseteq\Rd$.
By $\lfloor a\rfloor$ we denote the largest integer not greater than $a\in\R$.
For a set $A$, $\cl A$ denotes its closure, while $\Int A$ its interior (the respective topology will always be clear from the contest).
In particular, for a subset of complex numbers, this allows $\overline A$ to be used for the set of its complex conjugates.

%	\Mx{k\times l}\Cr  
%	\otimes
%	\lL(B) om. operatori

By $\hat{u}={\mathcal F}u$ we shall denote the Fourier transform of $u$, where for $u\in{\rm L}^1(\Rd)$ we use the formula:
$\hat{u}(\mxi)=\int_\Rd e^{-2\pi i\mx\cdot\mxi}u(\mx)d\mx$.
By $\check{u}=\bar{{\mathcal F}}u$ we denote it's inverse, while by ${\mathcal A}_\psi$ the Fourier multiplier operator with symbol
$\psi$: ${\mathcal A}_\psi(u) = (\psi \hat{u})^{\lor}$.

For $p\in[ 1,\infty]$, by ${\rm L}^p_{{\rm loc}}(\Rd)$ we denote the space of all distributions $u$ such that the following holds
$$
	(\forall\varphi\in{\rm C}^\infty_c\Rd)\quad \varphi u\in{\rm L}^p(\Rd)\;.
$$
Actually, ${\rm C}^\infty_c(\Rd)$ can be reduced to ${\mathcal G}$, its subset such that $(\forall \mx\in\Rd)(\exists\varphi\in{\mathcal G})\; {\rm Re}\;\varphi(\mx) > 0$,
which can be chosen to be countable. We endow it with the locally convex topology induced by a family of seminorms $|\cdot|_{\varphi,p}$ (for $\ph\in\mathcal{G}$)
$$
	|u|_{\varphi,p} := \|\varphi u\|_{{\rm L}^p(\Rd)}\;.
$$
It can be shown that neither the definition of ${\rm L}^p_{{\rm loc}}(\Rd)$ nor its topology depend on the choice of family ${\mathcal G}$
with the above property. 
This definition is equivalent to a definition where one requires that ${\rm L}^p_{{\rm loc}}$ functions have finite ${\rm L}^p$ norms
over every compact subset of $\Rd$ (indeed, we can take ${\mathcal G}$ to consist of all characteristic
functions $\chi_K$ of compacts $K\subseteq\Rd$ and notice that smoothness is not actually needed for the definition of ${\rm L}^p_{\rm loc}(\Rd)$ space).

We say that a sequence $(u_n)$ is bounded in ${\rm L}^p_{{\rm loc}}(\Rd)$ if for every seminorm $|\cdot|_{\varphi,p}$ there
exists $C_{\varphi,p} >0$ such that $|u_n|_{\varphi,p} < C_{\varphi,p}$ uniformly in $n$.
By choosing a countable ${\mathcal G} = \{\vartheta_l\; :\; l\in\N \}$ such that $0\le \vartheta_l\le 1$ and $\chi_{K_l}\le \vartheta_l\le\chi_{K_{l+1}}$,
where $K_l:=\Kz l\subseteq \Rd$ is a closed ball of radius $l$ centred around the origin, we can define a metric $d_p$ on ${\rm L}^p_{{\rm loc}}(\Rd)$ by
$$
	d_p(u,v) := \sup_{l\in\N} 2^{-l} { |u-v|_{\vartheta_l,p} \over 1 + |u-v|_{\vartheta_l,p}}\;.
$$
With this metric, ${\rm L}^p_{{\rm loc}}(\Rd)$ is a Fr\'{e}chet space for each $p\in[ 1,\infty]$, separable for $p\in[ 1,\infty)$
and reflexive for $p\in(1,\infty)$.
For $p\in[ 1,\infty)$ it is also valid that ${\rm C}^\infty_c(\Rd)$ is dense in ${\rm L}^p_{{\rm loc}}(\Rd)$, while ${\rm L}^{p'}_c(\Rd)$,
a subspace of ${\rm L}^{p'}(\Rd)$ consisting of all functions in that space having a compact support, equipped with the topology
of strict inductive limit
$$
{\rm L}^{p'}_c(\Rd) = \bigcup_{l\in\N} {\rm L}^{p'}_{K_l}(\Rd)\;,
$$
is the dual of ${\rm L}^p_{\rm loc}(\Rd)$.
Of course,  we define
$$
{\rm L}^{p'}_{K_l}(\Rd) := \Bigl\{f\in{\rm L}^{p'}(\Rd): \supp f\subseteq K_l\Bigr\} \sim {\rm L}^{p'}(K_l)\;,
$$
and equip it with the ${\rm L}^{p'}$ norm topology.
Let us just remark that we could have replaced $\Rd$ by any open set $\Omega\subseteq\Rd$ and all of the above definitions and conclusions would remain valid.
For omitted proofs and further references, we refer the interested reader to \cite{ABu}.

Our last point regarding the notation is that by $\langle\cdot,\cdot\rangle$ we shall denote various duality products, always assuming that it is bilinear.
In particular, this means that its relation to the ${\rm L}^2$ scalar product (which we take to be antilinear in the first variable) is as follows
$$
	\langle f|g\rangle_{{\rm L}^2} = \int \overline{f(\mx)} g(\mx) d\mx = \langle \bar f, g\rangle \;.
$$

%}

\subsection{Organisation of paper}

%\smallskip
%\noindent{\bf Organisation of paper. } 

The second section we start by defining the radial 
compactification of $\Rdz$ and then we study 
smooth functions on it. 
These results, 
together with the kernel theorem for anisotropic distributions on 
manifolds with boundary presented in the third section, 
are crucial in the construction of one-scale H-distributions,
which is conducted in the fourth section.
In the rest of the fourth section we study some properties of one-scale H-distributions, and then present two important 
examples of weakly (but not strongly) converging sequences: the concentration 
and the oscillation. The section is closed by an alternative 
construction of semiclassical distributions via the Wigner transform. 
The last section is devoted to the derivation of the localisation principle 
suitable for studying differential equations with a characteristic length.
The paper is closed by an example illustrating possible applications
of the localisation principle. 

We would like to mention that a preliminary version of some results
presented in this paper appeared as a part of Marko Erceg’s Ph.D.~thesis \cite{MEphd}.

\smallskip

%\clearpage

\section{Space of test functions in the dual space}

\subsection{Compactification}

In the construction of any microlocal defect functional it is fundamental to carefully specify spaces of test functions,
both in the physical and the Fourier (moment or dual) space. In this section we focus on the latter.

For the construction of one-scale H-measures (called \emph{variant of H-measures using a characteristic length}\/ there), 
Tartar \cite[Definition 32.5]{Tar2} chose the space of continuous 
functions on the suitable 
compactification of $\Rdz:=\Rd\setminus\{\vnul\}$. More precisely, 
let $\kob$ be the set obtained by 
adding two unit spheres, namely $\Sigma_0$ and $\Sigma_\infty$, to $\Rdz$, the first one at the origin and 
the second one at infinity, equipped with the topology making $\kob$ homeomorphic to a
$d$-dimensional spherical shell, as depicted in Figure \ref{fig:top_compt}. 
The points on these two spheres we 
denote by $0^\ve$ and $\infty^\ve$, respectively, where $\ve\in \Sdmj$ is the (radial) direction 
of the point (i.e.~$0^\ve$ and $\infty^\ve$ are the endpoints of a ray passing through $\ve$). 
In the plane, ${\rm K_{0,\infty}}(\R^2)$ is homeomorphic to the closed annulus.

\begin{figure}[h!]
	\centering
	\includegraphics{slika1.mps}
	\caption{Visualisation of $\kob$, the compactification of $\Rdz$.}
	\label{fig:top_compt}
\end{figure} 

Since for the construction of one-scale H-distributions we 
shall need the differential structure on 
$\kob$ (in order to apply the H\"ormander-Mihlin theorem; see Theorem \ref{thm:mihlin}), 
we present the compactification above in more detail.

For the compactifying map $\lJ$ we take the composition of 
the translation from the origin in the radial direction for $\rnul>0$:
\begin{equation}\label{eq:dfT}
\Rdz\ni\mxi\overset{\mathcal{T}}{\longmapsto} \frac{|\mxi|+\rnul}{|\mxi|}\mxi\in\Rd\setminus\Kz{\rnul} \;,	
\end{equation}
where $\Kz{\rnul}$ denotes the closed ball of radius $\rnul$ centred at the origin,
and a compactifying map of the radial compactification.
Concerning the radial compactification, we shall consider the standard construction obtained by 
the modified stereographic projection (cf.~\cite[Chapter 1]{RM_lect}), 
which we briefly sketch here. 

\begin{figure}[!ht]
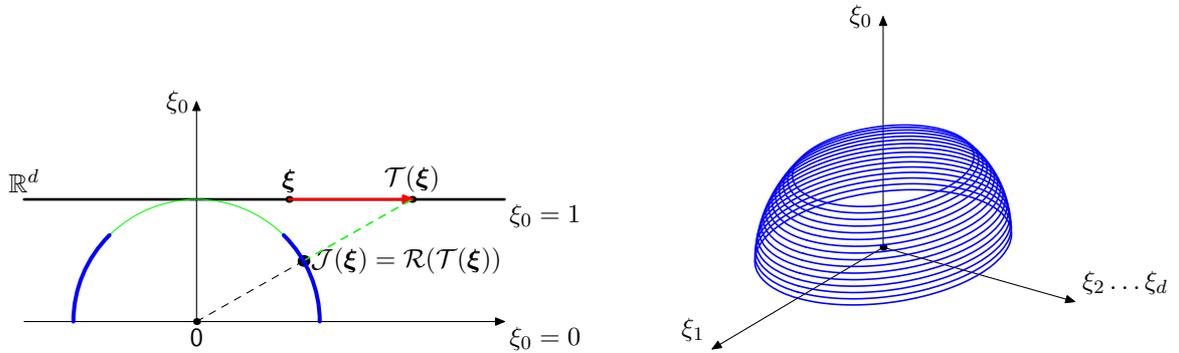

	\centering
	\includegraphics[scale=0.9]{slika2.mps} \qquad\quad
	\includegraphics[scale=0.9]{slika3.mps}
	\caption{On the left there is a graphical description of the mapping
		$\lJ$ for $\rnul=1$, where $\kvo$ is denoted in blue, while a
		visualisation of $\kvo$ in space is given on the right.}
	\label{fig:dif_compt}
\end{figure}

We first identify $\Rd$ with the hypersurface $\xi_0=1$ in $\R^{1+d}_{\xi_0,\mxi}$, 
and then apply the 
modified stereographic projection based on the line through the origin (instead of the
South Pole). More precisely, the radial compactification map $\mathcal{R}$ maps $\mxi$ to the 
intersection of $[0,1]\ni t\mapsto(t,t\mxi)$ (the line through $(1,\mxi)$ and $(0,\vnul)$ 
in $\R^{1+d}$) and the upper half of the unit sphere centred at the origin:
$\mathrm{S}^d_+:=\{(\zeta_0,\mzeta)\in\mathrm{S}^d : \zeta_0>0\}$
(see Figure \ref{fig:dif_compt} for a visual description). 
Since the intersection occurs at $t=(1+|\mxi|^2)^{-\frac{1}{2}}$, we have
that $\mathcal{R}:\Rd\to\mathrm{S}^d_+$
is given by
$$
\mathcal{R}(\mxi) = \Bigl(\frac{1}{\sqrt{1+|\mxi|^2}},\frac{\mxi}{\sqrt{1+|\mxi|^2}}\Bigr) \;.
$$ 

Finally, let us define
$\mathcal{J}:=\mathcal{R}\circ\mathcal{T}$, where $\mathcal{T}$ was defined in \eqref{eq:dfT}.
Since 
$$
\mathcal{R}(\Rd\setminus\Kz{\rnul})=\Bigl\{(\zeta_0,\mzeta)\in\mathrm{S}^d : 
	0<\zeta_0< r_1 \Bigr\} =: \kvo \,,
$$
where $r_1:=(1+\rnul^2)^{-1/2}$, we have that
$\lJ:\Rdz\to\kvo$ and has the explicit expression
\begin{equation}\label{eq:lj}
\lJ(\mxi) =  \biggl(\frac{1}{\sqrt{1+(|\mxi|+\rnul)^2}},\frac{|\mxi|+\rnul}
	{\sqrt{1+(|\mxi|+\rnul)^2}}\frac{\mxi}{|\mxi|}\biggr) \,.
\end{equation}
Here on $\kvo$ we consider the topology inherited from $\mathbb{R}^{1+d}$,
thus 
\begin{equation}\label{eq:cl-kvz}
\cl\kvo = \Bigl\{(\zeta_0,\mzeta)\in\mathrm{S}^d : 
	0\leq \zeta_0\leq r_1 \Bigr\} =: \kvz
\end{equation}
is a compact set. 
Furthermore, $\lJ$ is a $\mathrm{C}^\infty$-diffeomorphism (thus a homeomorphism as well), 
and its inverse $\lJ^{-1}:\kvo\to\Rdz$ is given by
\begin{equation*}
\lJ^{-1}(\zeta_0,\mzeta) =  \frac{\mzeta}{\zeta_0} - \rnul\frac{\mzeta}{|\mzeta|} \,.
\end{equation*}
Therefore, $(\kvz,\lJ)$ is a \emph{compactification} of $\Rdz$, i.e.~$\lJ$ is a smooth embedding 
(injective $\mathrm{C}^\infty$-diffeomorphism onto its image) such that the image
$\lJ(\Rdz)$ is dense in $\kvz$.

For future reference we introduce the following notation:
\begin{equation}\label{eq:set_Sd_I}
\mathrm{S}^d_{I} := \Bigl\{(\zeta_0,\mzeta)\in\mathrm{S}^d : \zeta_0\in I \Bigr\} \;,
\end{equation} 
where $I\subseteq [0,1]$, and the shorthand $\mathrm{S}^d_a:=\mathrm{S}^d_{\{a\}}$.
It is obvious that boundaries of $\kvz$ are given by $\kvn$ and $\kvb$.

In order to relate this compactification to the one used by Tartar \cite[Definition 32.5]{Tar2} mentioned above, we define 
$$
\Sigma_0:=\{0^\ve : \ve\in\Sdmj\} \qquad \hbox{and} \qquad
\Sigma_\infty:=\{\infty^\ve : \ve\in\Sdmj\} \,,
$$
and then we take (a disjoint union) $\kob:=\Rdz\dot\cup\Sigma_0\dot\cup\Sigma_\infty$ (see Figure \ref{fig:top_compt}).
Since the mapping 
\begin{equation}\label{eq:rational}
[0,+\infty)\ni x\mapsto \frac{1}{\sqrt{1+(x+\rnul)^2}}
\end{equation}
is strictly decreasing, it is easy to see that for any sequence $(\mxi_n)$
in $\Rdz$ we have 
\begin{equation*}%\label{TestFjlJ}
\begin{aligned}
\lim_n\Bigl|\lJ(\mxi_n)- \Bigl(\!r_1,\rnul r_1\frac{\mxi_n}{|\mxi_n|}\Bigr)\Bigr|=0
	&\quad\iff\quad \lim_n|\mxi_n|=0 \;,\\
\lim_n\Bigl|\lJ(\mxi_n)-\Bigl(\!0,\frac{\mxi_n}{|\mxi_n|}\Bigr)\Bigr|=0
&\quad\iff\quad \lim_n|\mxi_n|=+\infty \;.
\end{aligned}
\end{equation*}
This motivates the extension of $\lJ$ to the set $\kob$ (and consequently also $\lJ^{-1}$ to $\kvz$) by defining $\lJ(0^\ve):=
(r_1,\rnul r_1\ve)$ and $\lJ(\infty^\ve):=(0,\ve)$, implying that $\lJ(\Sigma_0)=\kvn$ and 
$\lJ(\Sigma_\infty)=\kvb$ (although the notation for the boundaries of
$\kvz$ we found natural, one needs to be careful here as the sphere at infinity is
mapped to $\kvb$).

Now we can define a metric on the compactification $\kob$ by pulling back the Euclidean metric from $\kvz$, 
i.e.
\begin{equation}\label{eq:metric_dob}
d_\kob(\mxi_1,\mxi_2):=|\lJ(\mxi_1)-\lJ(\mxi_2)| \,.
\end{equation}
With these definitions in our hands, it is evident that $(\kob,d_\kob)$ is isometrically isomorphic to 
$\kvz$ ($\lJ$ being the isometric isomorphism),
hence a complete metric space. Note that $\kob$ can be seen as a completion of $(\Rdz,d_\kob)$
and it represents a concrete realisation of the compactification \cite[Definition 32.5]{Tar2} discussed at the beginning of the 
section.

The metric structure allows us to define continuous functions on $\kob$, and it can easily be seen 
that $\psi\in\mathrm{C}(\kob)$ if and only if $\psi\circ\lJ^{-1}\in \mathrm{C}(\kvz)$.

\begin{remark}
The choice of $\rnul>0$ is completely irrelevant as long as we keep it fixed, 
which is the case in the paper. 
Of course, for $\rnul=0$ one can notice that we would have $\lJ=\lR$,
i.e.~$\kob$ would just be the standard radial compactification of $\Rd$.
\end{remark}

\subsection{Smooth test functions}

For simplicity, we shall often denote by
$\psi^\ast:=(\lJ^{-1})^\ast\psi=\psi\circ\lJ^{-1}$ the pullback of function $\psi\in\mathrm{C}(\kob)$ along
$\lJ^{-1}$. 

A common way for defining the space of smooth functions is 
by generalising the previously mentioned 
characterisation for continuous functions on $\kob$: for $\kappa\in\N_0\cup\{\infty\}$ we define
$$
\pC\kappa\kob := \Bigl\{ \psi\in\Cp\kob : \psi^\ast\in\pC\kappa\kvz \Bigr\} \,,
$$   
where $\pC\kappa\kvz$ is the space of those functions which can be extended to a ${\rm C}^\kappa$ 
function on some open neighbourhood of $\kvz$.
Moreover, with this definition $\kob$ is equipped with 
a unique smooth structure such that $\lJ$ is a \emph{smooth covering map} 
(cf.~\cite[Chapter 4]{Lee13}).

%By the mean value theorem we have that in the definition of $\pC\kappa\kob$ it is sufficient 
%(but, of course, not necessary) to require that all derivatives of the order $\kappa+1$ are
%bounded on $\kvz$. 

Since $\lJ:\Rdz\to\kvo$ is a diffeomorphism, the restriction to $\Rdz$ of any 
function from the space $\pC\kappa\kob$ is a function of class ${\rm C}^\kappa$.
%Furthermore, for $\psi\in\pC\kappa\kob$ we have that $\psi_0,\psi_\infty$ 
%(given by \eqref{psinib}) are contained in $\in\pC\kappa\Sdmj$.

For $\kappa\neq\infty$ and $\psi\in\pC\kappa\kob$ we define 
$\|\psi\|_{\pC\kappa\kob} := \|\psi^\ast\|_{\pC\kappa\kvz}$.
Since $\pC\kappa\kvz$ is a Banach algebra, under this definition 
$\pC\kappa\kob$ is also a Banach algebra. 
Moreover, by the compactness of $\kvz$, the space $\pC\kappa\kvz$ is separable, implying 
that $\pC\kappa\kob$ is separable as well.
In the case $\kappa=\infty$ we get a Fr\'echet space by taking a monotone sequence of norms
$\bigl(\|\cdot\|_{\pC{n}{\kob}}\bigr)_n$.

As it was commented above, any $\mathrm{C}^\kappa$ function on $\kob$ is 
in fact of class $\mathrm{C}^\kappa$ on $\Rdz$ (in the sense of the restriction
to $\Rdz$). In the following we study under which additional conditions the converse 
holds as well, i.e.~when a function in $\mathrm{C}^\kappa(\Rdz)$ can be extended to a
function in $\mathrm{C}^\kappa(\kob)$.
Of course, if such extension exists, it is unique (since $\Rdz$ is dense in $\kob$). 
In the rest of the paper we use this identification of $\mathrm{C}^\kappa(\kob)$ with a
subset of $\mathrm{C}^\kappa(\Rdz)$.
 
It is crucial to examine the behaviour of functions around zero and at infinity,
i.e.~in neighbourhoods of $\Sigma_0$ and $\Sigma_\infty$
in $\kob$.
Obviously, neighbourhoods of $\Sigma_0$ and $\Sigma_\infty$
in $\kob$ (with respect to the metric \eqref{eq:metric_dob})
can be identified with the corresponding neighbourhoods of
$\kvn$ and $\kvb$ in $\kvz$, respectively. 
Thus, 
it is sufficient to study functions on 
$\mathrm{S}^d_{\left((1+(1+\rnul)^2)^{-1/2}, \, r_1\right]}$
and
$\mathrm{S}^d_{\left[0, \, (1+(1+\rnul)^2)^{-1/2}\right)}$
(open sets in $\kvz$ containing $\kvn$ and $\kvb$, respectively).
Finally, we transfer to $[0,1)\times\Sdmj$
by using that 
\begin{equation}\label{eq:ljlj0}
\lJ\circ\lJ_0^{-1}:[0,1)\times\Sdmj
\longrightarrow \mathrm{S}^d_{\left((1+(1+\rnul)^2)^{-1/2}, \, r_1\right]}
\end{equation}
and
\begin{equation}\label{eq:ljljinfty}
\lJ\circ\lJ_\infty^{-1} : [0,1)\times\Sdmj \longrightarrow
\mathrm{S}^d_{\left[0, \, (1+(1+\rnul)^2)^{-1/2}\right)}
\end{equation}
are $\mathrm{C}^\infty$-diffeomorphisms, 
where
\begin{equation*}
\lJ_0(\mxi):=\Bigl(|\mxi|,\frac{\mxi}{|\mxi|}\Bigr) \quad \hbox{and} \quad
\lJ_\infty(\mxi):=\Bigl(\frac{1}{|\mxi|},\frac{\mxi}{|\mxi|}\Bigr)\;.
\end{equation*}
Let us prove the latter claim only for $\lJ\circ\lJ_0^{-1}$, 
while the same for $\lJ\circ\lJ^{-1}_\infty$ follows analogously.

We have $\lJ_0^{-1}(s,\meta)=s\meta$ and
$\lJ_0^{-1}([0,1)\times\Sdmj) = \Ko{1}$,
where $\Ko{1}$ is the open ball of radious $1$ centred at the 
origin. Moreover, $\lJ(\Ko{1})=\mathrm{S}^d_{\left[0, \, (1+(1+\rnul)^2)^{-1/2}\right)}$, hence 
the composition \eqref{eq:ljlj0} is well-defined.
Concerning the action we have
\begin{equation*}
\bigl(\lJ\circ\lJ_0^{-1}\bigr)(s,\meta) = \biggl(\frac{1}{\sqrt{1+(s+\rnul)^2}},
\frac{s+\rnul}{\sqrt{1+(s+\rnul)^2}}\meta\biggr) \,,
\end{equation*}
which is obviously a smooth bijection (see \eqref{eq:rational})
from $[0,1)\times\Sdmj$
to $\mathrm{S}^d_{\left((1+(1+\rnul)^2)^{-1/2}, \, r_1\right]}$.
Its inverse is smooth as well since
\begin{equation*}
\bigl(\lJ_0\circ\lJ^{-1}\bigr)(\zeta_0,\mzeta) = \Bigl(\frac{|\mzeta|}{\zeta_0}-\rnul,
\frac{\mzeta}{|\mzeta|}\Bigr)
\end{equation*}
and on its domain it holds $\zeta_0>\frac{1}{\sqrt{1+(1+\rnul)^2}}>0$
and $|\mzeta|\geq \sqrt{1-r_1^2}=\rnul r_1>0$.

The above discussion can be understood so that a (suitable small)
neighbourhood of $0^\ve\in\Sigma_0$ in $\kob$ (and thus of $\mzeta\in\kvn$ in $\kvz$ as well) 
can be identified with a neighbourhood of $(0,\ve)$
in $[0,1)\times\Sdmj$ using the variables $|\mxi|$ and $\frac{\mxi}{|\mxi|}$,
and analogously for $\Sigma_\infty$, but using variables
$\frac{1}{|\mxi|}$ and $\frac{\mxi}{|\mxi|}$.

Now we can prove the following important lemma.

\begin{lemma}\label{lm:smooth_on_kob}
Let $\kappa\in\N_0\cup\{\infty\}$.
The following two statements are equivalent:
\begin{itemize}
\item[i)] $\psi\in\pC\kappa\kob$;
\item[ii)] $\psi\in\pC\kappa\Rdz$ and there exist $\tilde\psi_0,\tilde\psi_\infty\in
	\pC\kappa{[0,1)\times\Sdmj}$ such that
	\begin{align}
	\psi(\mxi) &=\tilde\psi_0\Bigl(|\mxi|,\frac{\mxi}{|\mxi|}\Bigr) \ , \quad 0<|\mxi|<1 \;, 
		\label{eq:tildepsin}\\
	\psi(\mxi) &=\tilde\psi_\infty\Bigl(\frac{1}{|\mxi|},\frac{\mxi}{|\mxi|}\Bigr) \ , \quad 
		|\mxi|>1 \label{eq:tildepsib}\;.
	\end{align}
\end{itemize}
\end{lemma}

\begin{proof}
Let $\psi\in\pC\kappa\kob$. Then $\psi\in\pC\kappa\Rdz$ and 
$\psi=\psi^*\circ\lJ$, where $\psi^*\in\pC\kappa\kvz$
is the pullback. 
Since mappings \eqref{eq:ljlj0}
and \eqref{eq:ljljinfty} are smooth, we have that functions $\tilde\psi_0:=\psi^*\circ\lJ\circ
\lJ_0^{-1}$ and $\tilde{\psi}_\infty:=\psi^*\circ\lJ\circ\lJ_\infty^{-1}$ are
in $\pC\kappa{[0,1)\times\Sdmj}$.
Moreover, $\tilde\psi_0(|\mxi|,\frac{\mxi}{|\mxi|})=\tilde\psi_0(\lJ_0(\mxi))=\psi^*(\lJ(\mxi))=\psi(\mxi)$,
and analogously for $\tilde\psi_\infty$.

Assume now that condition (ii) is fulfilled, and let us prove that 
$\psi\circ\lJ^{-1}\in\pC\kappa\kvz$.
Since $\lJ$ is smooth from $\Rdz$ to $\kvo$, we obviously have $\psi\circ\lJ^{-1}\in\pC\kappa\kvo$.
Thus, it remains to prove the smoothness on a neighbourhood of the boundary. 
As the mapping \eqref{eq:ljlj0} is smooth, we have that function 
$\tilde\psi_0\circ(\lJ_0\circ\lJ^{-1})$ is of class $\mathrm{C}^\kappa$ on
$\mathrm{S}^d_{\left((1+(1+\rnul)^2)^{-1/2}, \, r_1\right]}$.
Furthermore, by \eqref{eq:tildepsin}, the functions $\psi\circ\lJ^{-1}$ and $\tilde\psi_0\circ(\lJ_0\circ\lJ^{-1})$ 
coincide on $\mathrm{S}^d_{\left((1+(1+\rnul)^2)^{-1/2}, \, r_1\right]}$, 
hence $\psi\circ\lJ^{-1}$ is smooth up to the boundary $\kvn$. 
Smoothness on a neighbourhood of $\kvb$ can be obtained in an analogous way by using 
the smoothness of \eqref{eq:ljljinfty} and condition \eqref{eq:tildepsib}.
\end{proof}

If in addition $\psi$ is smooth up to the origin, i.e.~$\psi\in\pC\kappa\Rd$, 
then the function $\tilde\psi(s,\meta):=\psi(s\meta)$ is smooth on $[0,1)\times\Sdmj$ 
and trivially satisfies \eqref{eq:tildepsin}. Thus, we have the following corollary.

\begin{corollary}\label{cor:smooth_on_kob}
Let $\psi\in\pC\kappa\Rd$ and let $\tilde\psi_\infty\in\pC\kappa{[0,1)\times\Sdmj}$ satisfy 
\eqref{eq:tildepsib}. Then $\psi\in\pC\kappa\kob$. 
\end{corollary}

If we focus only on $\Sigma_0$ and $\Sigma_\infty$,
for $\psi\in\pC{\kappa}{\kob}$, 
we can derive from Lemma \ref{lm:smooth_on_kob} 
a simplified necessary condition provided in the following corollary.
In the statement one can see that in the case $\kappa=0$ this is
also a sufficient condition, which was already noted in 
\cite[Definition 32.5]{tar_book} (cf.~\cite{AEL}).

\begin{corollary}\label{cor:psinib}
Let $\psi\in\pC\kappa\kob$. Then there exist unique functions $\psi_0,\psi_\infty\in\pC\kappa\Sdmj$
such that 
\begin{equation}\label{eq:psinib}
\begin{aligned}
\psi(\mxi) &- \psi_0\Bigl(\frac{\mxi}{|\mxi|}\Bigr) \longrightarrow 0\;, \quad |\mxi|\to 0\;, \\
\psi(\mxi) &- \psi_\infty\Bigl(\frac{\mxi}{|\mxi|}\Bigr) \longrightarrow 0\;, \quad |\mxi|\to \infty\;.
\end{aligned}
\end{equation}

Furthermore, if for $\psi\in\Cp\Rdz$ there exist $\psi_0,\psi_\infty\in\Cp\Sdmj$
such that \eqref{eq:psinib} holds, then $\psi\in\Cp\kob$. 
\end{corollary}

\begin{proof}
The uniqueness of functions $\psi_0$ and $\psi_\infty$ trivially follows from \eqref{eq:psinib}.	
	
For the existence, take $\psi\in\pC\kappa\kob$. Then by the previous lemma there exist $\tilde\psi_0,\tilde\psi_\infty
\in\pC\kappa{[0,1)\times\Sdmj}$ such that \eqref{eq:tildepsin}--\eqref{eq:tildepsib} hold. 
For $\meta\in\Sdmj$ we define $\psi_0(\meta):=\tilde\psi_0(0,\meta)$ and 
$\psi_\infty(\meta):=\tilde\psi_\infty(0,\meta)$, which are obviously functions contained 
in $\pC\kappa\Sdmj$. 
Relations \eqref{eq:psinib} now follow by \eqref{eq:tildepsin}--\eqref{eq:tildepsib}
and uniform continuity of functions $\tilde\psi_0$ and $\tilde\psi_\infty$ on $[0,1/2]\times\Sdmj$.

The second part holds by Lemma \ref{lm:smooth_on_kob}, once it is noted that functions
\begin{equation*}
\tilde\psi_0(s,\meta):= 
	\left\{\begin{array}{ll}\psi_0(\meta) \,, & s=0, \, \meta\in\Sdmj\\ 
	\psi(s\meta) \,, & s\in(0,1), \, \meta\in\Sdmj\end{array}\right.
	, \ \;  
	\tilde\psi_\infty(s,\meta):= 
	\left\{\begin{array}{ll}\psi_\infty(\meta) \,, & s=0, \, \meta\in\Sdmj\\ 
	\psi(\frac{\meta}s) \,, & s\in(0,1), \, \meta\in\Sdmj\end{array}\right.
\end{equation*}
are continuous on $[0,1)\times\Sdmj$.
\end{proof}

The characterisation of the space of continuous functions on $\kob$ from the previous 
corollary can easily be derived directly, with no reference to Lemma \ref{lm:smooth_on_kob}, 
which was the case in \cite[Definition 32.5]{tar_book} (see also \cite{AEL}). 
Unfortunately, this characterisation cannot be generalised to the case $\kappa\geq 1$
since smoothness of functions $\psi_0$ and $\psi_\infty$ is not sufficient to ensure 
the required smoothness of functions $\tilde\psi_0$ and $\tilde\psi_\infty$, defined in the proof
of the corollary. 
Indeed, in order to smoothly extend a function to the boundary, 
besides knowledge of tangential derivatives, one needs some information 
on the normal derivative as well.

In the rest of the paper we shall, for an arbitrary $\psi\in\pC\kappa\kob$, denote by 
$\psi_0$ and $\psi_\infty$ the functions defined by \eqref{eq:psinib},
i.e.~the restrictions of $\psi$ to $\Sigma_0$ and $\Sigma_\infty$ respectively.

Let us see how we can estimate the norm of $\psi\in\pC\kappa\kob$.
The interval $[0,r_1]$ is equal to the union of the following 
closed intervals:
\begin{align*}
I_1 &:=\left[0,\frac{1}{\sqrt{1+(2+\rnul)^2}}\right], \
	I_2:=\left[\frac{1}{\sqrt{1+(4+\rnul)^2}}, \frac{1}{\sqrt{1+(1/4+\rnul)^2}}\right], \\
I_3 &:=\left[\frac{1}{\sqrt{1+(1/2+\rnul)^2}}, r_1\right],
\end{align*}
where we intentionally made the intervals overlap.
Now, in order to derive an estimate of $\|\psi\|_{\pC{\kappa}{\kob}}$, 
it is enough to provide an estimate for the 
$\mathrm{C}^\kappa$ norm of $\psi^*$ on each of the compact sets
$\mathrm{S}^d_{I_i}$, $i\in\{1,2,3\}$. 
Indeed, it is obvious that $\kvz$ is equal 
to the union of sets $\mathrm{S}^d_{I_i}$, $i\in\{1,2,3\}$, 
while the fact that $\mathrm{S}^d_{I_2}$ overlaps with the other two 
on a \emph{thick} set is important for controlling derivatives 
globally, i.e.~on the whole $\kvz$.

For the estimate on $\mathrm{S}^d_{I_3}$ and $\mathrm{S}^d_{I_1}$
we use \eqref{eq:tildepsin} and \eqref{eq:tildepsib},
respectively, while on $\mathrm{S}^d_{I_2}$ 
the estimate can easily be obtained as we are on a compact set far from the boundary.
Before we proceed with calculations, let us note that 
$$
(\lJ_0\circ\lJ^{-1})(\mathrm{S}^d_{I_3})
=(\lJ_\infty\circ\lJ^{-1})(\mathrm{S}^d_{I_1})=[0,1/2]\times\Sdmj
$$
(hence Lemma \ref{lm:smooth_on_kob} is applicable),
and
$\lJ^{-1}(\mathrm{S}^d_{I_2})=\{\mxi\in\Rd : \frac{1}{4}\leq|\mxi|\leq 4\}$.

Let us start with the estimate on $\mathrm{S}^d_{I_3}$ 
(i.e.~close to $\Sigma_0$). 
By Lemma \ref{lm:smooth_on_kob} we have that the
equality $\tilde\psi_0\circ(\lJ_0\circ\lJ^{-1}) = \psi\circ\lJ^{-1}$ holds 
on $\mathrm{S}^d_{I_3}$.
Since $\lJ_0\circ\lJ^{-1}$ is a $\mathrm{C}^\infty$-diffeomorphism 
on that compact set, this function and all its derivatives are bounded. 
Now we apply the generalised chain rule formula, 
known under the name \emph{Fa\'a di Bruno formula}  (see
e.g.~\cite{FaadiBruno}):
{\it For sufficiently smooth functions $\vg:\Rd\to\R^r$ and 
	$f:\R^r\to\R$ one has}
\begin{equation}\label{FaaDiBruno}
\partial^\malpha(f\circ \vg)(\mxi) = |\malpha|! \sum_{\substack{\mbeta\in\N_0^r,\\
		1\leq |\mbeta|\leq |\malpha|}} C(\mbeta,\malpha) \,,
\end{equation}
{\it where}
$$
C(\mbeta,\malpha) = {(\partial^\mbeta f)(\vg(\mxi))\over \mbeta!} 
\sum_{\substack{\malpha_i\in\N_0^d,\\ 
		\sum_{i=1}^r\malpha_i = \malpha}}\; \prod_{j=1}^r \; 
\sum_{\substack{\mgamma_i\in\N_0^d\setminus\{\vnul\},\\
		\sum_{i=1}^{\beta_j}\mgamma_i=\malpha_j }} \; \prod_{s=1}^{\beta_j} 
{\partial^{\mgamma_s}g_j(\mxi)\over \mgamma_s!} \;.
$$
Therefore, there exists a constant $C_\kappa>0$ depending only on $\kappa$ such that
$$
\begin{aligned}
\|\psi\circ\lJ^{-1}\|_{\pC\kappa{\mathrm{S}^d_{I_3}}}
	&= \|\tilde\psi_0\circ(\lJ_0\circ\lJ^{-1})\|_{\pC\kappa{\mathrm{S}^d_{I_3}}} \\
&\leq C_\kappa \|\tilde\psi_0\|_{\pC\kappa{[0,\frac{1}2]\times\Sdmj}} \;.
\end{aligned}
$$ 
Applying the same reasoning on $\tilde\psi_0=(\psi\circ\lJ^{-1})\circ(\lJ\circ\lJ_0^{-1})$
we get 
$$
\begin{aligned}
\|\psi\|_{\pC\kappa\kob}&=\|\psi\circ\lJ^{-1}\|_{\pC\kappa\kvz}\\
&\geq \|\psi\circ\lJ^{-1}\|_{\pC\kappa{\mathrm{S}^d_{I_3}}}
	\geq c_\kappa \|\tilde\psi_0\|_{\pC\kappa{[0,\frac{1}2]\times\Sdmj}} \;.
\end{aligned}
$$ 

Completely analogous estimates could be derived on 
$\mathrm{S}^d_{I_1}$ for $\tilde\psi_\infty$. 
Finally, since $\lJ$ is a $\mathrm{C}^\infty$-diffeomorphism on the compact 
$\{\mxi\in\Rd : \frac{1}{4}\leq|\mxi|\leq 4\}$, we can estimate the norm of 
$\psi\circ\lJ^{-1}$ on 
$\mathrm{S}^d_{I_2}$
(both from below and above) by the norm of $\psi$ on 
$\{\mxi\in\Rd : \frac{1}{4}\leq|\mxi|\leq 4\}$ (up to a multiplicative constant).

The overall conclusion we give in the form of a lemma.  

\begin{lemma}\label{lm:psi_norm}
For any $\kappa\in\N_0$ there exist constants $c_\kappa, C_\kappa>0$,
depending only on $\kappa$, such that for any $\psi\in\pC\kappa\kob$ it holds
\begin{equation*}
\begin{aligned}
c_\kappa & \max\bigl\{\|\tilde\psi_0\|_{\pC\kappa{[0,\frac{1}2]\times\Sdmj}},
	\|\psi\|_{\pC\kappa{\{\mxi\in\Rd : \frac{1}{4}\leq|\mxi|\leq 4\}}},
	\|\tilde\psi_\infty\|_{\pC\kappa{[0,\frac{1}2]\times\Sdmj}}\bigr\} \\
&\leq \|\psi\|_{\pC\kappa\kob} \\
& \leq C_\kappa\max\bigl\{\|\tilde\psi_0\|_{\pC\kappa{[0,\frac{1}2]\times\Sdmj}},
\|\psi\|_{\pC\kappa{\{\mxi\in\Rd : \frac{1}{4}\leq|\mxi|\leq 4\}}},
\|\tilde\psi_\infty\|_{\pC\kappa{[0,\frac{1}2]\times\Sdmj}}\bigr\} \;,
\end{aligned}
\end{equation*}
where functions $\tilde\psi_0$ and $\tilde\psi_\infty$ are given by
\eqref{eq:tildepsin}--\eqref{eq:tildepsib}.
\end{lemma}

Let us now present some examples of functions inhabiting the space $\pC\kappa\kob$. 

\begin{corollary}\label{cor:examples}
Let $\mpi(\mxi):=\frac{\mxi}{|\mxi|}$ be the projection of $\Rdz$ along rays through origin to the 
unit sphere $\Sdmj$. Then the following statements hold:
\begin{itemize}
\item[i)] $\bigl\{\psi\circ\mpi : \psi\in\pC\kappa\Sdmj\bigr\}\subseteq \pC\kappa\kob$,
	$\kappa\in\N_0\cup\{\infty\}$,
\item[ii)] $\lS(\Rd)\subseteq\pC\infty\kob$,
\item[iii)] $(\forall m,l\in\N_0)(\forall \malpha\in\N_0^d)\quad l\leq|\malpha|\leq m \Longrightarrow
 \mxi\mapsto\frac{\mxi^\malpha}{|\mxi|^l+|\mxi|^m}\in\pC\infty\kob$, \quad and
\item[iv)] $(\forall m\in\N)$ $\ \mxi\mapsto\frac{1+|\mxi|^m}{(1+|\mxi|^2)^{\frac{m}2}}\,,\;
\mxi\mapsto\frac{(1+|\mxi|^2)^{\frac{m}{2}}}{1+|\mxi|^m}\in\pC\infty\kob$.
\end{itemize}
\end{corollary}

\begin{proof}
(i) For any $\psi\in\pC\kappa\Sdmj$ the function $\mxi\mapsto\psi(\frac{\mxi}{|\mxi|})$ trivially 
satisfies conditions of the second part of Lemma \ref{lm:smooth_on_kob} with 
$\tilde\psi_0(s,\meta)=\tilde\psi_\infty(s,\meta)=\psi(\meta)$, for $(s,\meta)\in[0,1)\times\Sdmj$,
i.e.~$\psi_0=\psi_\infty=\psi$.  
\smallskip

\noindent (ii) By Corollary \ref{cor:smooth_on_kob} we only have to show the existence of 
$\tilde\psi_\infty$. Let $\psi\in\lS(\Rd)$ and let us define 
$\tilde\psi_\infty(s,\meta):=\psi(\frac{\meta}{s})$, for $(s,\meta)\in(0,1)\times\Sdmj$.
It is obvious that $\tilde\psi_\infty\in\pC\infty{(0,1)\times\Sdmj}$ and that 
\eqref{eq:tildepsib} is satisfied. Hence, it remains to be proven that $\tilde\psi_\infty$
and all its derivatives can be extended by continuity on $s=0$. 
This holds since $\psi$ is a rapidly decreasing function and for an arbitrary 
$\malpha\in\N_0^d$ we have $\partial^\malpha \tilde\psi(s,\meta)
= (\partial^\malpha\psi)(\frac{\meta}{s}) P_\malpha(\frac{1}{s},\meta)$, 
where $P_\malpha$ is a polynomial (in $1+d$ variables).
Therefore, by letting $\partial^\malpha\tilde\psi_\infty(0,\meta):=0$, $\malpha\in\N_0^d$,
$\meta\in\Sdmj$, we get that $\tilde\psi_\infty$ is smooth up to the boundary $s=0$.
\smallskip

\noindent (iii)
Since $\frac{\mxi^\malpha}{|\mxi|^l+|\mxi|^m} = \frac{\mxi^\malpha}{|\mxi|^{|\malpha|}}
\frac{|\mxi|^{|\malpha|-l}}{1+|\mxi|^{m-l}}$ and $\pC\infty\kob$ is an algebra, 
by part (i) it is sufficient to prove that for any $m,k\in\N_0$, $k\leq m$, the
function $\psi^{m,k}(\mxi)=\frac{|\mxi|^k}{1+|\mxi|^m}$ is contained in $\pC\infty\kob$.
Let us now check the conditions of Lemma \ref{lm:smooth_on_kob}(ii). 

It is evident that $\psi^{m,k}$ is smooth on $\Rdz$. Moreover, since 
$$
\psi^{m,k}(\mxi) = \frac{|\mxi|^k}{1+|\mxi|^m} = 	
	\frac{(\frac{1}{|\mxi|})^{m-k}}{(\frac{1}{|\mxi|})^m+1} \ , \quad \mxi\in\Rdz \;,
$$
the functions $\tilde\psi^{m,k}_0(s,\meta)=\frac{s^k}{1+s^m}$
and $\tilde\psi_\infty(s,\meta)=\frac{s^{m-k}}{s^m+1}$ satisfy 
\eqref{eq:tildepsin}--\eqref{eq:tildepsib}, and they are clearly smooth on $[0,1)\times\Sdmj$. 
Thus, by Lemma \ref{lm:smooth_on_kob} we have that $\psi^{m,k}\in\pC\infty\kob$. 
\smallskip

\noindent (iv) By Lemma \ref{lm:smooth_on_kob} the statements trivially holds, following the same 
argument as in part (iii). Let us just remark that it is sufficient to prove the 
statement only for one function after noting that it is strictly positive (hence we can conclude that its
multiplicative inverse is in $\pC\infty\kob$ as well). 
\end{proof}

\begin{remark}
\item[i)] The compactification $\kob$ could be seen as the blow-up at $\vnul$ of the radial 
compactification of $\Rd$, or the blow-up at $\vnul$ and $\infty$ of the one-point 
compactification of $\Rd$ (see e.g.~\cite[Section 1+.5]{RM_lect}). 
If we denote the radial compactification by $\kb$ (this notation was used in 
\cite[Definition 32.5]{tar_book}), this means that the mapping 
$\mathcal{I}:\kob\to\kb$ given by
$$
\mathcal{I}(\mxi) = \left\{\begin{array}{ll}\mx \,, & \mx\not\in\Sigma_0 \\ \vnul \,, & \mx\in\Sigma_0\end{array}\right.
$$
is smooth. In particular, all smooth functions on the radial compactification
are smooth on $\kob$, which in fact we have already shown in Corollary \ref{cor:smooth_on_kob}
(in \cite[Chapter 1]{RM_lect} the smooth functions on the radial compactification are 
characterised precisely by the assumptions given in the corollary).

\item[ii)] Instead of \eqref{eq:lj}, for the compactifying map in \cite[Section II.1]{MEphd}  
the projection on the last $d$ variables of $\lJ(\mxi)$ was taken. 
That choice corresponds to a generalisation of the quadratic compactification 
(see e.g.~\cite[(L1.11)]{RM_lect}; for $\rnul=0$ it is precisely the quadratic compactification).
These two compactifications are not equivalent since the inverse of the projection 
of the upper half-sphere to the ball has a square-root singularity \cite[Chapter 1]{RM_lect}, 
thus in \cite[Section III.2]{MEphd} a smaller space of smooth functions was obtained.
In particular, the function $\mxi\mapsto(1+|\mxi|^2)^{-\frac{1}{2}}$ is not smooth with that choice, 
which caused some problems in deriving the localisation principle \cite[Section III.5]{MEphd}. 

\item[iii)] The idea of studying functions on a compactified space
could be seen as a way of controlling functions at infinity \cite{RM_lect}
(in our case of $\kob$ at the origin as well; see \eqref{eq:mihlin_psi} below).
Thus, with this approach the standard growth conditions of (pseduodifferential)
symbols could be set intrinsically by a suitable choice of compactification, which is a known approach in the literature.
\end{remark}

\subsection{Symbols of Fourier multipliers}

In the construction of one-scale H-distributions we shall need that test functions
on the dual space are \emph{Fourier multipliers}; i.e.~a $\psi\in\mathrm{L}^\infty(\Rd)$ is an
${\rm L}^p$ Fourier multiplier ($\psi\in\lM_p$) if the linear operator $\lA_\psi$ obtained
by multiplying the Fourier transform of a function by $\psi$, and then taking the inverse Fourier transform
$$
\lA_\psi\vu := (\psi\hat\vu)^\vee
$$
is a bounded operator on $\mathrm{L}^p(\Rd;\Cr)$ \cite[Section 2.5]{Grafakos} (see also \cite{IvanNenad}).
The norm of such $\psi\in\lM_p$ is by definition equal to the norm of the operator
$\lA_\psi$ on $\mathrm{L}^p(\Rd;\Cr)$.

Clearly, $\lM_2=\mathrm{L}^\infty(\Rd)$, and $\psi\mapsto \lA_\psi$ is a homomorphism of algebras (with norm equal to 1).
Furthermore, it is translation invariant, and such operators can be written in the convolution form $\check\psi\ast\vu$.
Similarly,  $\lM_1$ is the image by Fourier transform of the space of all bounded measures on $\Rd$.

It is a consequence of the Riesz-Thorin interpolation theorem that for $1\le p\le q\le2$ we have
$$
\lM_1\subseteq \lM_p \subseteq \lM_q \subseteq \lM_2 \;.
$$
For $2<r=p'<\infty$ we use the fact that $\lA_\psi$, for $\psi\in\lM_r = \lM_{p'}$, is also well defined
on $\mathrm{L}^p(\Rd;\Cr)$ (essentially, in convolution form, the adjoint operator is given by convolution
with $\psi$ replaced by $\overline{\tilde\psi}$, where $\tilde\psi(\mx)=\psi(-\mx)$), and that it has the same norm.

Our goal is to show that  there exists $\kappa$ such that for any 
$p\in(1,\infty)$ all functions from $\pC\kappa\kob$ are $\mathrm{L}^p$ 
Fourier multipliers (in the sense of the corresponding
restrictions to $\Rdz$). In the case $p=2$ this is trivially satisfied 
(see \cite[Theorem 2.5.10]{Grafakos}) since 
continuous functions on $\kob$ are bounded.
However, for $p\neq 2$ the complete answer to the question of characterising 
$\mathrm{L}^p$ Fourier multipliers is still open. Thus, one relies on 
a few available criteria providing sufficient conditions. 
Here we shall use the Mihlin condition \cite[Theorem 5.2.7]{Grafakos}
(note that $\lfloor r \rfloor$ denotes the greatest integer less than or equal to real $r$).

\begin{theorem}\label{thm:mihlin}
Let $\psi\in\pC{\lfloor\frac{d}{2}\rfloor+1}\Rdz$
be a bounded function. Assume that there exists $C>0$ such that for any 
$\malpha\in\N_0^d$, $|\malpha|\leq \lfloor\frac{d}{2}\rfloor+1$, it holds
\begin{equation}\label{eq:mihlin}
|\partial^\malpha\psi(\mxi)| \leq {C\over|\mxi|^{|\malpha|}} \ , \quad
	\mxi\in\Rdz \,.
\end{equation}
Then for any $p\in(1,\infty)$ the function $\psi$ is an $\mathrm{L}^p$ Fourier 
multiplier and the operator norm of $\lA_\psi$ satisfies
$$
\|\lA_\psi\|_{\lL(\mathrm{L}^p(\Rd))} \leq C_{d,p} C \,,
$$
where $C_{d,p}$ is a constant depending only on $d$ and $p$, while $C$ is the constant
appearing in \eqref{eq:mihlin}.
\end{theorem}

Now we finally prove the main result of this section. 

\begin{theorem}\label{thm:multipliers}
Any function from $\pC{\lfloor\frac{d}{2}\rfloor+1}{\kob}$ 
satisfies Mihlin's condition \eqref{eq:mihlin}. 

In particular, for any 
$p\in(1,\infty)$, it holds
\begin{equation*}
\|\lA_\psi\|_{\lL(\mathrm{L}^p(\Rd))} \leq C_{d,p}
	C_d\|\psi\|_{\pC{\lfloor\frac{d}{2}\rfloor+1}{\kob}}
	\ , \quad \psi\in\pC{\lfloor\frac{d}{2}\rfloor+1}{\kob} \;,
\end{equation*}
where $C_{d,p}$ is the constant from Theorem \ref{thm:mihlin}, and 
$C_d$ is a constant depending only on $d$.

\end{theorem}
\begin{proof}
Let $\psi\in\pC{\lfloor\frac{d}{2}\rfloor+1}{\kob}$.
By Theorem \ref{thm:mihlin} it is sufficient to prove
\begin{equation}\label{eq:mihlin_psi}
|\partial^\alpha\psi(\mxi)| \leq \frac{C_d 
	\|\psi\|_{\pC{\lfloor\frac{d}{2}\rfloor+1}{\kob}}}{|\mxi|^{|\malpha|}}
\ , \quad \mxi\in\Rdz \;,
\end{equation}
for any $|\malpha|\leq \lfloor\frac{d}{2}\rfloor+1$.

Let us consider first $\mxi\in\Rdz$ such that $\frac{1}{4}\leq|\mxi|\leq 4$.	
Since this set is compact and both functions $\partial^\malpha\psi$ and 
$\mxi\mapsto\frac{1}{|\mxi|^{|\malpha|}}$ are continuous on it
(hence bounded), 
one can easily see that there exists a constant $C_d>0$ such that 
\eqref{eq:mihlin_psi} holds for any $\mxi\in\Rdz$, $\frac{1}{4}\leq|\mxi|\leq 4$. 

For $|\mxi|\in(0,\frac{1}{2})$, by Lemma \ref{lm:smooth_on_kob} there exists 
$\tilde\psi_0\in\pC{\lfloor\frac{d}{2}\rfloor+1}{[0,1)\times\Sdmj}$ such that
$\psi(\mxi)=\tilde\psi_0(|\mxi|,\frac{\mxi}{|\mxi|})$.

Since $\mxi\mapsto\frac{\mxi}{|\mxi|}$ is homogeneous of degree zero,
each of its components satisfies Mihlin's condition \eqref{eq:mihlin}
(see e.g.~\cite[Section 5.2.2]{Grafakos}). 
For $\mxi\mapsto|\mxi|$ we have $\partial_j|\mxi|=\frac{\xi_j}{|\mxi|}$, 
thus by the previous observation for $\malpha\in\N_0^d$, $|\malpha|\geq 1$,
we have
$$
\partial^\malpha|\mxi| \leq \frac{C}{|\mxi|^{|\malpha|-1}} \leq \frac{C}{|\mxi|^{|\malpha|}} 
	\ , \quad  |\mxi|\in\Bigl(0,\frac{1}{2}\Bigr) \;,
$$
where in the second inequality we have used that $\frac{1}{|\mxi|}>1$ on 
the observed set. 
Therefore, each component of the function $\lJ_0(\mxi)=(|\mxi|,\frac{\mxi}{|\mxi|})$
satisfies Mihlin's condition. Then by the Fa\'a di Bruno formula \eqref{FaaDiBruno} for
$|\mxi|\in(0,\frac{1}{2})$ we get
\begin{align*}
|\partial^\malpha\psi(\mxi)|\leq \frac{\tilde C_d}{|\mxi|^{|\malpha|}} |\malpha|!
	\sum_{\substack{\mbeta\in\N_0^d,\\ 1\leq |\mbeta|\leq |\malpha|}}
	\frac{\bigl|(\partial^\mbeta\tilde\psi_0)\bigl(\lJ_0(\mxi)\bigr)\bigr|}{\mbeta!}
	\leq \frac{C_d 
		\|\psi\|_{\pC{\lfloor\frac{d}{2}\rfloor+1}{\kob}}}{|\mxi|^{|\malpha|}} \,,
\end{align*}
where in the last inequality we have used Lemma \ref{lm:psi_norm}, while 
$\tilde C_d$ and $C_d$ are constants dependending only on the dimension of space $d$.

The remaining case $|\mxi|\in (2,\infty)$ follows in the same manner since 
$\psi(\mxi)=\tilde\psi_\infty(\lJ_\infty(\mxi))$ in this case, and all components of 
$\lJ_\infty$ satisfy the Mihlin condition on this set (it is important to take into account that 
$|\mxi|>2$). 
\end{proof}

By the previous theorem we have that, for any $d\geq 1$ and 
$p\in (1,\infty)$, the linear mapping
$$
\pC{\lfloor\frac{d}{2}\rfloor+1}{\kob}\ni\psi\mapsto\lA_\psi\in\lL(\mathrm{L}^p(\Rd))
$$ 
is continuous.

\subsection{Commutation lemma}

For $\ph\in\mathrm{L}^\infty(\Rd)$, let us denote by 
$B_\ph \vu = \ph\vu$ the operator of multiplication by $\ph$.
This operator is clearly bounded on $\mathrm{L}^p(\Rd;\Cr)$, for any 
$p\in [1,\infty]$, with the norm being equal to the $\mathrm{L}^\infty$ 
norm of $\ph$.

Thus, by Theorem \ref{thm:multipliers} 
(see also Lemma \ref{lm:pom_mihlin} below),
for $\psi\in \pC{\lfloor{d\over 2}\rfloor+1}\kob$ 
and $\omega_n\to 0^+$ the sequence of commutators
$$
[B_\ph,\lA_{\psi_n}] := B_\ph\lA_{\psi_n}-\lA_{\psi_n}B_\ph \;,
$$ 
where $\psi_n(\mxi) := \psi(\omega_n\mxi)$,
is bounded in $\lL(\mathrm{L}^p(\Rd))$, for any $p\in (1,\infty)$.
The above sequence of operators (under an additional assumption on $\ph$)
has a special structure, which can be seen in \cite[Lemma 3]{AEL}
for the $\mathrm{L}^2$ case. 
Here we generalise that result by the interpolation argument to the
$\mathrm{L}^p$, $p\in(1,\infty)$, setting.

\begin{lemma}\label{lm:commutation}
	Let $\psi\in\pC{\lfloor{d\over 2}\rfloor+1}\kob$, 
	$\ph\in\mathrm{C}_0(\Rd)$, $\omega_n\to0^+$, and denote $\psi_n(\mxi) := \psi(\omega_n\mxi)$. 
	Then the commutator of multiplication and the Fourier multiplier can be expressed as a sum
	$$
	C_n := [B_\ph,\lA_{\psi_n}] = \tilde C_n +  K \,,
	$$
	where for any $p\in(1,\infty)$ we have that $K$ is a compact operator on $\mathrm{L}^p(\Rd)$, 
	while $\tilde C_n\to 0$ in the operator norm on $\lL(\mathrm{L}^p(\Rd))$.
\end{lemma}
\begin{proof} 
	Since 
	$$
	\lA_{\psi_n} = \lA_{\psi_n-\psi_0\circ\mpi}+\lA_{\psi_0\circ\mpi}
	$$
	(recall that $\psi_0$ is given by \eqref{eq:psinib}),
	we have that $C_n = \tilde C_n + K$, where $\tilde C_n:=[B_\ph,\lA_{\psi_n-\psi_0
		\circ\mpi}]$ and $K := [B_\ph, \lA_{\psi_0\circ\mpi}]$.
	By the standard First commutation lemma \cite[Lemma 1.7]{Tar} $K$
	is a compact operator on $\mathrm{L}^2(\Rd)$, while applying \cite[Lemma 32.4]{tar_book}
	$\tilde C_n\to 0$ in the operator norm on $\lL(\mathrm{L}^2(\Rd))$. 
	
	As $\psi_0\in\pC{\lfloor{d\over 2}\rfloor+1}\Sdmj$, we have
	$\psi_0\circ\mpi\in\pC{\lfloor{d\over 2}\rfloor+1}\kob$ 
	(see Corollary \ref{cor:examples}(i)). 
	Thus, $\psi_0\circ\mpi$ is an $\mathrm{L}^p$ Fourier multiplier for any $p\in(1,\infty)$
	(this is, of course, a well known fact; see e.g.~\cite[Section 5.2.2]{Grafakos}),
	which implies that for any $p\in(1,\infty)$ the operator $K$ is bounded on $\mathrm{L}^p(\Rd)$ 
	(the bound is not uniform in $p$, which does not cause any problems).
	On the other hand, since multipliers are invariant under dilations
	(see \cite[Proposition 2.5.14]{Grafakos} or Lemma \ref{lm:pom_mihlin} 
	below),
	by Theorem \ref{thm:multipliers}
	%by Lemma \ref{lm:pom_mihlin} (see below) 
	the sequence of operators 
	$(\tilde C_n)$ is uniformly (with respect to $n$) bounded on $\mathrm{L}^p(\Rd)$,
	$p\in(1,\infty)$.
	
	Let $p\in(1,\infty)$ be arbitrary, but fixed. Then there exist $r\in(1,\infty)$ and
	$\theta\in(0,1)$ such that $1/p = \theta/2 + (1-\theta)/r$. Since $K$ is compact on 
	$\mathrm{L}^2(\Rd)$ and bounded on $\mathrm{L}^{r}(\Rd)$, by 
	the result due to Krasnosel'skij (see \cite[Lemma 3]{AEM}, \cite{1stcommlemm} and 
	references therein) 
	we get that $K$ is compact on $\mathrm{L}^p(\Rd)$.
	
	Similarly, by using the Riesz-Thorin interpolation theorem (see 
	e.g.~\cite[Theorem 1.3.4]{Grafakos}) we obtain
	$$
	\|\tilde C_n\|_{\lL(\mathrm{L}^p(\Rd))} \leq \|\tilde C_n\|_{\lL(\mathrm{L}^2(\Rd))}^\theta
	\|\tilde C_n\|_{\lL(\mathrm{L}^{r}(\Rd))}^{1-\theta} \,,
	$$
	implying $\tilde C_n\to 0$ in the operator norm on $\mathrm{L}^p(\Rd)$.
\end{proof}

\begin{remark}\label{rem:comm}
\begin{itemize}
\item[i)] For $\psi = \tilde\psi\circ\mpi$, where 
$\tilde\psi\in\pC{\lfloor{d\over 2}\rfloor+1}\Sdmj$ (see Corollary 
\ref{cor:examples}(i)), it is obvious that $\tilde  C_n=0$. 
Thus, the previous lemma can be seen as a generalisation
of \cite[Corollary 3]{AEM} as well.
\item[ii)]  For $\psi\in\pC{\lfloor{d\over 2}\rfloor+1}{\kob}\cap\Cp{\Rd}$,
in the previous lemma we can apply
\cite[Lemma 32.4]{tar_book} immediately on $C_n$, 
since $\psi$ is bounded and uniformly continuous on $\Rd$. 
Thus, there is no need for subtracting $\psi_0\circ\mpi$. 
Therefore, in this case the result is that for any 
$p\in (1,\infty)$ we have  $C_n\to 0$ in the operator norm on $\lL(\mathrm{L}^p(\Rd))$.
\end{itemize}
\end{remark}

The following simple result we provide for completeness. 

\begin{lemma}\label{lm:pom_mihlin}
	Let $\psi\in\mathrm{L}^\infty(\Rd)\cap 
	\pC{\lfloor{d\over 2}\rfloor+1}{\Rdz}$ satisfy Mihlin's condition:
	$$
	(\forall\mxi\in\Rdz)(\forall\malpha\in\N_0^d) \quad |\malpha|\leq \Bigl\lfloor{d\over 2}\Bigr\rfloor+1 
	\implies |\partial^\malpha\psi(\mxi)| \leq {C\over|\mxi|^{|\malpha|}} \,.
	$$
	Then for any $a>0$ the function $\psi_a:=\psi(a\,\cdot)$ also satisfies the Mihlin condition 
	with the same constant $C$. 
\end{lemma}
\begin{proof}
	For $|\malpha|\leq\lfloor{d\over 2}\rfloor+1$ and $\mxi\in\Rdz$ we have
	$$
	|\partial^\malpha\psi_a(\mxi)| = a^{|\malpha|} |(\partial^\malpha\psi)(a\mxi)|
	\leq a^{|\malpha|}{C\over|a\mxi|^{|\malpha|}} = {C\over|\mxi|^{|\malpha|}} \,,
	$$
	proving the claim.
\end{proof}

\section{Kernel theorem}

%Before the existence result of one-scale H-distributions, we need to 
%ensure some auxiliary results. In fact, we need a suitable variant of the First 
%commutation lemma, as well as the right form of the Schwartz kernel theorem.

\subsection{Anisotropic distributions on manifolds without boundary}

In the next section we shall study certain Schwartz' distributions 
\cite{Schwartz} of finite order
in both variables $\mx$ (the physical variable) and $\mxi$ (the dual variable)
on a manifold with boundary. 
Although such distributions are not in general Radon measures 
(distributions of order 0), we shall be able to prove that at least the order 
with respect to $\mx$ is 0. 
In order to do so, we need to work with distributions of 
\emph{anisotropic} order, which have recently been introduced 
on manifolds without boundary in
\cite{AEM} as \emph{anisotropic distributions}.
Let us start by recalling the definition of these objects. 

For simplicity, let us explain the construction for an
open subset $\Omega\subseteq \mathbb{R}^d_\mx\times\mathbb{R}^r_\my$,
while the case where $\Rd$ and $\mathbb{R}^r$ are replaced by 
differentiable manifolds without boundary $X$ and $Y$ then easily follows 
using the local nature of distributions and the fact that every 
differentiable manifold is locally diffeomorphic to some Euclidean
space.

For $l,m\in \Nnb$ we introduce the following space of functions:
\begin{equation*}
\mathrm{C}^{l,m}(\Omega) := \Bigl\{f:\Omega\to\C : (\forall\malpha\in\N_0^d)
(\forall\mbeta\in\N_0^r) \ |\malpha|\leq l \,\&\, |\mbeta|\leq m \implies
\partial^\malpha_\mx \partial^\mbeta_\my f\in\Cp{\Omega}\Bigr\} \;,
\end{equation*} 
where the partial order on $\N_0$ is naturally extended by 
$k\leq\infty$, $k\in \N_0\cup\{\infty\}$. 
Thus, ${\rm C}^{0,0}(\Omega) = {\rm C}(\Omega)$ and 
${\rm C}^{\infty,\infty}(\Omega)={\rm C}^{\infty}(\Omega)$,
but in general we have only
$$
\mathrm{C}^{l+m}(\Omega)\subseteq \mathrm{C}^{l,m}(\Omega)
\subseteq \mathrm{C}^{\min\{l,m\}}(\Omega) \,,
$$
where we extend the addition to $\N_0\cup\{\infty\}$ by taking $k+\infty=\infty$, for any $k\in\Nnb$.

Take a sequence of nested compact sets $K_n$ in $\Omega$, such that $\Omega = \bigcup_{n\in\N}K_n$ 
and $K_n\subseteq\Int K_{n+1}$ (by $\Int A$ we denote the interior of set $A$), 
and for $l,m\in\N_0$ and
$f\in {\rm C}^{l,m}(\Omega)$ we define the seminorms
\begin{equation}\label{eq:seminormRd}
p^{l,m}_{K_n}(f):=\max_{|\malpha|\leq l, |\mbeta|\leq m}
\|\partial^{\malpha,\mbeta} f \|_{{\rm L}^\infty(K_n)} \;.
\end{equation}
If $l=\infty$, instead of seminorms in \eqref{eq:seminormRd}, 
we take the following sequence of seminorms on ${\rm C}^{\infty,m}(\Omega)$:
$$
p^{n,m}_{K_n}(f):=\max_{|\malpha|\leq n, |\mbeta|\leq m}\|\partial^{\malpha,\mbeta} 
f \|_{{\rm L}^\infty(K_n)} \;,
$$
and similarly for $m=\infty$ (of course, the choice $l=m=\infty$ reduces to the standard isotropic case of ${\rm C}^{\infty}(\Omega)$).
For $l,m\in\N_0\cup\{\infty\}$ these seminorms turn ${\rm C}^{l,m}(\Omega)$ 
into a separable Fr\'{e}chet space with
the topology of uniform convergence on compact sets
of functions and their derivatives up to order $l$ in $\mx$ and $m$ in $\my$,
while $\mathrm{C}^\infty_c(\Omega)$ is dense in it
(cf.~\cite[Theorem 1]{AEM}).
For a compact set $K\subseteq\Omega$ and finite $l$ and $m$, its subspace
$$
{\rm C}^{l,m}_K(\Omega):=\Big\{ f\in{\rm C}^{l,m}(\Omega):\; \supp{f}\subseteq K \Big\}
$$
is a Banach space, and the inherited topology from ${\rm C}^{l,m}(\Omega)$
is a norm topology determined by
$$
\| f\|_{l,m,K}:= p^{l,m}_K(f)\;.
$$
If $l=\infty$ or $m=\infty$, 
then we shall not get a Banach space, but a Fr\'{e}chet space.
Finally, the set of all $\mathrm{C}^{l,m}(\Omega)$ functions 
with compact support
$$
{\rm C}^{l,m}_c(\Omega):= \bigcup_{n\in\N}{\rm C}^{l,m}_{K_n}(\Omega)
$$
we equip with the topology of strict inductive limit, obtaining a complete 
locally convex topological vector space.

\begin{definition}\label{defn:aniso-distr}
	Any continuous linear functional on ${\rm C}^{l,m}_c(\Omega)$,
	where $\Omega\subseteq\Rd\times\Rr$ is an open set, we call 
	a \emph{distribution of anisotropic order}, and
	such functionals form a vector space denoted by
	$\lD'_{l,m}(\Omega):=({\rm C}^{l,m}_c(\Omega))'$.
	%
	%Similarly, any continuous linear functional on 
	%${\rm C}^{l,m}(\Omega)$ we call a 
	%\emph{distribution of anisotropic order with compact support}, 
	%and such functionals form a vector space 
	%$\mathcal{E}'_{l,m}(\Omega)=({\rm C}^{l,m}(\Omega))'$.
\end{definition}

A linear functional $T$ on
${\rm C}^{l,m}_c(\Omega)$ is continuous (at zero) if and only if its restriction to any space ${\rm C}^{l,m}_K(\Omega)$,
where $K\subseteq\Omega$ is a compact, %(denoted also by $K\in\lK(\Omega)$ below), 
is continuous (at zero) \cite{AEM}, meaning that
$$
(\forall K\Subset\Omega) (\exists{C>0}) 
(\forall \ph\in{\rm C}^{l,m}_c(\Omega)) \quad \supp\ph\subseteq K \implies |\langle T,\ph\rangle| \leq C p^{l,m}_K(\ph)\;.
$$
If either $l$ or $m$ is infinite, we have to modify the above in the obvious way.

Since inclusions $\mathrm{C}^\infty_c(\Omega)\hookrightarrow
\mathrm{C}^{l,m}_c(\Omega)$ are continuous and dense, 
any distribution of anisotropic order is a distribution.
In fact, by the equivalent definition of continuity above, we can characterise
anisotropic distributions as distributions of order $(l,m)$, i.e.
\begin{equation*}
T\in \lD'_{l,m}(\Omega) \iff
\left\{
\begin{array}{l}
T\in \lD'(\Omega) \,, \ \hbox{and} \\
(\forall K\Subset\Omega) (\exists{C>0}) 
(\forall \ph\in{\rm C}^\infty_K(\Omega)) \quad 
|\langle T,\ph\rangle| \leq C p^{l,m}_K(\ph) \;,
\end{array}
\right.
\end{equation*}
where, again, 
if either $l$ or $m$ is infinite, we have to modify the above in the obvious way.

Let us now consider the case when $\Omega$ is an open subset of 
a differentiable manifold without boundary.
To be precise, under the notion \emph{differential manifold without boundary
	of dimension $d$}
we mean a locally Euclidean (of the fixed dimension $d$, i.e.~locally 
diffeomorphic to $\Rd$) 
second countable Hausdorff topological space on which an equivalence class of 
${\rm C}^\infty$ smooth atlases is given. 
In particular, every differentiable manifold without boundary $X$ is 
paracompact, hemicompact (i.e.~admits a countable compact exhaustion), 
separable, metrisable, and for any open cover of $X$ there exists a 
smooth partition of unity subordinate to this cover (see e.g.~\cite[16.1--4]{Dieud}
or \cite[Chapter 1]{Lee13}).

Definition \ref{defn:aniso-distr} applies also when $\Omega$ is an open subset of 
$X\times Y$, where $X$ and $Y$ are differentiable manifolds without boundary 
of dimensions $d$ and $r$, respectively.
Indeed, the construction is completely standard by means of reducing the problem (locally) 
to charts (a detailed construction can be found in \cite[Subsection 2.2]{AEM}).

Since the main advantage of using anisotropic distributions is
to give a more precise description of the order of distributions,
one needs a more precise version of the Schwartz kernel theorem
(the classical kernel theorem for manifolds without boundary can be found 
in \cite[23.9.2]{Dieud}).
In \cite[Theorem 5]{AEM} the following theorem has been proven. 

\begin{theorem}\label{thm:kernel-no-bd}
	Let $X$ and $Y$ be differential manifolds, of dimension $d$ and $r$, and $l,m\in\N_0\cup\{\infty\}$. 
	Then the following statements hold:
	
	\begin{itemize}
		\item[i)] If $K\in\mathcal{D}'_{l,m}(X\times Y)$, then for each $\varphi\in{\rm C}^l_c(X)$ the linear form $K_\varphi$, defined by
		$\psi\mapsto\langle K,\varphi\otimes\psi\rangle$, is a distribution of order not more than $m$ on $Y$.
		Furthermore, the mapping $\varphi\mapsto K_\varphi$, taking ${\rm C}^l_c(X)$ with its strict inductive limit topology to $\mathcal{D}'_m(Y)$
		with weak $\ast$ topology, is linear and continuous.
		
		\item[ii)] Let $A:{\rm C}^l_c(X)\to\mathcal{D}'_m(Y)$ be a continuous linear operator, in the pair of topologies as in (i) above.
		%	 such that for every $L\subset X$ compact set, $A:{\rm C}^l_L(X)\to{\cal D}'_m(Y)$ is continuous.
		Then there exists a unique distribution of anisotropic order $K\in\mathcal{D}'_{l,r(m+2)}(X\times Y)$ such that for 
		any $\varphi\in{\rm C}^l_c(X)$ and $\psi\in{\rm C}^{r(m+2)}_c(Y)$ one has
		$$
		\langle K,\varphi\otimes\psi\rangle = \langle K_\varphi, \psi\rangle = \langle A\varphi,\psi\rangle \;.
		$$
		\vskip-2mm
	\end{itemize}
\end{theorem}

A nice property of the second statement in the previous theorem is the fact that
we can keep the order with respect to one space unchanged.
In particular, if we have $l=0$, then the resulting kernel distribution is still 
defined on merely continuous functions in $X$.

\subsection{Anisotropic distributions on manifolds with boundary}

For the construction of one-scale H-distributions (which is the main result of the 
following section) we need to study distributions on $\kob$. 
Since $\kob$ is not a manifold without boundary, but 
can be understood as a (compact) manifold with boundary, we 
need to generalise the notion of anisotropic distributions and 
the corresponding kernel theorem of the previous subsection.
In order to avoid any confusions, let us stress that the definition of
\emph{differential manifolds with boundary} used here differs from 
the notion of differential manifolds without boundary only in the fact that we allow 
that the former is locally diffeomorphic either to $\Rd$ or to the closed half-space
$\cl\R^d_+=\{\mx=(x_1,x_2,\dots,x_d)\in\Rd : x_d\geq 0\}$
(see e.g.~\cite[Chapter 1]{Lee13}).

For simplicity of the presentation, in this subsection 
we shall consider only the case $X=\Omega$, $\Omega\subseteq\Rd$ open, 
and $Y=\kob$, while one can generalise all the conclusions to the 
case where $X$ and $Y$ are differentiable manifolds without boundary and
with boundary, respectively (some remarks in this direction are addressed at the 
end of this subsection).

Since studying distributions on compact sets is not completely standard, 
let us start by fixing our notation. 

The space of distributions on $\kob$ of order $l\in\N\cup\{\infty\}$ 
we define by 
\begin{equation*}
\lD_l'(\kob) = \bigl(\mathrm{C}^l(\kob)\bigr)' \,,
\end{equation*}
where the case $l=\infty$ we shall also denote by $\lD'(\kob)$. 
Of course, for any $l\leq m$ we have the continuous inclusion $\lD_l'(\kob)\subseteq \lD_m'(\kob)$.

Although this definition is quite natural and simple, one needs to be 
careful how to define derivatives of these objects since the integration 
by parts leaves some nontrivial boundary terms
(as test functions do not vanish on the boundary).
Such difficulties will not occur in the present manuscript, so 
we shall not elaborate on the issue in more detail.

For $\Omega\subseteq \Rd$ open, the space of anisotropic distributions 
on $\Omega\times\kob$ of 
order $(l,m)\in (\N\cup\{\infty\})^2$
is defined by
\begin{equation}\label{eq:aniso-distr-kob}
\lD_{l,m}'(\Omega\times\kob) = \bigl(\mathrm{C}^{l,m}_c(\Omega\times\kob)\bigr)'\,. 
\end{equation}
Since $\kob$ is a compact set, the compact support of test functions matters 
only with respect to the first variable $\mx\in\Omega$. 

Let us note that it is sufficient to introduce distributions on 
$\Omega\times\kvz$ since by applying the pushforward $(\lJ^{-1})_*$
we have a one-to-one correspondence with distributions on $\Omega\times\kob$.
More precisely, the mapping $(\lJ^{-1})_*$ given by
$$
\bigl\langle (\lJ^{-1})_*\nu ,\Phi\bigr\rangle 
= \bigl\langle \nu,\Phi\bigl(\cdot,\lJ^{-1}(\cdot)\bigr)\bigr\rangle 
\;, \quad \nu\in\lD_{l,m}'(\Omega\times\kvz) \,,
$$
defines an isomorphism between spaces 
$\lD_{l,m}'(\Omega\times\kvz)$ and $\lD_{l,m}'(\Omega\times\kob)$.
Here the space $\lD_{l,m}'(\Omega\times\kvz)$ is defined in the same manner as
in \eqref{eq:aniso-distr-kob}, i.e.
$\lD_{l,m}'(\Omega\times\kvz) = \bigl(\mathrm{C}^{l,m}_c(\Omega\times\kvz)\bigr)'$.

It can easily be seen that $\kvz$ is an embedded submanifold with boundary in 
the unit sphere around the origin in $\R^{d+1}$, denoted by $\mathrm{S}^d$.
Indeed, it is even a \emph{regular domain} in $\mathrm{S}^d$ (cf.~\cite[Chapter 5]{Lee13})
with a \emph{defining function} $\mathrm{S}^d\ni (\zeta_0,\meta)\mapsto \zeta_0$ 
(see \eqref{eq:cl-kvz}).
Thus, there exists a simultaneous linear extension
$$
E : \mathrm{C}^m(\kvz) \to \mathrm{C}^m(\mathrm{S}^d)
$$
which is continuous and such that $E(\psi)\rest{\kvz}=\psi$.
This can easily be proved by using the Seeley extension theorem for the half line
and the local character of manifolds with boundary, or by the above defining function. 
A nice exposition on extension theorems can be found in \cite[Chapter 2]{BB12}.

The restriction operator $R: \mathrm{C}^m(\mathrm{S}^d) \to \mathrm{C}^m(\kvz)$
defined by $R\psi = \psi\rest{\kvz}$ is linear and continuous. 
Moreover, since for $\psi\in \mathrm{C}^m(\kvz)$ we have 
$\psi = R(E\psi)$, the restriction operator is surjective. 

We shall make a small abuse of notation and use the same symbol $E$
for the extension operator from  $\Omega\times\kvz$ to $\Omega\times\mathrm{S}^d$,
which is defined on tensor products by 
$E(\varphi\otimes\psi)=\varphi\otimes E\psi$.

Let us consider a continuous 
linear operator $A:{\rm C}^l_c(\Omega)\to\mathcal{D}'_m(\kvz)$,
where the inductive limit topology is assumed on ${\rm C}^l_c(\Omega)$
and the weak $\ast$ topology on $\mathcal{D}'_m(\kvz)$.
We want to show a statement analogous to Theorem \ref{thm:kernel-no-bd}(ii).
The strategy is to reduce the problem to the setting of Theorem \ref{thm:kernel-no-bd}. 

By
\begin{equation*}
\langle \tilde A\varphi,\tilde\psi\rangle = \langle A\varphi,R\tilde{\psi}\rangle
\;, \qquad \varphi\in {\rm C}^l_c(\Omega)\,,\, 
\tilde{\psi}\in \mathrm{C}^m(\mathrm{S}^d) \,,
\end{equation*}
a continuous linear operator $\tilde{A}:{\rm C}^l_c(\Omega)
\to\mathcal{D}'_m(\mathrm{S}^d)$ is defined. 
Applying Theorem \ref{thm:kernel-no-bd}(ii) we get the distributional kernel
$\tilde{K}\in\mathcal{D}'_{l,d(m+2)}(\Omega\times \mathrm{S}^d)$
such that for any $\varphi\in{\rm C}^l_c(\Omega)$ and 
$\tilde\psi\in{\rm C}^{d(m+2)}(\mathrm{S}^d)$ it holds
$$
\langle \tilde K,\varphi\otimes\tilde \psi\rangle	
= \langle \tilde A\varphi,\tilde\psi\rangle \;.
$$

Finally, we define $K\in\mathcal{D}'_{l,d(m+2)}(\Omega\times \kvz)$
by
\begin{equation*}
\langle K,\Phi\rangle = \langle\tilde K, E\Phi\rangle \;, \quad
\Phi\in \mathrm{C}^{l,d(m+2)}(\Omega\times\kvz) \,,
\end{equation*}
and it remains to be seen that $K$ is the kernel of operator $A$. 

Using the above, for any $\varphi\in{\rm C}^l_c(X)$ 
and $\psi\in{\rm C}^{d(m+2)}(\kvz)$ we have
\begin{align*}
\langle K,\varphi\otimes\psi\rangle	
= \langle \tilde{K},\varphi\otimes E\psi\rangle
= \langle \tilde{A}\varphi,E\psi\rangle
= \langle A\varphi,R(E\psi)\rangle
= \langle A\varphi,\psi\rangle \;.
\end{align*}

Therefore, we have proved the following corollary of Theorem \ref{thm:kernel-no-bd}.
\begin{corollary}\label{cor:kernel}
	Let $\Omega\subseteq\Rd$ be open and $l,m\in\N_0\cup\{\infty\}$. 
	Furthermore, let $A:{\rm C}^l_c(\Omega)\to\mathcal{D}'_m(\kob)$ be a continuous 
	linear operator, taking ${\rm C}^l_c(\Omega)$ with its inductive limit topology 
	and $\mathcal{D}'_m(\kob)$ with weak $\ast$ topology.
	
	Then there exists a unique distribution of anisotropic order $K\in\mathcal{D}'_{l,d(m+2)}(\Omega\times \kob)$ such that for 
	any $\varphi\in{\rm C}^l_c(X)$ and $\psi\in{\rm C}^{d(m+2)}(\kob)$ one has
	$$
	\langle K,\varphi\otimes\psi\rangle	= \langle A\varphi,\psi\rangle \;.
	$$
	\vskip-2mm
\end{corollary}

\begin{remark}
	\begin{itemize}
		\item[i)] The transition from $\kob$ to $\kvz$ in the analysis above is not essential. 
		Indeed, if we had an arbitrary differential manifold $Y$ instead of $\kob$, the
		existence of a differential  manifold without boundary $\tilde Y$ in which $Y$ is embedded 
		would be ensured by the standard theory of smooth manifolds (see e.g. Theorem 6.15 and Example
		9.32 in \cite{Lee13}). Furthermore, a simultaneous linear extension from $Y$ to $\tilde Y$
		exists even in this general case (see e.g.~\cite[Section 2.2]{BB12}), which can also be 
		proven constructively by using the local character of manifolds with boundary. 
		
		\item[ii)] In \cite{RM} one can find a more general theory of distributions on 
		manifolds with corners. In order to have well defined integrals in such generality
		(i.e.~to be able to \emph{identify}\/ functions with distributions), 
		one has to consider densities either at the level of distributions or at the level of
		test functions. 
		Since $\kvz$ is a regular domain in $\mathrm{S}^d$, it is orientable, hence we have a canonical way how to integrate
		differential $d$-forms. Therefore, in our situation explicit reference to densities can be avoided.

%		Since $\kvo$ is open in $\Rd$, on $\kvz$ we have a canonical way how to integrate 
%		differential $d$-forms, thus in our situation densities can be omitted. 
		
		Moreover, in \cite[Chapter 3]{RM} three different definitions of distributions 
		on manifolds with corners are presented.
		Let us briefly mention all of them. 
		For a smooth compact manifold with corners $X$ let us denote by $\Omega X$
		a $\mathrm{C}^\infty$ line bundle over $X$ consisting of densities 
		($1$-densities). The spaces
		$$
		\bigl(\mathrm{C}^\infty_c(\Int X;\Omega X)\bigr)' \ , \quad
		\bigl(\mathrm{C}^\infty(X;\Omega X)\bigr)' \  \hbox{\ and} \quad
		\bigl(\mathrm{C}^\infty_0(X;\Omega X)\bigr)' \;
		$$
		are called \emph{distributions in the interior}, \emph{supported distributions} and 
		\emph{extendible distributions}, respectively.
		Thus, our choice corresponds to the supported distributions. 
	\end{itemize}
\end{remark}

\section{One-scale H-distributions}

\subsection{Existence and first properties}

Now we have prepared all the prerequisites required for the construction of one-scale H-distributions.

\begin{theorem}[existence of one-scale H-distributions]\label{thm:jhd_exist}
Let $\Omega\subseteq\Rd$ be open.
If $u_n\rightharpoonup 0$ in $\mathrm{L}_\mathrm{loc}^p(\Omega)$ and $(v_n)$ is bounded in 
$\mathrm{L}_\mathrm{loc}^{q}(\Omega)$ (for some $p\in(1,\infty)$ and $q\ge p'$, where $1/p+1/p'=1$), and if\/ $\omega_n\to0^+$, 
then there exist subsequences $(u_{n'})$, $(v_{n'})$, 
and a complex valued (supported) distribution
$\jhd^{(\omega_{n'})}\in\lD'_{0,\kappa}(\Omega\times\kob)$,
where $\kappa:=d(\lfloor \frac{d}{2}\rfloor+3)$, such that
for any $\ph_1, \ph_2\in\mathrm{C}_c(\Omega)$ and $\psi\in\pC\kappa\kob$, one has:
\begin{equation}\label{eq:jhd_def}
\begin{aligned}
\lim_{n'\to\infty} \int\limits_\Rd \lA_{\psi_{n'}}(\ph_1 u_{n'})(\mx) \overline{(\ph_2 v_{n'})(\mx)} 
	\,d\mx &= \lim_{n'\to\infty} \int\limits_\Rd (\ph_1 u_{n'})(\mx) \overline{\lA_{\bar\psi_{n'}}
	(\ph_2 v_{n'})(\mx)} \,d\mx \\
&= \Dupp{\jhd^{(\omega_{n'})}}{\ph_1\bar\ph_2\otimes\psi} \,,
\end{aligned}
\end{equation}
where $\psi_n=\psi(\omega_n\cdot)$.
The distribution $\jhd^{(\omega_{n'})}$ we call the\/ {\rm one-scale H-distribution 
	(with the characteristic length $(\omega_{n'})$)} associated to (sub)sequences 
$(u_{n'})$ and $(v_{n'})$.

Moreover, for $p=2$ the one-scale H-distribution above is the one-scale H-measures 
with characteristic length $(\omega_{n'})$ associated to (sub)sequences $(u_{n'})$ and $(v_{n'})$. 
\end{theorem}

%The dual product in the statement of previous theorem is well defined as we shall 
%show that $(\ph,\psi)\mapsto\langle\jhd^{(\omega_{n'})},\ph\otimes\psi\rangle$ is a
%continuous bilinear functional in the topology of $\mathrm{C}_c(\Omega)\times\pC\kappa\kob$, 
%hence, we can extend it by continuity to the whole space. Therefore, 
%in the previous statement we have in addition provided the link between that 
%extension and the corresponding limit of integrals.
%
%Here we always take a dual product to be antilinear in the first, and linear in the 
%second argument, so the sequence of integrals we can equivalently write as
%$$
%\int_\Rd \lA_{\psi_{n'}}(\ph_1 u_{n'}) \overline{\ph_2 v_{n'}} \,d\mx 
%	= \Dupp{\ph_2 v_{n'}}{\lA_{\psi_{n'}}(\ph_1 u_{n'})} \,.
%$$

\begin{proof}
Let us first observe that for $p=q=2$, after applying the Plancherel theorem, \eqref{eq:jhd_def} 
reveals the definition of one-scale H-measures (\cite[Lemma 32.6]{tar_book}, \cite[Theorem 6]{AEL}), 
and since the space $\pC\kappa\kob$ is 
dense in $\Cp{\kob}$, we get the claim.

Returning to the general $p$, we should note that for $q\ge p'$, $(\ph_2 v_n)$ is bounded in $\mathrm{L}^{p'}(\Omega)$ as well,
so we can work as if $q=p'$.

The first equality in \eqref{eq:jhd_def} stems from the fact that $\lA_{\bar\psi_{n'}}$ is the adjoint of $\lA_{\psi_{n'}}$.

%Let us denote by $(K_m)$ a sequence of compacts in $\Omega$ which exhaust $\Omega$; more precisely, 
%such that $K_m\subseteq \Int K_{m+1}$ and $\bigcup_m K_m=\Omega$.

Let us choose arbitrary test functions $\ph_1,\ph_2\in\mathrm{C}_c(\Omega)$ and 
$\psi\in\mathrm{C}^{\lfloor\frac{d}{2}\rfloor+1}(\kob)$.
Then there exists a compact set $K\subseteq\Omega$ such that the supports of both $\ph_1$ and $\ph_2$
are contained in it.
By Theorem \ref{thm:multipliers} (bearing in mind the invariance of multipliers under dilations; see Lemma \ref{lm:pom_mihlin}),
the sequence $(\lA_{\bar\psi_n}(\ph_1 u_n))$ is bounded in $\mathrm{L}^{p}(\Rd)$.
Thus, by the H\"older inequality the following sequence of integrals is bounded (in $\C$):
$$
\biggl|\int_\Rd \lA_{\psi_{n}}(\ph_1 u_{n})(\mx) \overline{(\ph_2 v_{n})(\mx)} \,d\mx\biggr| 
	\leq C \|\ph_1\|_{\mathrm{L}^\infty(K)}\|\ph_2\|_{\mathrm{L}^\infty(K)}\|\psi\|_{\pC{\lfloor\frac{d}{2}\rfloor+1}\kob} \,,
$$
where $C=C_{d,p}C_d(\sup_n\|u_n\|_{\mathrm{L}^p(K)})(\sup_n\|v_n\|_{\mathrm{L}^{p'}(K)})$, 
while $C_{d,p}$ and $C_d$ are the constants given in Theorem \ref{thm:multipliers}. 
Hence, we can pass to a subsequence converging to a limit satisfying the same bound.
Moreover, since $\mathrm{C}^{\lfloor\frac{d}{2}\rfloor+1}(\kob)$ and
$$
\{\ph\in\mathrm{C}_c(\Omega) : \supp\ph\subseteq K\}
$$
are both separable Banach spaces, while $\mathrm{C}_c(\Omega)$ can be expressed as a strict
inductive limit of the latter spaces, this passage can be done in such a 
way that the obtained subsequence is good for any choice of test functions
(see the proof of \cite[Theorem 6]{AEL}).
The corresponding subsequences we denote by $(u_{n'})$ and $(v_{n'})$.

Therefore, by
\begin{equation}\label{eq:jhd_proof_L}
L(\ph_1,\ph_2,\psi) = \lim_{n'} \dupp{\mathrm{L}^{p'}}{\overline{\ph_2 v_{n'}}}{\lA_{\psi_{n'}}(\ph_1 u_{n'})}{\mathrm{L}^p}
\end{equation}
a trilinear form on $\mathrm{C}_c(\Omega)\times\mathrm{C}_c(\Omega)\times\mathrm{C}^{\lfloor\frac{d}{2}\rfloor+1}(\kob)$ 
is defined. 
Furthermore, $L(\ph_1,\ph_2,\psi)$ depends only on the product $\ph_1\bar\ph_2$ (and not on both $\ph_1$ and 
$\ph_2$ independently). Indeed, let us denote by $\zeta_i\in\mathrm{C}_c(\Omega)$ 
an arbitrary function which is equal to 1 on the support of $\ph_i$, $i=1,2$. Then we have
$$
\begin{aligned}
\lim_{n'}\dupp{\mathrm{L}^{p'}}{\overline{\ph_2 v_{n'}}}{\lA_{\psi_{n'}}(\ph_1 u_{n'})}{\mathrm{L}^p} 
	&= \lim_{n'}\dupp{\mathrm{L}^{p'}}{\overline{\ph_2 v_{n'}}}{\ph_1\lA_{\psi_{n'}}(\zeta_1 u_{n'})}{\mathrm{L}^p} \\
&= \lim_{n'}\dupp{\mathrm{L}^{p'}}{\overline{\bar\ph_1\ph_2 v_{n'}}}{\lA_{\psi_{n'}}(\zeta_1 u_{n'})}{\mathrm{L}^p} \\
&= \lim_{n'}\dupp{\mathrm{L}^{p'}}{\overline{\zeta_1\zeta_2 v_{n'}}}{\ph_1\bar\ph_2\lA_{\psi_{n'}}(\zeta_1 u_{n'})}{\mathrm{L}^p} \\
&= \lim_{n'}\dupp{\mathrm{L}^{p'}}{\overline{\zeta_1\zeta_2 v_{n'}}}{\lA_{\psi_{n'}}(\ph_1\bar\ph_2 u_{n'})}{\mathrm{L}^p} \,,
\end{aligned}
$$
where in the first and the last equality we have used Lemma \ref{lm:commutation} and the assumption 
that $u_n\rightharpoonup 0$ in $\mathrm{L}_\mathrm{loc}^{p}(\Omega)$. 
%Therefore, we have $L(\ph_1,\ph_2,\psi)=L(\ph_1\bar\ph_2,\zeta_1\zeta_2,\psi)$.
Since $\zeta_1\zeta_2$ depends only on $\ph_1\bar\ph_2$ (we only need $\zeta_1\zeta_2=1$ on the support of
$\ph_1\bar\ph_2$), we have the conclusion. 

For $\ph\in\mathrm{C}_c(\Omega)$ and $\psi\in\pC{\lfloor\frac{d}{2}\rfloor+1}\kob$ let us define 
$B(\ph,\psi):=L(\ph,\zeta,\psi)$, where $\zeta\in\mathrm{C}_c(\Omega)$ is equal to 1 on 
the support of $\ph$. By the computations above it is easy to see that $B$ does not depend 
on the choice od $\zeta$. Thus, $B$ is a well-defined continuous bilinear form on $\mathrm{C}_c(\Omega)\times\pC{\lfloor\frac{d}{2}\rfloor+1}\kob$ 
which satisfies
\begin{equation}\label{eq:jhd_proof_B}
B(\ph_1\bar\ph_2,\psi) = L(\ph_1,\ph_2,\psi) \;, \qquad \ph_1,\ph_2\in\mathrm{C}_c(\Omega), \, \psi\in \pC{\lfloor\frac{d}{2}\rfloor+1}\kob \,.
\end{equation}

Now we can apply Corollary \ref{cor:kernel} on the bilinear form $B$, which gives us the unique $\jhd^{(\omega_{n'})}
\in\lD'_{0,\kappa}(\Omega\times\kob)$. By \eqref{eq:jhd_proof_L} and \eqref{eq:jhd_proof_B},
this distribution satisfies the required equality \eqref{eq:jhd_def}. 
\end{proof}

When it will cause no ambiguities, we shall omit explicit writting of the characteristic length, 
e.g.~we shall use only $\jhd$ instead of $\jhd^{(\omega_{n'})}$, for simplicity.

Since the role of spaces in Corollary \ref{cor:kernel} can be interchanged, 
it is straightforward that one-scale H-distributions
are contained in $\lD'_{2d,\lfloor\frac{d}{2}\rfloor+1}(\Omega\times\kob)$ as well.
However, even the intersection $\lD'_{0,\kappa}(\Omega\times\kob)\cap\lD'_{2d,\lfloor\frac{d}{2}\rfloor+1}(\Omega\times\kob)$
does not provide an optimal description of one-scale H-distributions.
Indeed, in the case $p=2$ the resulting object 
is always a Radon measure (one-scale H-measure; cf.~\cite{Tar2, AEL}).

\begin{remark}\label{rem:jhd-def-bounded-seq}
In fact, in the proof of the previous theorem we need only that at least one of sequences 
$(u_n)$, $(v_n)$ converges weakly to zero (in the corresponding space).
This was needed for application of Lemma \ref{lm:commutation}. 
Thus, the statement of Theorem \ref{thm:jhd_exist} still holds if $(u_n)$ is merely 
bounded in $\mathrm{L}_{loc}^p(\Omega)$ and 
$v_n\rightharpoonup 0$ in $\mathrm{L}_\mathrm{loc}^{p'}(\Omega)$.

If $(u_n)$ and $(v_n)$ are both only bounded (in the corresponding spaces),
then we cannot (in general) apply Lemma \ref{lm:commutation}.
Nevertheless, one can still conclude that the limit (on a subsequence) 
\eqref{eq:jhd_def} is equal to the sum of the one-scale H-distribution 
associated to (sub)sequences $(u_{n'}-u)$ and $(v_{n'}-v)$, where 
$u$ and $v$ are the weak limits of $(u_{n'})$ and $(v_{n'})$, and the
\emph{corrector} term
\begin{equation*}
	\int_{\Rd} \lA_{\psi_0\circ\mpi}(\ph_1 u)(\mx) \overline{(\ph_2 v)(\mx)}\,d\mx  \,,
\end{equation*}
where $\psi_0$ is given by \eqref{eq:psinib} and $\mpi(\mxi)=\mxi/|\mxi|$
(this can be easily justifed by using Corollary \ref{cor:psinib} and Lemma \ref{lm:technical_symbol} given below).
Note that this term cannot be incorporated within the one-scale 
H-distribution since it is in general only a trilinear object. 

However, if $\psi_0$ is a constant function, i.e.~if $\psi\in\pC{\kappa}{\kob}$ can be
extended by continuity at the origin, then the above corrector term can be written
as 
\begin{equation*}
\Dupp{u\bar v \lambda \otimes \delta_\vnul}{\ph_1\bar\ph_2\otimes\psi} \,,
\end{equation*}
where $\lambda$ is the Lebesgue measure in $\mx$ and $\delta_\vnul$ is the Dirac 
mass at $\mxi=0$. This term is obviously a Radon measure.
Therefore, if we replace the space $\pC{\kappa}{\kob}$ in the statement of the previous theorem
by its subspace $\lS(\Rd)$ (see Corollary \ref{cor:examples}(ii)), 
one might weaken the assumptions 
and assume that $(u_n)$ and $(v_n)$ are only bounded (in the corresponding spaces),
while \eqref{eq:jhd_def} still holds (with an anisotropic distribution explicitly given 
by the previous discussion). This approach we consider when studying 
semiclassical distributions (see Definition \ref{def:semiclass-distr}).
Another point of view to this generalisation might be by a direct inspection 
of the proof of the previous theorem. 
Indeed, by Remark \ref{rem:comm}(ii) for $\psi\in\lS(\Rd)$ the whole commutator $[B_\ph, \lA_{\psi_n}]$ 
converges to zero in the operator norm. Hence, the part of the proof where it is shown that 
$L$ is a bilinear object (the only place of the proof where Lemma \ref{lm:commutation} is used) 
is in this case valid if $(u_n)$ and $(v_n)$ are only bounded. 
\end{remark}

Since one-scale H-distributions are defined at the level of subsequences, for a given 
pair of sequences the corresponding one-scale H-distribution is not necessarily unique.
Thus, we shall say that $(u_n)$ and $(v_n)$ form an \emph{$(\omega_n)$-pure pair of
sequences} if the associated one-scale H-distribution with characteristic length 
$(\omega_n)$ is unique for any choice of subsequences. 
The order of sequences in the definition of one-scale H-distributions is important. 
However, the following straightforward relation holds (for the same result in the case of
H-distributions, v.~\cite[Section 4]{AEM}).

\begin{lemma}\label{lm:jhd-bar}
Let $\omega_n\to 0^+$ and let $(u_n)$ and $(v_n)$ form an $(\omega_n)$-pure pair of
sequences. Denote by $\jhd$ the corresponding one-scale H-distribution. 

Then the pair $(v_n)$ and $(u_n)$ is also $(\omega_n)$-pure, and the one-scale 
H-distribution corresponding to $(v_n)$ and $(u_n)$ is $\jhds{\bar\nu}$, where 
$\Dupp{\jhds{\bar\nu}}{\Psi}=\overline{\Dupp{\jhd}{\bar\Psi}}$.
\end{lemma}

The following property is an immediate consequence of the local nature of the 
definition of one-scale H-distributions.

\begin{corollary}\label{cor:jhd_F_support}
Let $(u_n)$ and $(v_n)$ be sequences from Theorem \ref{thm:jhd_exist}.
If there exist closed sets $F_1$ and $F_2$ in $\Omega$ such that all $u_n$ keep their support 
in $F_1$ and $v_n$ in $F_2$, then the support of any one-scale H-distribution with any 
characteristic length corresponding to (sub)sequences of $(u_n)$ and $(v_n)$ is included 
in $(F_1\cap F_2)\times \kob$.
\end{corollary}

We have already seen in Theorem \ref{thm:jhd_exist} that one-scale H-distributions 
are an extension of one-scale H-measures to the $\mathrm{L}^p-\mathrm{L}^{q}$ setting. 
Now let us show that in the $\mathrm{L}^p-\mathrm{L}^{q}$ framework the
one-scale H-distributions are also an extension of H-distributions \cite{AM}.
%Let us recall that we say that sequences $(u_n)$ and $(v_n)$ form a \emph{pure pair} 
%if the associated H-distribution is unique for all subsequences (see \cite[Section 4]{AEM}).

\begin{corollary}\label{cor:jhd-hd}
Let $\omega_n\to 0^+$ and let $(u_n)$ and $(v_n)$ form an $(\omega_n)$-pure pair of
sequences. Denote by $\jhd^{(\omega_n)}$ the corresponding one-scale H-distribution.

Then $(u_n)$ and $(v_n)$ form a pure pair of sequences, and the corresponding H-distribution
is given by $\hd=\mpi_*\jhd^{(\omega_n)}$, i.e.
\begin{equation}\label{eq:jhd-hd}
\Dupp{\hd}{\Psi} = \Dupp{\jhd^{(\omega_{n})}}{\Psi(\cdot\,,\mpi(\cdot))} \;,
\quad \Psi\in \mathrm{C}^{0,\kappa}(\Omega\times\Sdmj) \;,
\end{equation}
where $\mpi(\mxi)=\mxi/|\mxi|$.
\end{corollary}

\begin{proof}
The claim follows by Corollary \ref{cor:examples}(i) and the fact that for $\psi=\tilde\psi\circ\mpi$, 
$\tilde\psi\in\mathrm{C}^\kappa(\Sdmj)$, \eqref{eq:jhd_def} reveals the definition of H-distributions 
(see \cite[Theorem 7]{AEM}, \cite[Theorem 2.1]{AM}). 
\end{proof}

Many other properties of H-distributions are also shared with their generalisation, the one-scale H-distributions.
Let us emphasise here only the one related to compactness of weakly converging sequences.
More precisely, we have the following result (cf.~\cite[Lemma 5]{AEM}).

\begin{lemma}\label{lm:strong_compactness}
Let $u_n\rightharpoonup 0$ in $\mathrm{L}_\mathrm{loc}^p(\Omega)$, for some $p\in (1,\infty)$.
Then the following statements are equivalent:
\begin{itemize}
\item[(a)] $u_n\to 0$ (strongly) in $\mathrm{L}_\mathrm{loc}^p(\Omega)$.
\item[(b)] For every bounded sequence $(v_n)$ in $\mathrm{L}_\mathrm{loc}^{p'}(\Omega)$ and 
every $\omega_n\to 0^+$, $(u_n)$ and $(v_n)$ form an $(\omega_n)$-pure pair and the corresponding
one-scale H-distribution is zero. 
\item[(c)] For $v_n=|u_n|^{p-2}u_n$ and some $\omega_n\to 0^+$,
$(u_n)$ and $(v_n)$ form an $(\omega_n)$-pure pair and the corresponding 
one-scale H-distribution is zero. 
\end{itemize}
\end{lemma}

\begin{proof}
It is trivial to see that (b) implies (c), while the implication from (c) to (a) follows
by \eqref{eq:jhd-hd} and \cite[Lemma 5 and Remark 7]{AEM}.
Finally, if $u_n\to 0$ in $\mathrm{L}_\mathrm{loc}^p(\Omega)$, then \eqref{eq:jhd_def}
is zero since $(\mathcal{A}_{\bar \psi_n}(\ph_2 v_n))$ is bounded in $\mathrm{L}^{p'}_{loc}(\Rd)$,
thus implying (b).
\end{proof}

\begin{remark}\label{rem:jhd-star-conv}
The form of the sequence $(v_n)$ in part (c) of the previous lemma is essential.
More precisely, for a general bounded sequence $(v_n)$ in $\mathrm{L}_\mathrm{loc}^{p'}(\Omega)$
and some $\omega_n\to 0^+$ for which $(u_n)$ and $(v_n)$ form an $(\omega_n)$-pure pair 
and the corresponding one-scale H-distribution is zero we can only deduce that 
$u_n\bar{v}_n\xrightharpoonup{\;*\;} 0$ (vaguely) as (unbounded) Radon measures. 
This immediately follows by taking $\psi=1$ in \eqref{eq:jhd_def}.
\end{remark}

For vector valued sequences $(\vu_n)$ and $(\vv_n)$ we define the corresponding matrix one-scale H-distribution 
$\jhdv$ as a matrix which $(i,j)$-component, $\jhd^{ij}$, is the one-scale H-distribution associated to 
the pair of sequences $(u_n^i)$ and $(v_n^j)$.

\subsection{Example I: Concentration}\label{subsec:conc}

For $p\in(1,\infty)$, $\mz\in\Rd$, and $\eps>0$, let us define a linear operator
$\zeta_{p,\eps}$ on $\mathrm{L}^p(\Rd)$ by
\begin{equation*}
\zeta_{p,\eps} u(\mx) = \eps^{-\frac{d}{p}} u\Bigl(\frac{\mx-\mz}{\eps}\Bigr) \;.
\end{equation*}
Of course, $\zeta_{p,\eps}$ depends on $\mz$ as well, but since it will be fixed, 
we choose to simplify the notation and omit the explicit writing of label $\mz$. 
A simple change of variables shows that $\zeta_{p,\eps}$ is a linear isometry on 
$\mathrm{L}^p(\Rd)$, i.e.~$\|\zeta_{p,\eps}u\|_{\mathrm{L}^p(\Rd)}
=\|u\|_{\mathrm{L}^p(\Rd)}$.
Moreover, for any $\eps_n\to 0^+$ and $u\in\mathrm{L}^p(\Rd)$ the sequence
$(\zeta_{p,\eps_n}u)$ converges weakly to 0 in $\mathrm{L}^p(\Rd)$ 
(cf.~\cite[Section 5]{AEM}). 

For arbitrary $\eps_n\to0^+$, $\omega_n\to0^+$, $u\in\mathrm{L}^p(\Rd)$, and 
$v\in\mathrm{L}^{p'}(\Rd)$, where
$1/p+1/p'=1$, we shall construct all one-scale H-distributions with the characteristic length
$(\omega_n)$ associated to the pair of sequences $(\zeta_{p,\eps_n}u)$ and
$(\zeta_{p'\!,\eps_n}v)$. 

For given $u$, a possible (canonical) choice for $v$ can be $v=|u|^{p-2}u$. 

First we list the following two simple lemmata that we use in the construction
(cf.~\cite[lemmata 7 and 8]{AEM}).

\begin{lemma}\label{lm:concentration_conv}
Let $p\in(1,\infty)$, $\mz\in\Rd$, and $\eps_n\to 0^+$. 
For any $u\in\mathrm{L}^p(\Rd)$ and $\ph\in\mathrm{C}_c(\Rd)$ it holds
$$
\ph\zeta_{p,\eps_n}u-\ph(\mz)\zeta_{p,\eps_n}u\longrightarrow 0 \ \; \hbox{in} \ \;
	\mathrm{L}^p(\Rd) \;. 
$$
\end{lemma}

\begin{lemma}\label{lm:concentration_comm}
Let $p\in(1,\infty)$, $\mz\in\Rd$, $\eps>0$, and let $\psi\in\mathrm{L}^\infty(\Rd)$ 
be an $\mathrm{L}^p$ Fourier multiplier. 
Then
$$
\lA_\psi \zeta_{p,\eps} = \zeta_{p,\eps}\lA_{\psi(\frac{\cdot}{\eps})} \;.
$$
\end{lemma}

For arbitrary $\ph_1,\ph_2\in\mathrm{C}_c(\Rd)$ and 
$\psi\in\mathrm{C}^{\lfloor\frac{d}{2}\rfloor+1}(\kob)$ we have
\begin{equation*}
\begin{aligned}
\lim_n\int_\Rd \lA_{\psi(\omega_n\cdot)}(\ph_1\zeta_{p,\eps_n}u)(\mx)
	& \overline{\ph_2(\mx)(\zeta_{p'\!,\eps_n}v)(\mx)}\,d\mx \\
&= \ph_1(\mz)\bar{\ph}_2(\mz) \lim_n\int_\Rd \lA_{\psi(\omega_n\cdot)}(\zeta_{p,\eps_n}u)(\mx)
	\overline{(\zeta_{p'\!,\eps_n}v)(\mx)}\,d\mx \\
&= \ph_1(\mz)\bar{\ph}_2(\mz) \lim_n\int_\Rd \bigl(\zeta_{p,\eps_n}
	\lA_{\psi(\frac{\omega_n}{\eps_n}\,\cdot\,)}u\bigr)(\mx)\overline{(\zeta_{p'\!,\eps_n}v)(\mx)}\,d\mx \\
&= \ph_1(\mz)\bar{\ph}_2(\mz) \lim_n\int_\Rd
	\lA_{\psi(\frac{\omega_n}{\eps_n}\,\cdot\,)}u(\my)\overline{v(\my)}\,d\my \;,
\end{aligned}
\end{equation*}
where we have used the preceding lemmata in the first two equalities, and the change of 
variables $\my=(\mx-\mz)/\eps_n$ in the last one.  
In order to explicitly obtain the limit, the crucial point is to apply Lemma \ref{lm:technical_symbol}
given below, which we elaborate below. 

Let us assume that the sequence of positive numbers $(\frac{\omega_n}{\eps_n})$
converges in $\R\cup\{+\infty\}$ (if not, we just need to pass to a converging subsequence),
and let us denote the limit by $c$. 
Then, for any $q\in (1,\infty)$, by Theorem \ref{thm:multipliers} and Lemma \ref{lm:pom_mihlin} the sequence 
$(\psi(\frac{\omega_n}{\eps_n}\,\cdot\,))$ is a bounded sequence of $\mathrm{L}^q$ Fourier multipliers. 
Moreover, since $\psi$ is continuous on $\Rdz$ and satisfies \eqref{eq:psinib}, 
$(\psi(\frac{\omega_n}{\eps_n}\,\cdot\,))$ converges pointwise. The limit is 
$\psi_c\circ\mpi$ if $c\in\{0,\infty\}$, and $\psi(c\,\cdot)$ for $c\in (0,\infty)$.
Thus, we can apply Lemma \ref{lm:technical_symbol} (given below) to the above sequence of integrals
and the result reads:
$$
\lim_n \int_\Rd
	\lA_{\psi(\frac{\omega_n}{\eps_n}\cdot\,)}u(\mx)\overline{v(\mx)}\,d\mx
	= \left\{\begin{array}{ll}
	\Dupp{\vartheta_0}{\psi}:=\int_\Rd \lA_{\psi_0\circ\mpi}u(\mx)\overline{v(\mx)} \,d\mx \;, & \lim_n\frac{\omega_n}{\eps_n}=0 \\
	\Dupp{\vartheta_c}{\psi}:=\int_\Rd \lA_{\psi(c\,\cdot)}u(\mx)\overline{v(\mx)} \,d\mx \;, & \lim_n\frac{\omega_n}{\eps_n}=c
		\in(0,\infty) \\
	\Dupp{\vartheta_\infty}{\psi}:=\int_\Rd \lA_{\psi_\infty\circ\mpi}u(\mx)\overline{v(\mx)} \,d\mx \;, & \lim_n\frac{\omega_n}{\eps_n}
		=\infty \,,\\
	\end{array}\right.
$$
where $\psi_0$ and $\psi_\infty$ are given by \eqref{eq:psinib}.
Note that $\vartheta_\gamma$, $\gamma\in [0,+\infty]$ defined above belong (at least) 
to the space $\lD'_{\lfloor\frac{d}{2}\rfloor+1}(\kob)$
(distributions on $\kob$ of order not more than $\lfloor\frac{d}{2}\rfloor+1$).

Thus, if $\lim_n\frac{\omega_n}{\eps_n}=:c\in [0,\infty]$, then $(\zeta_{p,\eps_n}u)$ and $(\zeta_{p'\!,\eps_n}v)$
form an $(\omega_n)$-pure pair of sequences and the corresponding 
one-scale H-distribution reads
\begin{equation*}
\Dupp{\jhd}{\Psi} = \Dupp{\vartheta_c}{\Psi(\mz,\cdot\,)} \,, \quad \Psi\in\mathrm{C}^{0,\lfloor\frac{d}{2}\rfloor+1}(\Omega\times\kob) \,.
\end{equation*}
Here we can see that the one-scale H-distribution is a distribution of order $(0,\lfloor\frac{d}{2}\rfloor+1)$, 
which is slightly better than the estimate on order provided by Theorem \ref{thm:jhd_exist} (cf.~\cite{AEM}).
However, this example shows that one-scale H-distributions are not Radon measures in general (cf.~\cite[Section 5.2]{AEM}).

The fact that $\jhd$ is, with respect to $\mx$, a constant multiple of the Dirac measure centered at $\mz$ is not surprising.
Indeed, it is easy to see that after multiplying sequences $(\zeta_{p,\eps_n}u)$ and $(\zeta_{p',\eps_n}v)$ by a smooth cut-off
function around an arbitrary neighbourhood of $\mz$, the distribution remains unchanged. Thus, by Corollary \ref{cor:jhd_F_support}
the distribution is supported with respect to $\mx$ in $\{\mz\}$, implying that it is a finite sum of multiples of derivatives 
of the Dirac mass at that point. However, by Theorem \ref{thm:jhd_exist} we know that $\jhd$ is of order $0$ with respect to 
$\mx$, which leaves the only possibility that it is a constant multiple of the Dirac measure at $\mz$. 

Applying Corollary \ref{cor:jhd-hd} we get that $(\zeta_{p,\eps_n}u)$ and $(\zeta_{p'\!,\eps_n}v)$ form a pure pair
and the corresponding H-distribution $\hd$ is given by
\begin{equation*}
\Dupp{\hd}{\Psi} = \Dupp{\vartheta_1}{\Psi(\mz,\mpi(\cdot))} 
	= \int_\Rd \lA_{\Psi(\mz,\mpi(\cdot))}u(\mx)\overline{v(\mx)} \,d\mx\,, \quad \Psi\in\mathrm{C}^{0,\lfloor\frac{d}{2}\rfloor+1}(\Omega\times\Sdmj) \,,
\end{equation*}
which agrees with the result obtained in \cite[Section 5.1]{AEM}.
We would get the same if we took $c\in\{0,\infty\}$ instead of $c=1$, since $(\psi\circ\mpi)_0=(\psi\circ\mpi)_\infty=\psi$, 
for $\psi\in \mathrm{C}^{\lfloor\frac{d}{2}\rfloor+1}(\Sdmj)$. 

In this example we can see that an advantage of one-scale H-distributions over H-distributions is in the fact that 
the former contain the information on the characteristic length of sequences. 
More precisely, by tuning the characteristic length of one-scale H-distribution $(\omega_n)$ we can detect the order of 
convergence of the characteristic length of the sequence $(\eps_n)$. On the other hand, the H-distribution 
is the same for any choice of $(\eps_n)$. 

Furthermore, we can also see that the opposite implication of the one given in Corollary \ref{cor:jhd-hd} does not 
hold. Indeed, here we have that sequences  $(\zeta_{p,\eps_n}u)$ and $(\zeta_{p',\eps_n}v)$ always form 
a pure pair (i.e.~for any choice of $\eps_n\to 0^+$ associated H-distribution
is unique). 
However, if $(\frac{\omega_n}{\eps_n})$ does not converge in $\R\cup\{+\infty\}$
(i.e., if it has two or more accumulation points), those sequences are not 
$(\omega_n)$-pure.

\begin{lemma}\label{lm:technical_symbol}
	Let $p\in(1,\infty)$ and let $(\psi_n)$ be a bounded sequence of\/ $\mathrm{L}^p$
	Fourier multipliers such that it converges almost everywhere to an $\mathrm{L}^p$
	Fourier multiplier $\psi$. 
	
	Then for any $u\in\mathrm{L}^p(\Rd)$ the sequence $\bigl(\lA_{\psi_n}u\bigr)$
	converges weakly to $\lA_\psi u$ in $\mathrm{L}^p(\Rd)$.
	
	Furthermore, if the assumptions above hold for any $p\in(1,\infty)$,
	then for any  $u\in\mathrm{L}^p(\Rd)$ the sequence $\bigl(\lA_{\psi_n}u\bigr)$
	converges strongly to $\lA_\psi u$ in\/ $\mathrm{L}^p(\Rd)$, for any $p\in(1,\infty)$.
\end{lemma} 

\begin{proof}
	Let $u\in\mathrm{L}^p(\Rd)$.
	By the assumption $(\lA_{\psi_n-\psi})$ is a bounded sequence of operators
	on $\mathrm{L}^p(\Rd)$,	implying that 
	$\bigl(\lA_{\psi_n-\psi}u\bigr)$ is a bounded sequence in $\mathrm{L}^p(\Rd)$. 
	Thus, for the first claim it is sufficient to show that
	\begin{equation}\label{eq:technical}
	\bigl(\forall \ph\in\mathrm{C}^\infty_c(\Rd)\bigr) \qquad
	\dup{\mathrm{L}^{p'}}{\ph}{\lA_{\psi_n-\psi}u}{\mathrm{L}^p} \longrightarrow 0 \;.
	\end{equation}
	
	For an arbitrary $\eps>0$, let us take $u_\eps\in\mathrm{C}^\infty_c(\Rd)$
	such that $\|u-u_\eps\|_{\mathrm{L}^p(\Rd)}<\eps$.
	Since $\|\lA_{\psi_n-\psi}(u-u_\eps)\|_{\mathrm{L}^p}=O(\eps)$
	(the bound is independent of $n$),
	we are left to prove \eqref{eq:technical} when $u$ is replaced by $u_\eps$. 
	
	By the Plancherel formula it holds
	$$
	\dup{\mathrm{L}^{p'}}{\ph}{\lA_{\psi_n-\psi}u_\eps}{\mathrm{L}^p}
	= \int_\Rd \lA_{\psi_n-\psi}u_\eps(\mx) \bar{\bar{\ph}}(\mxi) \,d\mx
	= \int_\Rd \bigl(\psi_n(\mxi)-\psi(\mxi)\bigr) \widehat{u_\eps}(\mxi) \overline{\hat{\bar\ph}(\mxi)}\,d\mxi \,.
	$$
	Therefore, for a fixed $\eps$ and $\ph$, by the Lebesgue dominated convergence theorem
	we have that $\lim_n\dup{\mathrm{L}^{p'}}{\ph}{\lA_{\psi_n-\psi}u_\eps}{\mathrm{L}^p}=0$.
	
	Assume now that the stronger condition holds, i.e.~that for any $p\in (1,\infty)$
	the sequence 
	$(\psi_n)$ is a bounded sequence of $\mathrm{L}^p$ Fourier multipliers
	and its pointwise limit $\psi$ is an $\mathrm{L}^p$ Fourier multiplier
	for any $p\in (1,\infty)$, as well.
	Let $\theta\in (0,1)$ and $r\in (1,\infty)$ be such that 
	$1/p = \theta/2 + (1-\theta)/r$ and let $u_\eps$ be as above.
	By the classical interpolation inequality we have
	\begin{align*}
	\|\lA_{\psi_{n}-\psi}u\|_{\mathrm{L}^p} &\leq 
		\|\lA_{\psi_{n}-\psi}(u-u_\eps)\|_{\mathrm{L}^p} + 
		\|\lA_{\psi_{n}-\psi}u_\eps\|_{\mathrm{L}^p} \\
	&\leq O(\eps) + \|\lA_{\psi_{n}-\psi}u_\eps\|_{\mathrm{L}^2}^\theta
		\|\lA_{\psi_{n}-\psi}u_\eps\|_{\mathrm{L}^r}^{1-\theta} \;.
	\end{align*}
	Similarly as above, by the Plancherel formula and the Lebesgue 
	dominated convergence theorem we have (for fixed $\eps$)
	$\lim_n \|\lA_{\psi_{n}-\psi}u_\eps\|_{\mathrm{L}^2} =0$,
	while $(\|\lA_{\psi_{n}-\psi}u_\eps\|_{\mathrm{L}^r})_n$ is bounded 
	since $(\psi_n-\psi)$ is a bounded sequence of $\mathrm{L}^r$
	Fourier multipliers, and this concludes the proof.
\end{proof}

\begin{remark}
It can easily be seen from the above proof that the assumption of the second part 
of the previous lemma can be weakened, while still preserving the strong 
convergence in the starting $\mathrm{L}^p(\Rd)$ space. 
Namely, if $(\psi_n)$ is a bounded
sequence of $\mathrm{L}^q$ Fourier multipliers, for $q\in\{p,r\}$, such that 
it converges almost everywhere to an $\mathrm{L}^q$ 
Fourier multiplier $\psi$ and there exists $\theta \in (0,1)$ such that 
$1/p = \theta/2 + (1-\theta)/r$, then for any $u\in\mathrm{L}^p(\Rd)$ 
the sequence $\bigl(\lA_{\psi_n}u\bigr)$
converges strongly to $\lA_\psi u$ in $\mathrm{L}^p(\Rd)$.
\end{remark}

\subsection{Example II: Oscillation}\label{subsec:oscill}

Let $p\in(1,\infty)$, $\mathsf{k}\in\Rdz$, $\eps_n\to 0^+$,
and $u\in\mathrm{L}^p_\mathrm{loc}(\Rd)$. 
We define two sequences 
$$
u_n(\mx)=u(\mx)e^{2\pi i \frac{\mx}{\eps_n}\cdot\mathsf{k}} \qquad\hbox{and}\qquad
	v_n(\mx)=|u(\mx)|^{p-2}u(\mx)e^{2\pi i \frac{\mx}{\eps_n}\cdot\mathsf{k}}\,,
$$
which are weakly converging to 0 in $\mathrm{L}^p_\mathrm{loc}(\Rd)$
and $\mathrm{L}^{p'}_\mathrm{loc}(\Rd)$, respectively. 

For an arbitrary $\omega_n\to0^+$ let us derive the formula for the one-scale H-distribution
with characteristic length $(\omega_n)$ associated to the pair of sequences
$(u_n)$, $(v_n)$. 

Take $\ph_1,\ph_2\in\mathrm{C}_c(\Rd)$ and $\psi\in\mathrm{C}^{\lfloor\frac{d}{2}\rfloor+1}(\kob)$.
Using the following property of Fourier multipliers
\begin{equation}
\label{eq:prFtr}
\lA_\psi(e^{2\pi i \,.\,\cdot\mathsf{k}}f) = e^{2\pi i \,.\,\cdot\mathsf{k}}\lA_{\psi(\,\cdot\,+\mathsf{k})}f
\end{equation}
we get
$$
\int_\Rd \lA_{\psi(\omega_n\cdot)}(\ph_1 u_n)(\mx) \overline{\ph_2(\mx) v_n(\mx)}\,d\mx
	= \int_\Rd \lA_{\psi(\omega_n\cdot\,+\frac{\omega_n}{\eps_n}\mathsf{k})}(\ph_1 u)(\mx)
	\overline{\ph_2(\mx) |u(\mx)|^{p-2}u(\mx)}\,d\mx \;.
$$

For the moment assume that the sequence of positive numbers $(\frac{\omega_n}{\eps_n})$
converges in $\R\cup\{+\infty\}$. Then by Theorem \ref{thm:multipliers}, Corollary \ref{cor:psinib},
Lemma \ref{lm:pom_mihlin}
and the fact that Fourier multipliers are invariant under translations 
(see \cite[Proposition 2.5.14]{Grafakos}) we can apply Lemma \ref{lm:technical_symbol} and
pass to the limit in the expression above. 
More precisely, by Corollary \ref{cor:psinib} for fixed $\mxi\in\Rdz$ we have
$$
\psi\bigl(\omega_n\mxi+\frac{\omega_n}{\eps_n}\mathsf{k}\bigr) \longrightarrow
	\left\{\begin{array}{ll}
	\psi_0\bigl(\frac{\mathsf{k}}{|\mathsf{k}|}\bigr) \;, & \lim_n\frac{\omega_n}{\eps_n}=0 \\
	\psi(c\mathsf{k})\;, & \lim_n\frac{\omega_n}{\eps_n}=c\in(0,\infty) \\
	\psi_\infty\bigl(\frac{\mathsf{k}}{|\mathsf{k}|}\bigr) \;, & \lim_n\frac{\omega_n}{\eps_n}=\infty
	\end{array}\right. \;.
$$
Therefore, we can conclude that the unique one-scale H-distribution is given by
$$
\jhd^{(\omega_n)} = |u|^p\lambda\otimes
	\left\{\begin{array}{ll}
	\delta_{\vnul^\frac{\mathsf{k}}{|\mathsf{k}|}} \;, & \lim_n\frac{\omega_n}{\eps_n}=0 \\
	\delta_{c\mathsf{k}} \;, & \lim_n\frac{\omega_n}{\eps_n}=c\in(0,\infty) \\
	\delta_{\infty^\frac{\mathsf{k}}{|\mathsf{k}|}} \;, & \lim_n\frac{\omega_n}{\eps_n}=\infty \,,
	\end{array}\right.
$$
where by $\lambda$ we denoted the Lebesgue measure on $\Rd$. 
Thus, in this example the one-scale H-distributions are in fact Radon measures
(distributions of order zero), which is  \emph{not} the case for all one-scale H-distributions
(see the previous example). 
Moreover, their form is the same for any $p\in(1,\infty)$ 
(cf.~\cite[Example 4]{AEL} for the case $p=2$).

In general, the sequence $(\frac{\omega_n}{\eps_n})$ can have more than one accumulation point (in $\R\cup\{+\infty\}$). 
Then the pair $(u_n)$, $(v_n)$ is not $(\omega_n)$-pure, and the corresponding one-scale H-distribution
is of one of the above three forms.

However, by Corollary \ref{cor:jhd-hd} we have that the pair $(u_n)$, $(v_n)$ is pure, and 
the corresponding H-distribution $\hd$ reads
$$
\hd = |u|^p\lambda\otimes \delta_{\frac{\mathsf{k}}{|\mathsf{k}|}} \,,
$$
so it is a Radon measure as well.

\subsection{Semiclassical pseudodifferential operators and the Wigner transform}

In the $\mathrm{L}^2$-setting ($p=2$) one gets semiclassical measures when 
the space of test functions in the dual space $\mathrm{C}^\kappa(\kob)$
is replaced by $\lS(\Rd)$ \cite[Corollary 6]{AEL}.
More precisely, semiclassical measures are restrictions of one-scale 
H-measures to $\Omega\times\Rd$. 
This motivates the introduction of the following generalisation of semiclassical measures
(for simplicity, we take the second sequence only in ${\mathrm L}^{p'}$ based spaces, and
we do not consider a minor generalisation obtained by taking $q\ge p'$ here).

\begin{definition}\label{def:semiclass-distr}
	For $\Omega\subseteq\Rd$ open and $p\in(1,\infty)$ let $(u_n)$ and $(v_n)$ 
	be sequences bounded in $\mathrm{L}^p_\mathrm{loc}(\Omega)$ and 
	$\mathrm{L}^{p'}_\mathrm{loc}(\Omega)$, respectively, and let $\omega_n\to0^+$.
	
	A distribution $\pkd^{(\omega_{n'})}\in\lD_{0,\kappa}'(\Omega\times\Rd)$, 
	where $\kappa:=d(\lfloor\frac{d}{2}\rfloor+3)$, 
	we call the \emph{semiclassical distribution 
		(with characteristic length $(\omega_{n'})$)} associated to (sub)sequences 
	$(u_{n'})$ and $(v_{n'})$ if 
	for any $\ph_1,\ph_2\in\mathrm{C}_c(\Omega)$
	and $\psi\in\lS(\Rd)$ we have
	\begin{equation*}
	\lim_{n'\to\infty} \int\limits_\Rd \lA_{\psi_{n'}}(\ph_1 u_{n'})(\mx) 
	\overline{(\ph_2 v_{n'})(\mx)} \,d\mx
	= \Dupp{\pkd^{(\omega_{n'})}}{\ph_1\bar\ph_2\otimes\psi} \,,
	\end{equation*}
	where $\psi_n := \psi(\omega_n\cdot)$.
	
	When it will cause no ambiguities, we shall omit explicit writting of the characteristic length, 
	e.g.~we shall use only $\pkd$ instead of $\pkd^{(\omega_{n'})}$, for simplicity.
\end{definition}

The existence of subsequences referred to above follows from Theorem \ref{thm:jhd_exist}
and Remark \ref{rem:jhd-def-bounded-seq}.

It should be mentioned that the semiclassical measures were first introduced by G\'erard \cite{GerMsc}
using semiclassical pseudodifferential operators, while an alternative 
construction was obtained by Lions and Paul \cite{LP93} in terms of the Wigner 
transform. Moreover, they are presumably better known
nowadays as Wigner measures, since the Wigner transform is a very popular 
tool in mathematical physics. 
However, both these approaches were confined to the $\mathrm{L}^2$ setting, 
hence here we present their generalisation to the $\mathrm{L}^p-\mathrm{L}^{p'}$
framework. 
Moreover, we shall see that the resulting object indeed coincides with the 
above introduced semiclassical distributions, thus justifying the name. 

We start by recalling some well known results form the theory of 
(semiclassical) pseudodifferential operators for which we mainly 
follow \cite{Zworski} (see also \cite{Zhang}).

For $a\in\lS(\Rd\times\Rd)$, $\omega\in (0,1]$ and $t\in[0,1]$ 
we define the semiclassical pseudodifferential operator in 
$t$-quantisation to be the operator $\operatorname{Op}_t^{(\omega)}(a)$
acting on $u\in\lS(\Rd)$ by formula
\begin{equation*}
\begin{aligned}
\operatorname{Op}_t^{(\omega)}(a)u(\mx)
	&= \int_\Rd \int_\Rd e^{2\pi i (\mx-\my)\cdot\mxi} 
	a(t\mx+(1-t)\my,\omega\mxi)u(\my)\,d\my d\mxi\\
&= \omega^{-d}\int_\Rd \int_\Rd e^{\frac{2\pi i}{\omega} (\mx-\my)\cdot\meta} 
	a(t\mx+(1-t)\my,\meta)u(\my)\,d\my d\meta \,,
\end{aligned}
\end{equation*}
where the second equality is due to the change of variables $\meta=\omega\mxi$. 
Let us note that in this paper we use a different definition of the Fourier 
transform than in \cite{Zworski}, and consequentially the above definition 
slightly differs as well.
For $t=\frac{1}{2}$ we get the Weyl quantisation and for $t=1$ 
the standard (Kohn-Nirenberg) quantisation, while the common notation for these cases is 
\begin{equation*}
a^\mathrm{w}(\mx,\omega D) = \operatorname{Op}_\frac{1}{2}^{(\omega)} \qquad
	\hbox{and} \qquad
	a(\mx,\omega D) = \operatorname{Op}_1^{(\omega)} \;.
\end{equation*}
If the symbol $a$ does not depend on $\mx$, i.e.~$a=a(\mxi)$, 
then the above operators are just Fourier multipliers.
More precisely, for any $t\in [0,1]$ we have 
$\operatorname{Op}_t^{(\omega)}(a)=\lA_{a(\omega\,\cdot)}$.
On the other hand, if $a$ does not depend on $\mxi$, i.e.~$a=a(\mx)$,
then the above operator for $t=1$ simplifies to the operator of 
multiplication by $a$, i.e.~$\operatorname{Op}_1^{(\omega)}(a)u(\mx)
=a(\mx)u(\mx)$.

A very useful conversion formula between different quantisations holds:
if 
$$
\operatorname{Op}_t^{(\omega)}(a_t) = \operatorname{Op}_s^{(\omega)}(a_s)
$$
for $t,s\in [0,1]$ and $a_t, a_s\in\lS(\Rd\times\Rd)$, then
$$
a_t(\mx,\mxi) = e^{\frac{2\pi i}{\omega} (t-s)(\omega D_\mx)\cdot 
	(\omega D_\mxi)} a_s(\mx,\mxi) \,. 
$$
By \cite[theorems 3.17(i) and 4.8(ii)]{Zworski},
from the above one can deduce that there exists 
$b^{(\omega)}\in\lS(\Rd\times\Rd)$,
which is bounded with respect to $\omega$, such that 
\begin{equation}\label{eq:transition-quant}
a_t = a_s + \omega b^{(\omega)} \;.
\end{equation}

For any $a\in\lS(\Rd\times\Rd)$ and $t\in [0,1]$, 
$\operatorname{Op}_t^{(\omega)}(a)$ can be extended 
to a continuous linear operator from $\lS'(\Rd)$ to $\lS'(\Rd)$. 
Here we are particularly interested in their restrictions to $\mathrm{L}^p$ spaces, $p\in(1,\infty)$.
More precisely, for any $p\in (1,\infty)$ 
there exists a $C>0$ and multiindices $\malpha,\mbeta\in \N_0^{2d}$, 
depending only on $p$ and $d$,
such that for any $a\in\lS(\Rd\times\Rd)$, $t\in [0,1]$ 
and $\omega\in(0,1]$ it holds
\begin{equation}\label{eq:semiclass-bdd-lp}
\|\operatorname{Op}_t^{(\omega)}(a)\|_{\lL(\mathrm{L}^p(\Rd))} 
	\leq C |a|_{\malpha,\mbeta}  + O(\omega) \;,
\end{equation}
as $\omega\to 0^+$, where $|a|_{\malpha,\mbeta}=\sup_{\R^{2n}}|\mx^\malpha\partial^\mbeta a|$
are the standard seminorms used in the definition of the Schwartz space $\lS(\Rd\times\Rd)$.

A very simple proof of this statement can be found in 
\cite[Lemma 2.2]{KTZ07} for the standard quantisation, which can be used, 
together with the above conversion formula, to 
prove the result for general $t$. 
Moreover, by \eqref{eq:transition-quant} we have
\begin{equation}\label{eq:semiclass-diff}
\lim_{\omega\to 0^+} \|\operatorname{Op}_t^{(\omega)}(a) 
	- \operatorname{Op}_s^{(\omega)}(a)\|_{\lL(\mathrm{L}^p(\Rd))}=0 \;.
\end{equation}

The commutator of two semiclassical pseudodifferential 
operators (in the same quantisation) vanishes as $\omega$ approaches zero, 
i.e.~for any $a,b\in\lS(\Rd\times\Rd)$ it holds
\begin{equation}\label{eq:semiclass-comp}
\lim_{\omega\to 0^+} \bigl\|[\operatorname{Op}_t^{(\omega)}(a),
	\operatorname{Op}_t^{(\omega)}(b)]\bigr\|_{\lL(\mathrm{L}^p(\Rd))}=0 \;.
\end{equation}
This is due to the fact that the symbol of composition 
$\operatorname{Op}_t^{(\omega)}(a)\operatorname{Op}_t^{(\omega)}(b)$
is equal to $ab+\omega c^{(\omega)}$, where 
$c^{(\omega)}\in\lS(\Rd\times\Rd)$ is bounded with respect to $\omega$.

Now we use the approach of semiclassical pseudodifferential operators 
in the study of weakly converging sequences in the $\mathrm{L}^p-\mathrm{L}^{p'}$
setting, $1/p+1/p'=1$ (as it has already been done in Definition \ref{def:semiclass-distr}).

For simplicity of the exposition, in the remaining part of this subsection, we shall consider 
a more specific situation than 
in Definition \ref{def:semiclass-distr}:
let $(u_n)$ and $(v_n)$ be bounded sequences in $\mathrm{L}^p(\Rd)$ 
and $\mathrm{L}^{p'}(\Rd)$, respectively, and let $\omega_n\to0^+$.

Take any $t\in [0,1]$.
First we consider an arbitrary $a$ from a countable dense subset $\mathcal{G}$ 
of $\lS(\Rd\times\Rd)$ (the separability of the Schwartz space is a well-known result). 
By \eqref{eq:semiclass-bdd-lp} the sequence 
$\dup{\mathrm{L}^{p}}{\operatorname{Op}_t^{(\omega_n)}(a)u_n}
{\bar v_n}{\mathrm{L}^{p'}}$ is bounded, thus admitting a converging 
subsequence, which by the diagonal argument we can take to be the same 
for any $a\in\mathcal{G}$. 
For simplicity, we do not relabel the subsequence. 
At this moment we have that 
$\lim_n \dup{\mathrm{L}^{p}}{\operatorname{Op}_t^{(\omega_n)}(\cdot)u_n}
{\bar v_n}{\mathrm{L}^{p'}}$ is a linear functional on $\mathcal{G}$.
Again, by using \eqref{eq:semiclass-bdd-lp} we can see that it is in fact 
bounded on $\lS(\Rd\times\Rd)$. 
Thus, $\tilde\nu_{sc}:=\lim_n \dup{\mathrm{L}^{p}}{\operatorname{Op}_t^{(\omega_n)}(\cdot)u_n}
{\bar v_n}{\mathrm{L}^{p'}} \in \lS'(\Rd\times\Rd)$.
It is important to notice that by \eqref{eq:semiclass-diff}
the distribution $\tilde\nu_{sc}$ does not depend on $t$.

We shall prove that $\tilde\nu_{sc}$ coincides with the distribution given 
in Definition \ref{def:semiclass-distr} (for given sequences $(u_n)$, $(v_n)$ and $(\omega_n)$).

Let us assume that defining subsequences for distributions 
$\pkd$ and $\tilde{\nu}_{sc}$ are the same (otherwise we first pass to
subsequences in order to define $\pkd$ and then we apply the above 
procedure for $\tilde{\nu}_{sc}$ on that pair of subsequences). 
Consider those subsequences and denote them by $(u_n)$, $(v_n)$ and $(\omega_n)$.
For any $\varphi_1,\varphi_2\in\mathrm{C}^\infty_c(\Rd)$ and $\psi\in\lS(\Rd)$
we have
\begin{align*}
\Dupp{\pkd}{\ph_1\bar\ph_2\otimes\psi} &=
\lim_{n\to\infty} \int\limits_\Rd \lA_{\psi_{n}}(\ph_1 u_{n})(\mx) 
	\overline{(\ph_2 v_{n})(\mx)} \,d\mx\\ 
&= \lim_{n\to\infty}\int \operatorname{Op}_1^{(\omega_n)}(\bar\varphi_2)\operatorname{Op}_1^{(\omega_n)}
	(\psi)\operatorname{Op}_1^{(\omega_n)}
	(\varphi_1)(u_n)(\mx)\overline{v_n(\mx)} \,d\mx\\
&=\lim_{n\to\infty} \int \operatorname{Op}_1^{(\omega_n)}(\varphi_1\bar{\varphi}_2\otimes\psi)(u_n)(\mx)
\overline{v_n(\mx)} \,d\mx \\
&=\Dupp{\tilde\nu_{sc}}{\ph_1\bar\ph_2\otimes\psi} \;,
\end{align*}
where in the third equality we have used \eqref{eq:semiclass-comp}.
Hence, since tensor products are dense in $\lS(\Rd\times\Rd)$, 
we get $\pkd=\tilde\nu_{sc}$.

Another approach in studying microlocal defect objects associated to weakly converging sequences is to 
consider the Wigner transform, which in $t$-quantisation we define by
\begin{equation*}
W_t^{(\omega)}(u,v)(\mx,\mxi) = \int_\Rd e^{-2\pi i \my\cdot\mxi}
	u(x+\omega t\my)\overline{v(\mx-\omega (1-t)\my)} \,d\my \;,
\end{equation*}
where we take $u\in \mathrm{L}^p(\Rd)$ and $v\in\mathrm{L}^{p'}(\Rd)$, 
$p\in(1,\infty)$. The case $t=\frac{1}{2}$ is the most common, and often
only that case is considered in the literature.
For any $a\in\lS(\Rd\times\Rd)$ we have
\begin{align*}
\dupp{\lS'}{W_t^{(\omega)}(u,v)}{a}{\lS} &= \iiint_{\R^{3d}}
e^{-2\pi i \my\cdot\mxi}
u(x+\omega t\my)\overline{v(\mx-\omega (1-t)\my)} a(\mx,\mxi)\,d\my d\mx d\mxi \\
&= \iiint_{\R^{3d}} e^{2\pi i(\mx'-\my')\cdot\mxi'} u(\my')
\overline{v(\mx')} a\bigl(t\mx'+(1-t)\my',\omega\mxi'\bigr) 
\, d\my' d\mx' d\mxi'\\
&= \dupp{\mathrm{L}^{p}}{\operatorname{Op}_t^{(\omega)}(a)u}
	{\bar v}{\mathrm{L}^{p'}} \,,
\end{align*}
where in the second equality we have applied the change of variables:
$\mx'=\mx-\omega(1-t)\my$, $\my'=\mx+\omega t\my$, $\mxi'=\frac{1}{\omega}\mxi$,
having the Jacobian equal to 1. 
Thus, the Wigner transform in $t$-quantisation is connected to the
$t$-quantisation of semiclassical pseudodifferential operators.

Therefore, for the above choice of sequences $(u_n)$, $(v_n)$ and $(\omega_n)$,
we have obtained
\begin{equation*}
\begin{aligned}
\Dupp{\pkd}{a}
	&= \lim_n \dupp{\mathrm{L}^{p}}{\operatorname{Op}_t^{(\omega_n)}(a)u_n}
	{\bar v_n}{\mathrm{L}^{p'}} \\
&= \lim_n \dupp{\lS'}{W_t^{(\omega_n)}(u_n,v_n)}{a}{\lS} \,,
\end{aligned}
\end{equation*}
for any $a\in\lS(\Rd\times\Rd)$, showing that all the above described approaches 
lead to the same object, as it was the case in the $\mathrm{L}^2$ setting. 

A short note on the definition of semiclassical distributions in terms
of the Wigner transform can be found in \cite{AEP}.

\begin{remark}
In this subsection the study of the limit of $W^{(\omega_n)}(u_n,v_n)$
was conducted using the relation between the Wigner transform and 
the semiclassical pseudodifferential operators. 
Historically \cite{LP93}, in the $\mathrm{L}^2$ case the study was 
made directly by
using the Husimi transform, 
which has some benefits of its own (see also \cite{Zhang}). 
\end{remark}

\section{Localisation principle for one-scale H-distributions}

We have seen that one-scale H-distributions are a suitable tool for studying strong compactness of weakly 
converging sequences, where sequences determining zero one-scale H-distribution have a desirable property 
(see Lemma \ref{lm:strong_compactness}).
Our goal is to extract as much information as possible about one-scale H-distributions
from the observed sequences. 

In this section we study weakly converging sequences which satisfy some linear partial
differential relations with variable coefficients. More precisely, we
consider a sequence of systems of partial differential equations in divergence form of order $m\in\N$ 
with a characteristic length $(\varepsilon_n)$
\begin{equation}\label{eq:loc_princ_eq}
\sum_{|\malpha|\leq m} \eps_n^{|\malpha|} \partial_\malpha \bigl(\mA^\malpha_n \vu_n\bigr) = \vf_n \,,
\end{equation}
where $\mA^\malpha_n$ are given matrices of (variable) coefficients. 
Since we assume that a sequence of solutions $(\vu_n)$ to \eqref{eq:loc_princ_eq} is given, the relation above could serve 
as a definition of the right hand side $\vf_n$.
Moreover, as it will be seen below, for the sequence $(\vf_n)$ it is only important that it satisfies certain 
compactness properties.

Equations of the form \eqref{eq:loc_princ_eq} often appear when modelling
highly oscillating problems (e.g.~the passage from micro- to meso- or macro-scale,
the semiclassical limit or some scaling suitable for a particular physical model).
A more precise list of related problems with selected references can be found in the Introduction. 

In \cite[Section 4]{AEL} this problem was treated in the $\mathrm{L}^2$ framework, 
while here we present a generalisation to the $\mathrm{L}^p$ setting, 
for $p\in(1,\infty)$.

\subsection{Right hand side}

Let $s\in\R$, $p\in(1,\infty)$ and $d,r\in\N$ be given. 
The Sobolev space $\W{s}{p}{\Rd;\Cr}$ we define by
\begin{equation*}
\W{s}{p}{\Rd;\Cr} := \Bigl\{ \vf\in \lS'(\Rd;\Cr) : \lA_{(1+|\mxi|^2)^{\frac s2}}\vf \in \mathrm{L}^p(\Rd;\Cr)\Bigr\} \,.
\end{equation*}
It is a vector space, and equipped with the norm $\|\vf\|_{\W{s}{p}{\Rd;\Cr}} := \|\lA_{(1+|\mxi|^2)^{\frac s2}}f\|_{\mathrm{L}^p(\Rd;\Cr)}$ it becomes a Banach 
space. 
This definition via Fourier multipliers coincides with the standard one in terms of weak derivatives 
(cf.~\cite[Section 6.2]{GrafakosModern}).

Moreover, in both definitions, of the space and of the norm,  the symbol $\mxi\mapsto (1+|\mxi|^2)^{\frac s2}$ can be replaced 
by $\mxi\mapsto(1+|\mxi|^{|s|})^{\operatorname{sign}s}$, which is the form we shall use in the rest of the paper. 
One just has to be aware that the latter symbol is 
not smooth around the origin (cf.~\cite[Section 10.1]{HormanderII}).
Another way to see this equivalence is by applying Corollary \ref{cor:examples}(iv) and Theorem \ref{thm:multipliers}.

The local Sobolev spaces are naturally defined in the following manner:
\begin{equation*}
\Wloc{s}{p}{\Omega;\Cr} := \Bigl\{ \vf\in \lD'(\Omega;\Cr) : (\forall\ph\in\mathrm{C}^\infty_c(\Omega)) \ \ph \vf \in \W{s}{p}{\Rd;\Cr}\Bigr\} \,,
\end{equation*}
where we assume identification of functions on $\Omega$ with their extension by 
zero to the whole $\Rd$.
When equipped with the weakest topology in which every mapping $\vf\mapsto \ph\vf$ 
from $\Wloc{s}{p}{\Omega;\Cr}$ to $\W{s}{p}{\Rd;\Cr}$ is continuous, $\Wloc{s}{p}{\Omega;\Cr}$ becomes a Fr\'echet space
(cf.~\cite[Chapter 31]{Treves}). In fact, the seminorms are of the form $\|\ph\vf\|_{\W{s}{p}{\Rd;\Cr}}$,
and it is enough to take $\ph$ only from a suitable countable family \cite{ABu}.

Let us return to \eqref{eq:loc_princ_eq}. Since on the left hand side we have a differential operator of order $m$, for $u\in{\rm L}^p$
the best we can expect in general for
the right hand side is that it belongs to $\mathrm{W}^{-m,p}$. However, instead of the standard strong topology, 
it is natural to include the characteristic length $(\eps_n)$ in the compactness condition for $(\vf_n)$. 
Namely, for a sequence of positive numbers $(\eps_n)$ we say that a sequence $(\vf_n)$ in $\W{-m}{p}{\Rd;\Cr}$ 
satisfies the \emph{$(\eps_n)$-compactness} condition if the following strong convergence holds:
\begin{equation}\label{ConvRhsLp}
\lA_{{1\over 1+|\eps_n\mxi|^{m}}}\vf_n \longrightarrow 
\vnul \quad \hbox{in} \quad \mathrm{L}^p(\Rd;\Cr) \,.
\end{equation}
If all $\vf_n$ belong only to $\Wloc{-m}{p}{\Omega;\Cr}$, then the condition above is naturally replaced by its local version:
\begin{equation}\label{ConvRhsLocLp}
(\forall \ph\in\mathrm{C}_c^\infty(\Omega)) \qquad 
	\lA_{\frac{1}{1+|\eps_{n}\mxi|^m}}(\ph\vf_n) 
	\longrightarrow \vnul \quad \hbox{in} \quad \mathrm{L}^p(\Rd;\Cr) \,,
\end{equation}
which we baptise the \emph{$(\eps_n)$-local compactness condition}.
In both conditions \eqref{ConvRhsLp} and \eqref{ConvRhsLocLp}
we can replace the symbol $\mxi\mapsto (1+|\eps_n\mxi|^m)^{-1}$ by
$\mxi\mapsto (1+|\eps_n\mxi|^2)^{-m/2}$, similarly as we did in the already discussed 
situation regarding the Sobolev norms. 
This can be justified by means of Corollary \ref{cor:examples}(iv), 
Theorem \ref{thm:multipliers} and Lemma \ref{lm:pom_mihlin}.

\begin{remark}
In the case $p=2$ Plancherel's theorem holds, so the conditions 
\eqref{ConvRhsLp} and \eqref{ConvRhsLocLp} can be expressed merely in 
terms of the Fourier transform. 
This was precisely the approach used in 
\cite[Chapter 32]{tar_book} and \cite{AEL}
(see conditions (32.25) and (4.3) in the aforementioned references, respectively).
\end{remark}

It is evident that for a stationary sequence $\eps_n=\eps>0$ the condition \eqref{ConvRhsLp}
and \eqref{ConvRhsLocLp} are equivalent to the strong convergences to zero of $(\vf_n)$ 
in $\W{-m}{p}{\Rd;\Cr}$ and $\Wloc{-m}{p}{\Omega;\Cr}$, respectively. 

A more precise description of $(\eps_n)$-compactness condition \eqref{ConvRhsLp} is given in the following lemma.

\begin{lemma}\label{lm:rhs-omegan}
Let $m\in\N$ and $p\in (1,\infty)$, and let $(\eps_n)$ and $(\omega_n)$ be sequences 
of positive real numbers.
\begin{itemize}
\item[i)] If\/ $\liminf\limits_{n\to\infty} \frac{\omega_n}{\eps_n}>0$, then the $(\eps_n)$-compactness
	implies $(\omega_n)$-compactness.
\item[ii)] If\/ $\limsup\limits_{n\to\infty} \frac{\omega_n}{\eps_n} <\infty$, then the $(\omega_n)$-compactness
	implies $(\eps_n)$-compactness.
\end{itemize}
\end{lemma}

\begin{proof}
Let us notice first that it is sufficient to prove only one statement (say (i)),
since by interchanging the roles of $(\eps_n)$ and $(\omega_n)$ 
we infer the second statement from it. 

Let $\liminf\limits_{n\to\infty} \frac{\omega_n}{\eps_n}>0$ and assume that for a sequence 
$(\vf_n)$ in $\W{-m}{p}{\Rd;\Cr}$ the $(\eps_n)$-compactness holds.
We set $c_0= \liminf\limits_{n\to\infty} \frac{\omega_n}{\eps_n}$
if $\liminf\limits_{n\to\infty} \frac{\omega_n}{\eps_n}<\infty$, 
and $c_0=1$ in the case 
$\liminf\limits_{n\to\infty} \frac{\omega_n}{\eps_n}=\infty$.
Then there exits $n_0\in\N$ such that for any 
$n\geq n_0$ one has $\omega_n\geq \frac{c_0}{2}\eps_n$.

By using 
$$
\frac{1}{1+|\omega_n\mxi|^m} = \frac{1+|\eps_n\mxi|^m}{1+|\omega_n\mxi|^m}\frac{1}{1+|\eps_n\mxi|^m} \,,
$$
and the composition rule for Fourier multipliers  
$\lA_{\psi_1\psi_2} = \lA_{\psi_1}\lA_{\psi_2}$,
it is sufficient to prove that $\bigl(\frac{1+|\eps_n\mxi|^m}{1+|\omega_n\mxi|^m}\bigr)_n$
is a bounded sequence of $\mathrm{L}^p$ Fourier multipliers.
Indeed, we have
\begin{align*}
\frac{1+|\eps_n\mxi|^m}{1+|\omega_n\mxi|^m} 
	= \Bigl(\frac{\eps_n}{\omega_n}\Bigr)^m
	+ \frac{1-\bigl(\frac{\eps_n}{\omega_n}\bigr)^m}{1+|\omega_n\mxi|^m} \,.
\end{align*}
Now the claim is obvious since $\frac{\eps_n}{\omega_n}\leq \frac{2}{c_0}$,
$\bigl|1-\bigl(\frac{\eps_n}{\omega_n}\bigr)^m\bigr|
\leq \max\bigl\{1,(\frac{2}{c_0})^m-1\bigr\}$, and by 
Corollary \ref{cor:examples}(iii), Theorem \ref{thm:multipliers} 
and Lemma \ref{lm:pom_mihlin} we have that $\bigl(\frac{1}{1+|\omega_n\mxi|^m}\bigr)_n$
is a bounded sequence of $\mathrm{L}^p$ Fourier multipliers.
\end{proof}

For the stationary sequence $\omega_n=1$ we have that $(\omega_n)$-compactness is 
equivalent to the strong convergence to zero in $\mathrm{W}^{-m,p}(\Rd;\Cr)$, 
from the previous lemma we immediately get the following corollary.

\begin{corollary}\label{cor:rhs_W-1}
Let $m\in\N$ and $p\in (1,\infty)$, and  let $(\eps_n)$ be a sequence 
of positive real numbers.
\begin{itemize}
	\item[i)] If\/ $\limsup\limits_{n\to\infty} \eps_n =\eps_\infty<\infty$, then $(\eps_n)$-compactness
	implies the strong convergence to zero in $\W{-m}{p}{\Rd;\Cr}$.
	\item[ii)] If\/ $\liminf\limits_{n\to\infty} \eps_n =\eps_0>0$, then 
	the strong convergence to zero in $\W{-m}{p}{\Rd;\Cr}$
	implies $(\eps_n)$-compactness.
\end{itemize}	
\end{corollary}

We are particularly interested in the case when $\eps_n\to 0^+$. 
By Corollary \ref{cor:rhs_W-1}, in this case \eqref{ConvRhsLp} is a stronger condition than 
the strong convergence to zero in $\mathrm{W}^{-m,p}(\Rd;\Cr)$.

On the other hand, since $\bigl(\frac{1}{1+|\eps_n\mxi|^m}\bigr)_n$
is a bounded sequence of $\mathrm{L}^p$ Fourier multipliers, 
the strong convergence to zero in $\mathrm{L}^p$ always implies
\eqref{ConvRhsLp}.
Here we present a procedure how to modify a strongly 
convergent sequence in a 
Sobolev space to satisfy \eqref{ConvRhsLp} in the regime of a bounded 
characteristic length $(\eps_n)$.

\begin{lemma}\label{lm:rhs_convergence}
Let $(\eps_n)$ be a sequence of positive real numbers bounded from above,
let  $k\in\{0,1,\dots,m\}$ and let $(\vf_n)$ be a sequence
of vector valued functions in  $\W{-k}p{\Rd;\Cr}$ which converges strongly to zero.
Then $(\eps_n^k\vf_n)$ satisfies \eqref{ConvRhsLp}.
\end{lemma}
\begin{proof}
Let us first note that by an application of Corollary \ref{cor:examples}(iii,iv) 
and Lemma \ref{lm:pom_mihlin}, we can replace the symbol in \eqref{ConvRhsLp}
by $\mxi\mapsto (1+|\eps_n\mxi|^2)^{-m/2}$ .

Now we have
\begin{align*}
\frac{\eps_n^k}{(1+|\eps_n \mxi|^2)^{\frac{m}{2}}} 
	&= \Bigr(\frac{\eps_n^2+|\eps_n\mxi|^2}{1+|\eps_n\mxi|^2}\Bigr)^{\frac k2} 
	\frac{1}{(1+|\eps_n\mxi|^2)^{\frac{m-k}{2}}} 
	\frac{1}{(1+|\mxi|^2)^{\frac k2}} \\
&=: \psi_1(\eps_n\mxi)\psi_2(\eps_n\mxi)\psi_3(\mxi) \;.
\end{align*}
By the assumption we have that $\lA_{\psi_3}\vf_n$ converges strongly to zero 
in $\mathrm{L}^p$, while $\psi_2(\eps_n\,\cdot)$ is a bounded sequence of 
$\mathrm{L}^p$ Fourier multipliers (see Corollary \ref{cor:examples}(iii,iv), 
Theorem \ref{thm:multipliers} and Lemma \ref{lm:pom_mihlin}). 
Since 
$\lA_{\psi_1(\eps_n\cdot)\psi_2(\eps_n\cdot)\psi_3}
= \lA_{\psi_1(\eps_n\cdot)}\lA_{\psi_2(\eps_n\cdot)}\lA_{\psi_3}$,
it is only left to see that $\psi_1(\eps_n\,\cdot)$ is 
a bounded sequence of $\mathrm{L}^p$ Fourier multipliers as well.

After expanding $\psi_1$ into factors
\begin{align*}
\psi_1(\eps_n\mxi)&=\frac{(1+|\mxi|^2)^{\frac{k}{2}}}{1+|\mxi|^k}
	\frac{1+|\eps_n\mxi|^k}{(1+|\eps_n\mxi|^2)^{\frac{k}{2}}}
	\frac{\eps_n^k+|\eps_n\mxi|^k}{1+|\eps_n\mxi|^k} \\
&= \frac{(1+|\mxi|^2)^{\frac{k}{2}}}{1+|\mxi|^k}
\frac{1+|\eps_n\mxi|^k}{(1+|\eps_n\mxi|^2)^{\frac{k}{2}}}
\biggl(1+\frac{\eps_n^k-1}{1+|\eps_n\mxi|^k}\biggr) \,,
\end{align*}
the claim follows by Corollary \ref{cor:examples}(iii, iv), 
Theorem \ref{thm:multipliers}, Lemma \ref{lm:pom_mihlin},
and the boundedness of $(\eps_n)$.
\end{proof}

Let us close this subsection by two examples of weakly (but not strongly) converging 
sequences in an $\mathrm{L}^p$ space,
satisfying the $(\eps_n)$-(local) compactness condition.

\begin{example}\label{exa:rhs-conv}
Let $(\eps_{n})$ be a sequence of positive real numbers,
$p\in (1,\infty)$, $m\in\N$ and $r=1$.
\begin{itemize}
\item[i)] For $\mz\in\Rd$ let $\zeta_{p,\eps}$ be defined as 
in Subsection \ref{subsec:conc} (recall that it depends on $\mz$, even though it is not explicitly included in the notation).

For an arbitrary $u\in\mathrm{L}^p(\Rd)$, $u\neq 0$ and $\omega_n\to 0^+$,
let us see under which conditions on $(\omega_n)$ the sequence 
$(\zeta_{p,\omega_n}u)$ satisfies the $(\eps_n)$-local compactness condition.

By Lemma \ref{lm:concentration_conv} our task is equivalent to examining the 
$(\eps_n)$-compactness condition, while by 
Lemma \ref{lm:concentration_comm} we have
\begin{equation*}
\Bigl\|\lA_{\frac{1}{1+|\eps_n\mxi|^m}}\zeta_{p,\omega_n}u\Bigr\|_{\mathrm{L}^p} =
	\Bigl\|\lA_{\frac{1}{1+|\frac{\eps_n}{\omega_n}\mxi|^m}}u
	\Bigr\|_{\mathrm{L}^p} \;.
\end{equation*}

If $c:=\lim_n\frac{\eps_n}{\omega_n}$ exists in $[0,\infty]$, then 
for a fixed $\mxi\in\Rdz$ we have
\begin{equation*}
\lim_n \frac{1}{1+|\frac{\eps_n}{\omega_n}\mxi|^m} =
	\left\{\begin{array}{ll}
	\frac{1}{1+|c\mxi|^m} \;, & c\in [0,\infty) \\
	0\;, & c=\infty
	\end{array}\right. \;.
\end{equation*}
Therefore, by Corollary \ref{cor:examples}(iii), 
Theorem \ref{thm:multipliers} and Lemma \ref{lm:pom_mihlin},
the assumptions of Lemma \ref{lm:technical_symbol} are fulfilled, 
implying 
\begin{equation*}
\lim_n \Bigl\|\lA_{\frac{1}{1+|\frac{\eps_n}{\omega_n}\mxi|^m}}u
	\Bigr\|_{\mathrm{L}^p} =
	\left\{\begin{array}{ll}
		\|\lA_{\frac{1}{1+|c\mxi|^m}}u\|_{\mathrm{L}^p} \;, & c\in [0,\infty) \\
	0\;, & c=\infty
	\end{array}\right. \;.
\end{equation*}
Hence, the sequence $(\zeta_{p,\omega_n}u)$ satisfies $(\eps_n)$-(local)
compactness condition if and only if $\lim_n\frac{\eps_n}{\omega_n}=\infty$. 

In a special case when $0<\liminf_n\eps_n\leq\limsup_n\eps_n<\infty$,
this result implies that $(\zeta_{p,\omega_n}u)$
converges strongly to zero in $\W{-m}{p}{\Rd}$ 
(see Corollary \ref{cor:rhs_W-1}), which is a well known fact (under the assumption of uniformly bounded support,
or if one uses the local spaces), since 
the sequence converges weakly to zero in $\mathrm{L}^p(\Rd)$.

\item[ii)] Let $\mathsf{k}\in\Rdz$, $\omega_n\to 0^+$ and 
$u\in\mathrm{L}^p_\mathrm{loc}(\Rd)$, $u\neq 0$, and consider
$$
u_n(\mx) = u(\mx) e^{2\pi i \frac{\mx}{\omega_n}\cdot\mathsf{k}} \;.
$$
By \eqref{eq:prFtr} we have
\begin{equation*}
\Bigl\|\lA_{\frac{1}{1+|\eps_n\mxi|^m}}(\varphi u_n)\Bigr\|_{\mathrm{L}^p} =
	\Bigl\|\lA_{\frac{1}{1+|\eps_n\mxi+\frac{\eps_n}{\omega_n}\mathsf{k}|^m}}
	(\varphi u) \Bigr\|_{\mathrm{L}^p} \;,
\end{equation*}
where $\varphi\in\mathrm{C}^\infty_c(\Rd)$ is arbitrary.
As in the first example, our aim is to apply Lemma \ref{lm:technical_symbol}. 

Firstly, let us notice that for any $r\in (1,\infty)$ the sequence of functions
$\frac{1}{1+|\eps_n\mxi+\frac{\eps_n}{\omega_n}\mathsf{k}|^m}$
defines a bounded sequence of\/ $\mathrm{L}^r$ Fourier multipliers.
Indeed, the claim follows by Corollary \ref{cor:examples}(iii), 
Theorem \ref{thm:multipliers}, Lemma \ref{lm:pom_mihlin}, and the 
fact that Fourier multipliers are invariant under translations
(see \cite[Proposition 2.5.14]{Grafakos}). 

It remains to study the pointwise limit of the sequence of symbols. 
To this end, let us assume that $c:=\lim_n \frac{\eps_n}{\omega_n}$
exists in $[0,\infty]$. Thus, for a fixed $\mxi\in\Rd$ we have
\begin{equation*}
\lim_n \frac{1}{1+|\eps_n\mxi+\frac{\eps_n}{\omega_n}\mathsf{k}|^m} =
\left\{\begin{array}{ll}
\frac{1}{1+|c\mathsf{k}|^m} \;, & c\in [0,\infty) \\
0\;, & c=\infty
\end{array}\right. \;.
\end{equation*}
Therefore, by Lemma \ref{lm:technical_symbol} we get
\begin{equation*}
\lim_n \Bigl\|\lA_{\frac{1}{1+|\eps_n\mxi+\frac{\eps_n}{\omega_n}\mathsf{k}|^m}}
	(\varphi u) \Bigr\|_{\mathrm{L}^p} =
\left\{\begin{array}{ll}
\frac{\|\ph u\|_{\mathrm{L}^p}}{1+|c\mathsf{k}|^m} \;, & c\in [0,\infty) \\
0\;, & c=\infty
\end{array}\right. \;,
\end{equation*}
implying that the sequence $(u_n)$ satisfies $(\eps_n)$-local compactness
condition if and only if $\lim_n\frac{\eps_n}{\omega_n}=\infty$.
\end{itemize}
\end{example}

\begin{remark}
In both previous examples we can see that the crucial relation between 
the characteristic length of the sequence $(\omega_n)$ and 
the characteristic length of the condition $(\eps_n)$ is 
$\lim_n\frac{\eps_n}{\omega_n} = \infty$. 
Thus, in terms of \cite{EL18}, the characteristic lengths of sequences 
should be \emph{faster}
than the characteristic lengths of the condition. 
One can notice that this is a stronger requirement than what is needed 
for the $(\eps_n)$-oscillatory property (see \cite[examples 1 and 2]{EL18}).
\end{remark}

\subsection{Localisation principle}

After completing all these preliminaries, we are ready to present 
the localisation principle for one-scale H-distributions. 

\begin{theorem}\label{thm:loc_princ_jhd}
Take $\Omega\subseteq\Rd$ open and $p\in(1,\infty)$. Let $\vu_n\rightharpoonup 0$ in $\mathrm{L}_\mathrm{loc}^{p}(\Omega;\Cr)$
satisfy \eqref{eq:loc_princ_eq},
%\begin{equation}\label{locprEqJhdA}
%\sum_{|\malpha|\leq m} \eps_n^{|\malpha|} \partial_\malpha(\mA^\malpha\vu_n) = \vf_n \,,
%\end{equation}
where $(\eps_n)$ is a sequence of positive real numbers, 
$\mA^\malpha_n\in\mathrm{C}(\Omega;\Mx{q\times r}{\C})$,
such that for any $\malpha\in\N_0^d$ the sequence $\mA^\malpha_n\to \mA^\malpha$ 
in the space $\mathrm{C}(\Omega;\Mx{q\times r}{\C})$ 
(in other words, $\mA^\malpha_n$ converges locally uniformly to $\mA^\malpha$),
while $(\vf_n)$
is a sequence of functions in $\Wloc{-m}p{\Omega;\Cr}$ satisfying 
$(\eps_n)$-local compactness condition \eqref{ConvRhsLocLp}.
Moreover, let $(\vv_n)$ be a bounded sequence in\/ $\mathrm{L}_\mathrm{loc}^{p'}(\Omega;\Cr)$ and let 
$\omega_n\to 0^+$ be a sequence of positive reals such that 
$c:=\lim_n\frac{\omega_n}{\eps_n}$ exists (in $[0,\infty]$).

Then any one-scale H-distribution $\jhdv^{(\omega_n)}$ associated 
to (sub)sequences (of) $(\vu_n)$ and $(\vv_n)$
with characteristic length $(\omega_n)$ satisfies:
\begin{equation}\label{eq:loc_princ_jhd}
\mmp_c(\mx,\mxi)\jhdv = \mnul \,,
\end{equation}
where, with respect to the value of $c$, we have
\begin{itemize}
\item[i)] $c=0:$
$$
\mmp_0(\mx,\mxi) = \sum_{|\malpha|= m} (2\pi i)^m 
{\mxi^\malpha\over 1+|\mxi|^{m}}\mA^\malpha(\mx) \,,
$$
\item[ii)] $c\in(0,\infty):$
$$
\mmp_c(\mx,\mxi) = \sum_{|\malpha|\leq m} \Bigl(\frac{2\pi i}{c}\Bigr)^{|\malpha|} 
{\mxi^\malpha\over 1+|\mxi|^{m}}\mA^\malpha(\mx) \,,
$$
\item[iii)] $c=\infty:$
$$
\mmp_\infty(\mx,\mxi) = {1\over 1+|\mxi|^{m}}\mA^\vnul(\mx) \,.
$$
\end{itemize}
\end{theorem}

\begin{proof}
\noindent\textbf{I.}
Let us first show that we can 
replace $\mA^\malpha_n$ by $\mA^\malpha$ in \eqref{eq:loc_princ_eq}. 
It is sufficient to show that 
$$
\sum_{|\malpha|\leq m} \eps_n^{|\malpha|} \partial_\malpha 
	\Bigl(\bigl(\mA^\malpha_n-\mA^\malpha\bigr) \vu_n\Bigr) 
$$
satisfies the $(\eps_n)$-local compactness condition \eqref{ConvRhsLocLp}, as then the new right hand side
\begin{equation*}
\tilde{\vf}_n := \vf_n - 
	\sum_{|\malpha|\leq m} \eps_n^{|\malpha|} \partial_\malpha 
	\Bigl(\bigl(\mA^\malpha_n-\mA^\malpha\bigr) \vu_n\Bigr)
\end{equation*}
still satisfies the same $(\eps_n)$-local compactness condition. 

Since $(\mA^\malpha_n-\mA^\malpha) \vu_n$ strongly converges to 
$\vnul$ in $\mathrm{L}^p_\mathrm{loc}(\Omega;\Cr)$, 
the sequence of functions $\partial_\malpha((\mA^\malpha_n-\mA^\malpha)\vu_n)$
strongly converges to $\vnul$ in $\Wloc{-|\malpha|}{p}{\Omega;\Cr}$.
Thus, the claim follows by Lemma \ref{lm:rhs_convergence}.

For simplicity of notation, we shall continue to use $\vf_n$ instead of $\tilde{\vf}_n$
for denoting the right hand side. 
\smallskip

\noindent\textbf{II.}
In the next step we localise \eqref{eq:loc_princ_eq} 
by multiplying it by a test function 
$\ph\in\mathrm{C}_c^\infty(\Omega)$, and after applying the Leibniz rule (cf.~\cite[Lemma 8]{AEL} and \cite[Lemma 9]{MEphd} for details) we have
$$
\sum_{|\malpha|\leq m}\sum_{\vnul\leq\mbeta\leq\malpha}(-1)^{|\mbeta|}
	{\malpha\choose \mbeta}\eps_n^{|\malpha|}
	\partial_{\malpha-\mbeta}\Bigl((\partial_\mbeta\ph)\mA^\malpha\vu_n \Bigr) 
	= \ph\vf_n \,.
$$
Our aim is to show that the terms with $\mbeta\ne\vnul$ on the left hand side of
this equality satisfy the $(\eps_n)$-compactness condition \eqref{ConvRhsLp}, 
i.e.~that we can insert them into the right hand side.
	
For given $\malpha$ and $\mbeta$ the sequence $\partial_{\malpha-\mbeta}
\Bigl((\partial_\mbeta\ph)\mA^\malpha\vu_n \Bigr)$ is supported in a 
fixed compact ($\supp \ph$), hence by the Rellich compactness theorem 
for $\mbeta\ne\vnul$ it is strongly precompact in $\W{-|\malpha|}p{\Rd;\Cr}$.
Therefore, by Lemma \ref{lm:rhs_convergence} we get that $\eps_n^{|\malpha|}\partial_{\malpha
	-\mbeta}\Bigl((\partial_\mbeta\ph)\mA^\malpha\vu_n \Bigr)$ satisfies 
\eqref{ConvRhsLp}. 
	
By the arbitrariness of $\malpha$ and $\mbeta$ we get the claim, so the equality 
above can be rewritten in the form 
\begin{equation}\label{eq:loc_princ_pom_eq}
\sum_{|\malpha|\leq m}\eps_{n}^{|\malpha|} \partial_\malpha
	(\mA^\malpha\ph\vu_n) = \tilde\vf_n \,,
\end{equation}
where the sequence $(\tilde\vf_n)$ satisfies \eqref{ConvRhsLp}.
\smallskip

\noindent\textbf{III.}
Let us first assume that $\omega_n=\eps_n$, i.e.~$c=1$, while 
the proof of the general case we postpone to the last step. 
Of course, here we implicitly assume that $\eps_n\to 0^+$. 

By applying $\lA_{\frac{1}{1+|\eps_n\mxi|^m}}$ to \eqref{eq:loc_princ_pom_eq}
we get
$$
\sum_{|\malpha|\leq m} \lA_{(2\pi i)^{|\malpha|}\psi^{m,\malpha}_n}
	(\ph\mA^\malpha\vu_n) \longrightarrow\vnul \quad \hbox{in} \quad  
	\mathrm{L}^p(\Rd;\Cr) \,,
$$
where $\psi^{m,\malpha}(\mxi)=\frac{\mxi^\malpha}{1+|\mxi|^m}$
(this function is in $\pC{\infty}{\kob}$; see 
Corollary \ref{cor:examples}(iii)) and we denote
$\psi^{m,\malpha}_n = \psi^{m,\malpha}(\eps_n\,\cdot)$. 
Further on, for an arbitrary $\psi\in\pC{d(\lfloor\frac{d}{2}\rfloor+3)}{\kob}$,
we apply $\lA_{\psi_n}$ (here $\psi_n=\psi(\eps_n\,\cdot)$) to the above
sum and form a complex tensor product with $\ph_1\vv_n$, 
for an arbitrary $\ph_1\in\mathrm{C}_c^\infty(\Omega)$.
After using $u_n^s$, $v_n^l$ and $A^\malpha_{ks}$ to denote components of
$\vu_n$, $\vv_n$ and $\mA^\malpha$ respectively,  
for the $(k,l)$ component of the above 
sum on the limit we get
$$
\begin{aligned}
0 &= \sum_{|\malpha|\leq m}\sum_{s=1}^d \lim_n \int_\Rd \lA_{(2\pi i)^{|\malpha|}
	\psi_n\psi_n^{m,\malpha}}(\ph A^{\malpha}_{ks}u_n^s) \overline{\ph_1 v_n^l} \,d\mx \\
&= \sum_{|\malpha|\leq m}\sum_{s=1}^d \Dupp{\jhd^{sl}}{\ph\bar{\ph_1} 
	A^\malpha_{ks}\otimes(2\pi i)^{|\malpha|}\psi\psi^{m,\malpha}} \\
&= \sum_{|\malpha|\leq m}\sum_{s=1}^d \Dupp{(2\pi i)^{|\malpha|}\psi^{m,\malpha}
	A^\malpha_{ks}\jhd^{sl}}{\ph\bar\ph_1 \otimes\psi} \\
&= \Dupp{\sum_{|\malpha|\leq m} \frac{(2\pi i\mxi)^\malpha}{1+|\mxi|^m}
	[\mA^{\malpha}\jhdv]_{kl}}{\ph\bar\ph_1\otimes\psi} \,,
\end{aligned}
$$
where $\jhdv$ is the one-scale H-distribution with characteristic scale 
$(\eps_n)$ associated to the chosen subsequences of $(\vu_n)$ and $(\vv_n)$.

Since test functions $\ph, \ph_1$ and $\psi$ are arbitrary, and tensor products
are dense in the space 
$\mathrm{C}_c^{0,d(\lfloor\frac{d}{2}\rfloor+3)}(\Omega\times\kob)$, 
we get the claim.
\smallskip

\noindent\textbf{IV.}
Let us consider now the case of a general sequence $\omega_n\to 0^+$. 
The aim is to transform the equation to the form in which the characteristic
length $\eps_n$ is replaced by $\omega_n$ and in which the right hand side 
satisfies the $(\omega_n)$-local compactness condition. 
Then we can apply steps I--III in order to get the result.

If $c\in (0,\infty]$, we rewrite \eqref{eq:loc_princ_eq} in the form
$$
\sum_{|\malpha|\leq m} \omega_n^{|\malpha|} \partial_\malpha 
	\bigl(\mB^\malpha_n \vu_n\bigr) = \vf_n \,,
$$
where $\mB^\malpha_n:= \bigl(\frac{\eps_{n}}{\omega_{n}}\bigr)^{|\malpha|}
\mA^\malpha_n$. Since $\mB^\alpha_n$ converges to 
$c^{-|\malpha|}\mA^\malpha$ in $\mathrm{C}(\Omega;\Mx{q\times r}{\C})$
and by Lemma \ref{lm:rhs-omegan}(i) the sequence $(\vf_n)$ satisfies the
$(\omega_n)$-local compactness condition, we are in a position 
to apply steps I--III. 
Let us note that for $c=\infty$ non-trivial terms are only present for $\malpha=\vnul$.

For $c=0$ we first multiply \eqref{eq:loc_princ_eq} by 
$(\frac{\omega_n}{\eps_n})^m$ and then rewrite it in the form
$$
\sum_{|\malpha|\leq m} \omega_n^{|\malpha|} \partial_\malpha 
	\bigl(\mathbf{C}^\malpha_n \vu_n\bigr) = \vg_n \,,
$$
where $\mathbf{C}^\malpha_n:= \bigl(\frac{\omega_{n}}{\eps_{n}}\bigr)^{m-|\malpha|}
\mA^\malpha_n$ and $\vg_n:= \bigl(\frac{\omega_{n}}{\eps_{n}}\bigr)^{m}\vf_n$.
It is obvious that $\mathbf{C}^\malpha_n$ converges 
in $\mathrm{C}(\Omega;\Mx{q\times r}{\C})$
and the limit is $\mnul$ for $|\malpha|\leq m-1$ and $\mA^\malpha$ for
$|\malpha|=m$.
It remains to be shown that $(\vg_n)$ satisfies the $(\omega_n)$-local compactness
condition. Indeed, we have
\begin{align*}
\Bigl(\frac{\omega_{n}}{\eps_n}\Bigr)^m\frac{1}{1+|\omega_n\mxi|^m}
	&= \frac{(\frac{\omega_{n}}{\eps_n})^m + 
		|\omega_{n}\mxi|^m}{1+|\omega_{n}\mxi|^m}
	\frac{1}{1+|\eps_n\mxi|^m} \\
&= \biggl(1+\frac{(\frac{\omega_{n}}{\eps_n})^m-1}{1+|\omega_m\mxi|^m}\biggr)
	\frac{1}{1+|\eps_n\mxi|^m} \;.
\end{align*}
Now we take into account the fact that $(\vf_n)$ satisfies the $(\eps_n)$-local compactness 
condition (for 
which we use the second factor), while the first factor defines 
a bounded sequence of $\mathrm{L}^p$ Fourier multipliers, which can be shown 
following the same argument as it was done in the proof of Lemma \ref{lm:rhs-omegan}.
\end{proof}

\begin{remark}
\begin{itemize}
\item[i)] The application of Theorem \ref{thm:loc_princ_jhd}, together with 
Lemma \ref{lm:strong_compactness}, could be seen as an $\mathrm{L}^p$ 
variant of \emph{compactness by compensation}\/ (in literature also referred to as 
\emph{compensated compactness}\/) suitable for problems with a characteristic 
length. In this context it is important to notice that we allow for variable 
coefficients of the corresponding partial differential equations, 
which need to be merely continuous. 
\item[ii)] In this paper we chose to work with bilinear dual products,
while in \cite{AEL} sesquilinear dual products were used. 
Because of that there is a small difference in formulae
\eqref{eq:loc_princ_jhd} and \cite[(4.10)]{AEL} (when compared for 
$p=2$): here we do not have the transpose of the one-scale H-distribution.

\item[iii)] The case $\eps_n\to 0^+$ is of the utmost importance in the 
problems with characteristic scales,
but let us briefly address some remarks concerning the case 
$0<\eps_0\leq\eps_n\leq\eps_\infty<\infty$, which is the case 
$c=0$ of Theorem \ref{thm:loc_princ_jhd}. Since $\mmp_0$ is equal to zero 
on $\Sigma_0$, our result \eqref{eq:loc_princ_jhd} does not provide 
information about the distribution on the whole domain. 

This problem can be solved by showing that the localisation principle 
\eqref{eq:loc_princ_jhd} still holds if we replace symbol $\mmp_0$ by
$$
\tilde{\mmp}_0(\mx,\mxi) = \sum_{|\malpha|= m} (2\pi i)^m 
{\mxi^\malpha\over |\mxi|^{m}}\mA^\malpha(\mx) \;,
$$
i.e.~one needs to remove '$1$' in the denominator.
In \cite{AEL} that was done in the $\mathrm{L}^2$ setting 
by applying a technical lemma (see Lemma 9 in the aforementioned reference).
In the $\mathrm{L}^p$ framework this is a bit more demanding since one 
needs to work with the Fourier multiplier. 
The standard approach is to use Riesz potentials in this case, 
as it was done in \cite[Theorem 4.1]{AM} for the first order differential equations. 
However, this approach requires to have a better integrability property
of the sequence $(\vv_n)$.

Of course, this adjustment, i.e.~the removal of '$1$', cannot be done (in general)
when $c=0$, but $\liminf_n\eps_n=0$. 

\item[iv)] The previous theorem shows that $(\eps_n)$-local compactness
of the right hand side is sufficient to ensure \eqref{eq:loc_princ_jhd}.
For the case $p=2$ it was proven in \cite[Theorem 9]{AEL} that it is also
a necessary condition.
\end{itemize}
\end{remark}

\subsection{Applications}

In this subsection we present a possible application of Theorem 
\ref{thm:loc_princ_jhd} in a simple example. 

Let $\Omega\subseteq\R^2$ be open and let $p\in (1,\infty)$. 
Assume that $u_n\rightharpoonup 0$ in $\mathrm{L}_\mathrm{loc}^p(\Omega)$
and $v_n\rightharpoonup 0$ in $\mathrm{L}_\mathrm{loc}^{p'}(\Omega)$
satisfy
\begin{equation}\label{eq:loc-princ_exa_eq}
\left\{
\begin{aligned}
u_n + \eps_n \partial_{x_1}(a_1 u_n) &= f_n \\
v_n + \eps_n^\alpha \partial_{x_2}(a_2 v_n) &= g_n \,,
\end{aligned}\right.
\end{equation}
where $\eps_n\to 0^+$, $\alpha\in(0,1]$, $f_n\in\Wloc{-1}{p}\Omega$ 
and $g_n\in\Wloc{-1}{p'}\Omega$ satisfy the
$(\eps_n)$-local compactness condition for $m=1$ ($f_n$ in $\mathrm{L}^p$ and 
$g_n$ in $\mathrm{L}^{p'}$), and $a_1, a_2\in \Cp{\Omega;\R}$ such 
that $a_1, a_2\neq 0$ everywhere. 

Let $\mu$ be a one-scale H-distribution associated to (sub)sequences (of)\/
$(u_n)$ and $(v_n)$ with characteristic length $(\eps_n)$. 
For simplicity, here we do not relabel the subsequences, but 
one needs to be aware that all conclusions extracted from one-scale 
H-distributions are at the level of the chosen subsequences (unless the sequences are 
$(\eps_n)$-pure). This is a common situation when studying weakly converging sequences. 

After applying Theorem \ref{thm:loc_princ_jhd} to the first equality in 
\eqref{eq:loc-princ_exa_eq} we get
\begin{equation*}
\Bigl(\frac{1}{1+|\mxi|} + \frac{2\pi i a_1(\mx)\xi_1}{1+|\mxi|}\Bigr)\mu = 0 \;.
\end{equation*}
The symbol
$$
p(\mx,\mxi):= \frac{1}{1+|\mxi|} + \frac{2\pi i a_1(\mx)\xi_1}{1+|\mxi|}
$$
is equal to zero only when both functions $\mxi\mapsto \frac{1}{1+|\mxi|}$
and $\mxi\mapsto\frac{a_1(\mx)\xi_1}{1+|\mxi|}$ are equal to zero,
which occurs only for $(\mx,\mxi)\in\Omega\times \{\infty^{(0,-1)},\infty^{(0,1)}\}$.
Since for any $\Phi\in\pC{0,d(\lfloor\frac{d}{2}\rfloor+3)}{\Omega\times\kob}$
equal to zero on a neighbourhood of both points 
$\infty^{(0,-1)}$ and $\infty^{(0,1)}$ in $\kob$ we have 
$\frac{1}{p}\Phi \in \pC{0,d(\lfloor\frac{d}{2}\rfloor+3)}{\Omega\times\kob}$,
the above localisation principle implies 
\begin{equation*}
\Dupp{\mu}{\Phi} = \Dupp{\mu}{p(\Phi/p)} = \Dupp{p\mu}{\Phi/p} = 0 \;.
\end{equation*}
Therefore, the support of $\mu$ is localised within the set 
$\Omega\times \{\infty^{(0,-1)},\infty^{(0,1)}\}$.

Let us now use the second equality in \eqref{eq:loc-princ_exa_eq}.
A direct application of Theorem \ref{thm:loc_princ_jhd} 
requires from us to work with the one-scale H-distribution associated to (sub)sequences (of)
$(v_n)$ and $(u_n)$ (the opposite order of sequences than it was above) 
with characteristic length $(\eps_n)$.
However, by Lemma \ref{lm:jhd-bar} this one-scale H-distribution 
is equal to $\bar{\mu}$.
Furthermore, depending on the value of $\alpha$, we have two possible cases. 
Let us consider only the case $\alpha\in(0,1)$ (the analysis for the 
remaining case $\alpha=1$ is completely analogous).
Thus, we have $\frac{2\pi i a_2(\mx)\xi_2}{1+|\mxi|}\bar\mu = 0$,
implying
\begin{equation*}
\frac{a_2(\mx)\xi_2}{1+|\mxi|}\mu = 0 \;.
\end{equation*}
Since the function $(\mx,\mxi)\mapsto \frac{a_2(\mx)\xi_2}{1+|\mxi|} \in 
\pC{0,\infty}{\Omega\times\kob}$ has no zeros in a neighbourhood
of the set $\Omega\times \{\infty^{(0,-1)},\infty^{(0,1)}\}$, 
as in the above argument we can conclude that 
$\supp\mu\cap \bigl(\Omega\times \{\infty^{(0,-1)},\infty^{(0,1)}\}\bigr)
=\emptyset$. 
Therefore, the overall conclusion is $\mu=0$. 

By the definition of the one-scale H-distribution (see Remark \ref{rem:jhd-star-conv})
this means that there exist subsequences $(u_{n'})$ and $(v_{n'})$ such that 
$u_{n'}\bar{v}_{n'}\xrightharpoonup{\;*\;} 0$ in the sense of (unbounded) Radon measures.

Let us close this subsection by presenting a concrete choice for 
data in \eqref{eq:loc-princ_exa_eq}.

\begin{example}
Let $\Omega = \R^2$ and let $p\in (1,\infty)$. 
For $u_1, u_2\in \mathrm{L}^p(\R)$ and  
$v_1, v_2\in\mathrm{L}^{p'}(\R)$ we define
\begin{align*}
u_n(x_1,x_2) &:= u_1(x_1)e^{2\pi i n x_1} + n^{2/p}u_2(n^2x_2) \,,\\
v_n(x_1,x_2) &:= n^{2/p'}v_1(n^2x_1) + v_2(x_2)e^{2\pi i n x_2} \,.
\end{align*}

The sequence of functions
\begin{align*}
f_n(x_1,x_2)&:=u_n(x_1,x_2) - \frac{1}{n}\partial_{x_1}u_n(x_1,x_2) \\
&= u_1(x_1)e^{2\pi i n x_1} + n^{2/p}u_2(n^2x_2) -
	\frac{1}{n}u_1'(x_1)e^{2\pi i nx_1} - u_1(x_1)e^{2\pi i nx_1} \\
&= n^{2/p}u_2(n^2x_2) - \frac{1}{n}u_1'(x_1)e^{2\pi i nx_1} 
\end{align*}
satisfies the $(1/n)$-local compactness condition in $\mathrm{L}^p(\R^2)$ 
(with $m=1$).
Indeed, for the first term in the last equality we apply
Example \ref{exa:rhs-conv}(i). On the other hand, for the second term we proceed by
%definition and reduce it (using \eqref{eq:prFtr}) to the application of Lemma \ref{lm:technical_symbol}.
definition, transform it by using \eqref{eq:prFtr}, and then apply Lemma \ref{lm:technical_symbol}.
However, note that in the case $u_1\in\W{1}{p}{\R}$ we can apply Lemma \ref{lm:rhs_convergence}
directly since then $(\frac{1}{n} u_1'(x_1)e^{2\pi i nx_1})$
converges strongly to zero in $\mathrm{L}^p(\R)$, hence also in $\mathrm{L}^p_{\mathrm{loc}}(\R^2)$ (since it is independent of $x_2$).

Analogously,  
\begin{align*}
g_n(x_1,x_2)&:= v_n(x_1,x_2) - \frac{1}{n}\partial_{x_2}v_n(x_1,x_2) \\
&= n^{2/p'}v_1(n^2x_1) - \frac{1}{n}v_2'(x_2)e^{2\pi i nx_2} 
\end{align*}
satisfies the $(1/n)$-local compactness condition in $\mathrm{L}^{p'}(\R^2)$  
(again with $m=1$).

Thus, sequences $(u_n)$ and $(v_n)$ satisfy \eqref{eq:loc-princ_exa_eq}
with $\eps_n=\frac{1}{n}$, $\alpha=1$ and $a_1=a_2=-1$, so that
$u_{n}\bar{v}_{n}\xrightharpoonup{\;*\;} 0$.
This can also be checked by a direct inspection of sequences. 

In this example the whole sequence converges since sequences 
$(u_n)$ and $(v_n)$ are $(1/n)$-pure (see subsections 
\ref{subsec:conc} and \ref{subsec:oscill}).
\end{example}

%========================================================================
%===============================================================================
\section*{Acknowledgements}%\label{sec:thanks}
% =========================================================================
The authors wish to thank Eduard Nigsch for numerous helpful
discussions and the referees
for their insightful comments that helped to improve the presentation and the quality
of the paper.

\section*{Statements and Declarations}

This work is supported in part by the Croatian Science Foundation under the project 
%IP-2013-11-9780 (WeConMApp) and 
IP--2018--01--2449 (MiTPDE) and by the Croatian-Austrian bilateral project
{\sl Anisotropic distributions and H-distributions}.

The authors have no relevant financial or non-financial interests to disclose.

Data sharing not applicable to this article as no datasets were generated or analysed during the current study.

\end{document}